\input amstex
\magnification=1200
\loadmsbm
\loadeufm
\loadeusm
\UseAMSsymbols
\input amssym.def

\font\BIGtitle=cmr10 scaled\magstep3
\font\bigtitle=cmr10 scaled\magstep1
\font\boldsectionfont=cmb10 scaled\magstep1
\font\section=cmsy10 scaled\magstep1

\def\scr#1{{\fam\eusmfam\relax#1}}
\def\scrA{{\scr A}}
\def\scrB{{\scr B}}
\def\scrC{{\scr C}}
\def\scrD{{\scr D}}
\def\scrE{{\scr E}}
\def\scrF{{\scr F}}
\def\scrG{{\scr G}}
\def\scrH{{\scr H}}
\def\scrI{{\scr I}}
\def\scrK{{\scr K}}
\def\scrJ{{\scr J}}
\def\scrL{{\scr L}}
\def\scrM{{\scr M}}
\def\scrN{{\scr N}}
\def\scrO{{\scr O}}
\def\scrP{{\scr P}}

\def\scrQ{{\scr Q}}
\def\scrS{{\scr S}}

\def\scrV{{\scr V}}
\def\scrZ{{\scr Z}}
\def\scrW{{\scr W}}
\def\scrR{{\scr R}}
\def\scrT{{\scr T}}
\def\scrX{{\scr X}}
\def\scrY{{\scr Y}}
\def\gr#1{{\fam\eufmfam\relax#1}}

\def\grA{{\gr A}}	
	
\def\grC{{\gr C}}	\def\grc{{\gr c}}

\def\grG{{\gr G}}	
	
\def\grI{{\gr I}}	
\def\grJ{{\gr J}}	
\def\grK{{\gr K}}	
	
\def\grl{{\gr l}}

\def\grP{{\gr P}}	 
	
	\def\grr{{\gr r}}
	\def\grs{{\gr s}}
\def\grT{{\gr T}}	\def\grt{{\gr t}}
\def\grU{{\gr U}}	\def\gru{{\gr u}}

\def\grX{{\gr X}}	
	
	\def\grz{{\gr z}}

\def\db#1{{\fam\msbfam\relax#1}}

\def\dbA{{\db A}} 
\def\dbC{{\db C}} 
 \def\dbF{{\db F}}
\def\dbG{{\db G}} \def\dbH{{\db H}}

 \def\dbN{{\db N}}
 
\def\dbQ{{\db Q}} \def\dbR{{\db R}}

 \def\dbZ{{\db Z}}

\def\eps{{\varepsilon}}

\def\Ker{\text{Ker}}
\def\der{\text{der}}
\def\Sh{\hbox{\rm Sh}}

\def\sc{\text{sc}}
\def\Res{\text{Res}}

\def\ad{\text{ad}}

\def\c{\text{c}}
\def\an{\text{an}}

\def\Gal{\text{Gal}}
\def\GL{\text{GL}}

\def\Hom{\text{Hom}}
\def\End{\text{End}}

\def\Spec{\text{Spec}\,}
\def\Spf{\text{Spf}\,}

\def\Lie{\text{Lie}}

\def\leaderfill{\leaders\hbox to 1em
     {\hss.\hss}\hfill}
\def\nspace{\lineskip=1pt\baselineskip=12pt\lineskiplimit=0pt}

\def\endproof{$\hfill \square$}
\def\finishproclaim{\par\rm
     \ifdim\lastskip<\medskipamount\removelastskip
     \penalty55\medskip\fi}
\def\proof{\par\noindent {\it Proof:}\enspace}
\def\references#1{\par
  \centerline{\boldsectionfont References}\smallskip
     \parindent=#1pt\nspace}
\def\Ref[#1]{\par\hang\indent\llap{\hbox to\parindent
     {[#1]\hfil\enspace}}\ignorespaces}
\def\Item#1{\par\smallskip\hang\indent\llap{\hbox to\parindent
     {#1\hfill$\,\,$}}\ignorespaces}
\def\ItemItem#1{\par\indent\hangindent2\parindent
     \hbox to \parindent{#1\hfill\enspace}\ignorespaces}

\def\arrowsim{\,\smash{\mathop{\to}\limits^{\lower1.5pt
  \hbox{$\scriptstyle\sim$}}}\,}

\def\doublemaprights#1#2#3#4{\raise3pt\hbox{$\mathop{\,\,\hbox to
     #1pt{\rightarrowfill}\kern-30pt\lower3.95pt\hbox to
     #2pt{\rightarrowfill}\,\,}\limits_{#3}^{#4}$}}

\def\rightcapdownarrow{\raise9pt\hbox{$\ssize\cap$}\kern-7.75pt
     \Big\downarrow}

\def\rcapmapdown#1{\rightcapdownarrow\kern-1.0pt\vcenter{
     \hbox{$\scriptstyle#1$}}}

\def\rmapdown#1{\Big\downarrow\kern-1.0pt\vcenter{
     \hbox{$\scriptstyle#1$}}}
\def\rightsubsetarrow#1{{\ssize\subset}\kern-4.5pt\lower2.85pt
     \hbox to #1pt{\rightarrowfill}}
\def\longtwoheadedrightarrow#1{\raise2.2pt\hbox to #1pt{\hrulefill}
     \!\!\!\twoheadrightarrow}

\def\Gal{\operatorname{\hbox{Gal}}}
\def\Hom{\operatorname{\hbox{Hom}}}

\def\Im{\hbox{Im}}

\NoBlackBoxes
\parindent=25pt
\document
\footline={\hfil}

\null
\vskip 0.2in
\centerline{\BIGtitle On the Tate and Langlands--Rapoport conjectures}

\centerline{\BIGtitle  for special fibres of integral canonical models}

\centerline{\BIGtitle of Shimura varieties of abelian type}

\vskip 0.2in 
\centerline{\bigtitle Adrian Vasiu}

\vskip 0.1in
\centerline{October 17, 2012}
\footline={\hfill}
\null

\smallskip\noindent
{\bf ABSTRACT}.  We prove the isogeny property for special fibres of integral canonical models of compact Shimura varieties of $A_n$, $B_n$, $C_n$, and $D_n^{\dbR}$ type. The approach used also shows that many crystalline cycles on abelian varieties over finite fields which are specializations of Hodge cycles, are algebraic. These two results have many applications. First, we prove a variant of the conditional Langlands--Rapoport conjecture for these special fibres. Second, for certain isogeny sets we prove a variant of the unconditional Langlands--Rapoport conjecture (like for many basic loci). Third, we prove that integral canonical models of compact Shimura varieties of Hodge type that are of $A_n$, $B_n$, $C_n$, and $D_n^{\dbR}$ type, are closed subschemes of integral canonical models of Siegel modular varieties. 



\footline={\hss\tenrm \folio\hss}
\pageno=1

\bigskip
\noindent
{\boldsectionfont 1. Introduction}

\bigskip
Let $p\in\dbN$ be a prime. Let $\dbF$ be an algebraic closure of the field $\dbF_p$ with $p$ elements. Let $W(\dbF)$ be the ring of Witt vectors with coefficients in $\dbF$ and let $B(\dbF)$ be the field of fractions of $W(\dbF)$. In this paper we study abelian varieties over $\dbF$. In all that follows, by an {\it algebraic cycle} on an {\it abelian variety} $\bigtriangleup$ we mean a $\dbQ$--linear combination of irreducible subvarieties of $\sqcup_{n\in\dbN^*} \bigtriangleup^n$. The main two goals of the paper are:

\medskip
{\bf (a)} to show that many {\it crystalline cycles} on abelian varieties over $\dbF$ that are specializations of {\it Hodge cycles}, are algebraic (i.e., are crystalline realizations of algebraic cycles);

\smallskip
{\bf (b)} to get combinatorial descriptions of certain {\it isogeny classes} of principally polarized abelian varieties over $\dbF$ endowed with families of crystalline cycles.   

\medskip
Behind (a) and (b) are the {\it Tate/Hodge} and {\it Langlands--Rapoport conjectures} on special fibres of {\it integral canonical models} of {\it Shimura varieties} of {\it Hodge type}. In order to describe these two conjectures and our results pertaining to them, we will first review basic things on integral canonical models of Shimura varieties. 

\bigskip\noindent
{\bf 1.1. Shimura varieties of Hodge type.}  Let $(W,\psi)$ be a symplectic space over $\dbQ$. Let $W^{\vee}:=\Hom(W,\dbQ)$. Let $\scrS$ be the set of all monomorphisms $\Res_{\dbC/\dbR} \dbG_m\hookrightarrow \pmb{\text{GSp}}(W\otimes_{\dbQ} {\dbR},\psi)$ that define Hodge $\dbQ$--structures on $W$ of type $\{(-1,0),(0,-1)\}$ and that have either $2\pi i\psi$ or $-2\pi i\psi$ as polarizations. Let $(G,X)\hookrightarrow (\pmb{GSp}(W,\psi),S)$ be an injective map of Shimura pairs; thus $G$ is a reductive group over $\dbQ$ and $X$ is a $G(\dbR)$-conjugacy class of monomorphisms $\Res_{\dbC/\dbR} \dbG_m\hookrightarrow G_{\dbR}$ that define elements of $S$. Each such Shimura pair $(G,X)$ is called of Hodge type. Let $d\in \dbN^*$ be such that $\dim(W)=2d$. Let $(v_{\alpha})_{\alpha\in\scrJ}$ be a family of tensors of the tensor algebra $\scrT(\End(W))$ of $\End(W)=W\otimes_{\dbQ} W^{\vee}$ such that $G$ is the subgroup of $\pmb{\GL}_W$ that fixes each $v_{\alpha}$ with $\alpha\in\scrJ$, cf. [De3, Prop. 3.1 c)] and the fact that $G$ contains the center of $\pmb{\GL}_W$. Let $L$ be a $\dbZ$-lattice of $W$ such that we have a perfect alternating form $\psi:L\times L\to \dbZ$. We will assume that:

\medskip\noindent
{\bf (*)} {\it The schematic closure $G_{\dbZ_{(p)}}$ of $G$ in $\pmb{\GL}_{L\otimes_{\dbZ} \dbZ_{(p)}}$ is a reductive group scheme over $\dbZ_{(p)}$}. 

\medskip
Let $E(G,X)$ be the {\it reflex field} of $(G,X)$. Let $\Sh(G,X)$ be the {\it canonical model} over $E(G,X)$ of the complex Shimura variety defined by $(G,X)$, cf. [De1]. Let $\dbA_f:=\widehat{\dbZ}\otimes_{\dbZ}\dbQ$ be the ring of finite ad\`eles of $\dbQ$. Let $\dbA_f^{(p)}$ be the ring of finite ad\`eles of $\dbQ$ with the $p$-component omitted; we have $\dbA_f=\dbQ_p\times \dbA_f^{(p)}$. The group $G(\dbA_f)$ acts naturally on $\Sh(G,X)$. For a compact subgroup $\flat$ of $G(\dbA_f)$, let $\Sh_{\flat}(G,X)$ be the quotient of $\Sh(G,X)$ by $\flat$.

Let $K_p:=\pmb{GSp}(L,\psi)(\dbZ_p)$ and $H:=K_p\cap G(\dbQ_p)=G_{\dbZ_{(p)}}(\dbZ_p)$; they are {\it hyperspecial subgroups} of $\pmb{GSp}(L,\psi)(\dbQ_p)$ and $G(\dbQ_p)$ (respectively). We have identities $\Sh(G,X)(\dbC)=G(\dbQ)\backslash [X\times G(\dbA_f)]$ and $\Sh(\pmb{GSp}(W,\psi),S)(\dbC)=\pmb{GSp}(W,\psi)(\dbQ)\backslash [S\times \pmb{GSp}(W,\psi)(\dbA_f)]$ (cf. [De2, Cor. 2.1.11]) and a functorial closed embedding $\Sh(G,X)\hookrightarrow\Sh(\pmb{GSp}(W,\psi),S)_{E(G,X)}$ (cf. [De1, Cor. 5.4]). We also have an identity $\Sh_H(G,\scrX)(\dbC)=G_{\dbZ_{(p)}}(\dbZ_{(p)})\backslash (\scrX\times G(\dbA_f^{(p)}))$ (cf. [Mi3, Prop. 4.11]); from this and its analogue for $\Sh_{K_p}(\pmb{\text{GSp}}(W,\psi),S)(\dbC)$ we get that we also have a functorial closed embedding 
$$\Sh_{H}(G,X)\hookrightarrow \Sh_{K_p}(\pmb{GSp}(W,\psi),S)_{E(G,X)}.$$ 
\noindent
\medskip\noindent
{\bf 1.1.1. Hodge cycles.} We will use the terminology of [De3] on Hodge cycles on an abelian scheme $B$ over a reduced $\dbQ$--scheme $Y$. Thus we write each Hodge cycle $v$ on $B$ as a pair $(v_{\text{dR}},v_{\acute et})$, where $v_{\text{dR}}$ and $v_{\acute et}$ are the de Rham and the \'etale component of $v$ (respectively). The \'etale component $v_{\acute et}$ as its turn has an $l$-component $v_{\acute et}^l$, for each rational prime $l$. In what follows we will be interested only in Hodge cycles on $B$ that involve no Tate twists and that are tensors of different tensor algebras $\scrT(\sharp)$ (with $\sharp$ a finite dimensional vector space). Accordingly, if $Y$ is the spectrum of a field $E$, then in applications $v_{\acute et}^p$ will be a suitable $\Gal(\overline{E}/E)$-invariant tensor of $\scrT(\End(H^1_{\acute et}(B_{\overline{Y}},\dbQ_p)))$, where $\overline{Y}:=\Spec \overline{E}$. If moreover $\overline{E}$ is a subfield of $\dbC$, then we will also use the Betti realization $v_B$ of $v$: it is a tensor of $\scrT(\End(H^1((B\times \Spec\dbC)^{\an},\dbQ)))$ that corresponds to $v_{\text{dR}}$ (resp. to $v_{\acute et}^l$) via the canonical isomorphism that relates the Betti cohomology of the complex manifold $(B\times_Y \Spec \dbC)^{\an}$ with $\dbQ$--coefficients with the de Rham (resp. the $\dbQ_l$ \'etale) cohomology of $B_{\overline{Y}}$ (see [De3, Sect. 2]). We recall that $v_B$ is also a tensor of the $F^0$-filtration of the Hodge filtration of $\scrT(\End(H^1((B\times_Y \Spec \dbC)^{\an},\dbC)))$.

\medskip\noindent
{\bf 1.1.2. Integral canonical models.} Let $v$ be a prime of $E(G,X)$ that divides $p$. Let $O_{(v)}$ be the localization of the ring of integers of $E(G,X)$ at $v$. Property 1.1 (*) implies that $O_{(v)}$ is an \'etale $\dbZ_{(p)}$-algebra, cf. [Mi4, Cor. 4.7 (a)]. 

Let $\scrM$ be {\it Mumford's moduli scheme} over $\dbZ_{(p)}$ that parametrizes isomorphism classes of principally polarized abelian schemes of relative dimension $d$ over $\dbZ_{(p)}$-schemes that have compatibly level-$N$ symplectic similitude structures for all $N\in\dbN\setminus p\dbN$, cf. [MFK, Thms. 7.9 and 7.10] applied to symplectic similitude level structures instead of only to level structures. One identifies naturally $\Sh_{K_p}(\pmb{GSp}(W,\psi),S)=\scrM_{\dbQ}$ (cf. [De1, Prop. 4.17]) and thus one can speak about the schematic closure $\scrN^{\text{cl}}$ of $\Sh_{H}(G,X)$ in $\scrM_{O_{(v)}}$. 

Let $\scrN$ be the normalization of $\scrN^{\text{cl}}$. It is known that $\scrN$ is a regular, formally smooth $O_{(v)}$-scheme on which $G(\dbA_f^{(p)})$ acts continuously from the right and that $\scrN$ is the integral canonical model of $\Sh_{H}(G,X)$ over $O_{(v)}$ or of the Shimura quadruple $(G,X,H,v)$ in the strongest sense of [Va1, Def. 3.2.3 6)] (see [Va1, Thm. 6.4.1 and Subsubsect. 3.2.12] and [Va5] for $p\ge 5$ and see [Va14, Parts I and II] for $p$ arbitrary; the case $p\ge 3$ is also claimed in [Ki]). The morphism $\scrN\to\scrM_{O_{(v)}}$ (equivalently, $\scrN\to\scrN^{\text{cl}}$) is finite (see [Va1, proof of Prop. 3.4.1], [Va6, Prop. 3 (b)], or [Va14, Part I, Prop. 2.2.1 (b)]).

\medskip\noindent
{\bf 1.1.3. Moduli interpretation.} Let $(\grA,\lambda_{\grA})$ be the principally polarized abelian scheme over $\scrN$ which is the natural pull-back of the universal principally polarized abelian scheme over $\scrM$. For each $N\in\dbN\setminus p\dbN$, one has a natural level-$N$ symplectic similitude structure on $(\grA,\lambda_{\grA})$ that is defined by an isomorphism $\eta_N:(L/NL)_{\scrN}\arrowsim \grA[N]$ of constant \'etale group schemes over $\scrN$. If $N_1\in N\dbN\setminus p\dbN$, then $\eta_{N}$ is induced naturally by $\eta_{N_1}$.  

The choice of the $\dbZ$-lattice $L$ of $W$ and of the family of tensors $(v_\alpha)_{\alpha\in\scrJ}$ allows a moduli interpretation of $\Sh(G,X)$, cf. [De1,2], [Mi4], and [Va1, Subsect. 4.1 and Lem. 4.1.3]. For instance, $\Sh(G,X)(\dbC)=G(\dbQ)\backslash [X\times G(\dbA_f)]$ is the set of isomorphism classes of principally polarized abelian
varieties over $\dbC$ of dimension $d$, that carry a family of Hodge
cycles indexed by the set $\scrJ$, that have compatible level-$N$ symplectic similitude structures for all $N\in\dbN^*$, and that satisfy some additional axioms. This interpretation endows the abelian scheme $\grA_{E(G,X)}$ with a family $(w_{\alpha}^{\grA})_{\alpha\in\scrJ}$ 
of Hodge cycles; all realizations of pull-backs of $w_{\alpha}^{\grA}$ via $\dbC$-valued points of $\scrN_{E(G,X)}$ correspond naturally to $v_{\alpha}$.

\medskip\noindent
{\bf 1.1.4. On $\dbF$-valued points.} 
We consider a point $y:\Spec\dbF\to\scrN$. Let $(A,\lambda_A):=y^*(\grA,\lambda_{\grA})$. Let $(M,\phi,\psi_M)$ be the principally quasi-polarized (contravariant) {\it Dieudonn\'e module} of the principally quasi-polarized $p$-divisible group of $(A,\lambda_A)$.  For each $N\in\dbN\setminus p\dbN$, let $\eta_{N,y}:=y^*(\eta_N)$; it is an isomorphism $(L/NL)_{\dbF}\arrowsim A[N]$ of constant \'etale group schemes over $\dbF$ that defines a level-$N$ symplectic similitude structure on $(A,\lambda_A)$. 

For each $\alpha\in\scrJ$, let $t_{\alpha}$ be the tensor of the tensor algebra $\scrT(\End(M[{1\over p}]))$ of $\End(M[{1\over p}])$ which is the crystalline realization of the Hodge cycle $\tilde z_{V[{1\over p}]}^*(w_{\alpha}^{\grA})$ on the generic fibre of the lift $\tilde z ^*(\grA)$ of $A$ defined by a (any) lift $\tilde z:\Spec V\to\scrN$ of $y$, with $V$ as a  finite, discrete valuation ring extension of $W(\dbF)$. Each $t_{\alpha}$ depends only on $y$ and not on its lift $\tilde z$ (for instance, cf. [Va10, Sect. 5, p. 69]). We conclude that to the point $y:\Spec\dbF\to\scrN$ one associates naturally a septuple $$(A,\lambda_A,M,\phi,(t_{\alpha})_{\alpha\in\scrJ},\psi_M,(\eta_{N,y})_{N\in \dbN\setminus p\dbN}).$$ 
If $l\in\dbN^*$ is a prime different from $p$, then the sequence of isomorphisms $(\eta_{l^n,y})_{n\in\dbN^*}$ defines naturally a $\dbQ_l$-linear isomorphism $\eta_{\dbQ_l,y}:W\otimes_{\dbQ} \dbQ_l\arrowsim (\text{proj.}\text{lim.}_{n\in\dbN^*} A[l^n])\otimes_{\dbZ_l} \dbQ_l$ which maps $v_{\alpha}$ to the $l$-adic component of the \'etale component of the Hodge cycle $\tilde z_{V[{1\over p}]}^*(w_{\alpha}^{\grA})$ on the generic fibre of $\tilde z^*(\grA)$. If $a$ is a $\dbQ$--endomorphism of $A$, then the $l$-adic realization of $a$ will be denoted also by $a$ and it will be identified with the element $\eta_{\dbQ_l,y}^{-1} \circ a\circ \eta_{\dbQ_l,y}\in\pmb{GSp}(W,\psi)(\dbQ_l)$. In particular, it makes sense to say that $a$ is or is not an element of $G(\dbQ_l)$ (or of $G(\dbA_f^{(p)})$). 

Let $L^{\vee}:=\Hom(L,\dbZ)$. There exist isomorphisms $H^1_{\acute et}(\tilde z^*(\grA)_{V[{1\over p}]},\dbZ_p)\arrowsim L^{\vee}\otimes_{\dbZ} \dbZ_p$ under which the $p$-component of $\tilde z_{V[{1\over p}]}^*(w_{\alpha}^{\grA})$ is mapped to $v_{\alpha}$ for all $\alpha\in\scrJ$ and under which the perfect bilinear form on $H^1_{\acute et}((\tilde z^*(\grA)_{V[{1\over p}]},\dbZ_p)$ defined by $\tilde z^*(\lambda_{\grA})$ is mapped to a $\dbG_m(\dbZ_p)$-multiple of the perfect bilinear form $\psi^{\vee}$ on $L^{\vee}\otimes_{\dbZ} \dbZ_p$ defined by $\psi$ (see [Va1] or [Va14, Part I, Lem. 2.3.4 (a)]). Thus, if either $p>2$ or $p=2$ and the $2$-rank of $A$ is $0$, then there exists a $W(\dbF)$-linear isomorphism (cf. [Va11, Thm. 1.2 and Ex. 4.4.1]) 
$$(M,(t_{\alpha})_{\alpha\in\scrJ},\psi_M)\arrowsim (L^{\vee}\otimes_{\dbZ} W(\dbF),(v_{\alpha})_{\alpha\in\scrJ},\psi).\leqno (1)$$
In what follows, we will assume that:

\medskip\noindent
{\bf (**)} If $p=2$, then such an isomorphism (1) exists for each point $y:\Spec\dbF\to\scrN$ for which $A$ has positive $2$-rank.${}^1$ $\vfootnote{1}{Property (**) is proved in [Va14, Part II, Thm. 1.7 (b)]. Here we will only use it in the case when the Shimura pair $(G,X)$ is compact in the usual sense recalled before Theorem 1.5, cf. Lemma 2.5.6.}$

\medskip
Let $\scrG$ be the schematic closure in $\pmb{\GL}_M$ of the subgroup of $\pmb{\GL}_{M[{1\over p}]}$ that fixes each $t_{\alpha}$ with $\alpha\in\scrJ$. From (*), (1), and (**) we get that $\scrG$ is a reductive group scheme over $W(\dbF)$ isomorphic to $G_{\dbZ_{(p)}}\times_{\Spec \dbZ_{(p)}} \Spec W(\dbF)$. It is well known that under the natural action of $\phi$ on $\scrT(\End(M[{1\over p}]))$, each $t_{\alpha}$ is fixed by $\phi$ (see [Va10, Cor. 5.1.7]). Thus the Lie algebra $\Lie(\scrG)[{1\over p}]$ is left invariant by $\phi$ (i.e., $(\Lie(\scrG)[{1\over p}],\phi)$ is a Lie sub-$F$-isocrystal of $(\End(M)[{1\over p}],\phi)$). Each lift $z:\Spec W(\dbF)\to\scrN$ of $y$ defines naturally a direct summand $F^1$ of $M$ such that we have $\phi({1\over p}F^1+M)=M$ and $\psi(F^1,F^1)=0$ i.e., such that the quadruple $(M,F^1,\phi,\psi_M)$ is a principally quasi-polarized $F$-crystal over $\dbF$. More precisely, $F^1$ is the Hodge filtration of $(M,\phi)$ defined by $z^*(\grA)$. The inverse of the canonical split cocharacter of $(M,\phi)$ defined in [Wi, p. 512] fixes each $t_{\alpha}$ (cf. the functorial aspects of [Wi, p. 513]) and thus it factors through a cocharacter $\mu:\dbG_m\to \scrG$ which is a {\it Hodge cocharacter} of $(M,\phi,\scrG)$ in the sense of [Va10, Def. 1.2.1]. One gets that the triple $(M,\phi,\scrG)$ (resp. quadruple $(M,F^1,\phi,\scrG)$) is a {\it Shimura $F$-crystal} (resp. {\it filtered Shimura $F$-crystal}) over $\dbF$ in the sense of [Va8] and [Va10, Def. 1.2.1], to be said to be attached to $y:\Spec\dbF\to\scrN$. One says that $y:\Spec\dbF\to\scrN$ is a {\it basic point} if all {\it Newton polygon slopes} of $(\Lie(\scrG)[{1\over p}],\phi)$ are $0$. 

The following conjecture is only either an ad\`elic variation of the classical Tate conjecture for abelian varieties over $\dbF$ or a specialization form of the classical Hodge conjecture for abelian varieties over $\dbC$.

\bigskip\noindent
{\bf 1.2. Conjecture (Tate/Hodge conjecture for $\dbF$-valued points of integral canonical models of Shimura varieties of Hodge type).} {\it Let $y:\Spec\dbF\to\scrN$ and we use the above notations. We assume that $(v_{\alpha})_{\alpha\in\scrJ}$ is the family of all tensors of $\scrT(\End(W))$ fixed by $G$. Then  each $t_{\alpha}$ is the crystalline realization of an algebraic cycle $c_{\alpha,A}$ on $A$.}

\bigskip
Part of the philosophy of [Va13, Sect. 2.5]) can be captured as follows.

\bigskip\noindent
{\bf 1.3. Philosophy.} {\it To prove Conjecture 1.2 for all $\dbF$-valued points of integral canonical models of Shimura varieties of Hodge type, it suffices to prove Conjecture 1.2 for all $\dbF$-valued basic points of integral canonical models of Shimura varieties of Hodge type.}

\bigskip\noindent
{\bf 1.4. The isogeny set of $y$.} Let $\grI(y)$ be the set of those elements $h\in \scrG(B(\dbF))$ such that the pair $(h(M),\phi)$ is a Dieudonn\'e module over $k$ and there exists an element $h_1\in\scrG(W(\dbF))$ for which we have $\phi^{-1}(ph(M))=hh_1(\phi^{-1}(pM))$. For each element $h\in \grI(y)$, let $A(h)$ be the abelian variety over $\dbF$ which is $\dbZ[{1\over p}]$-isomorphic to $A$ and whose Dieudonn\'e module is canonically identified with $(h(M),\phi)$. We denote by $\lambda_{A(h)}$ the principal polarization of $A(h)$ defined naturally by $\lambda_A$; its crystalline realization $\psi_{h(M)}$ is $p^{n(h)}\psi_M$ for some integer $n(h)\in\dbZ$ (for instance, if $h$ fixes $\psi_M$, then $n(h)=0$). 

For each prime $l\neq p$, the Tate modules of $A(h)$ and $A$ are canonically identified and therefore each $\eta_{N,y}$ is as well a level-$N$ symplectic similitude structure on $(A(h),\lambda_{A(h)})$. Let $y(h):\Spec\dbF\to\scrM$ be the $\dbF$-valued point defined by $(A(h),\lambda_{A(h)})$ and the compatible level-$N$ symplectic similitude structures $\eta_{N,y}$. The following isogeny property was defined and announced to hold in [Va13] (Kisin also announced a proof of it in May 2008). 

\medskip\noindent
{\bf 1.4.1. Definition.} We say that the {\it isogeny property} holds for the point $y:\Spec\dbF\to\scrN$, if for each element $h\in\grI(y)$ the point $y(h):\Spec\dbF\to\scrM$ factors naturally  through a morphism $y(h):\Spec\dbF\to\scrN$ in such a way that for a (any) lift $z(h):\Spec W(\dbF)\to\scrN$ of it, every tensor $t_{\alpha}\in\scrT(\End(M[{1\over p}]))$ is the crystalline realization of the Hodge cycle $z(h)_{B(\dbF)}^*(w_{\alpha}^{\grA})$ on $z(h)^*(\grA)_{B(\dbF)}$.

\medskip 
We emphasize that if such a factorization  $y(h):\Spec\dbF\to\scrN$ exists, then it is unique (cf. [Va1, Rm. 5.6.4]; see Lemma 2.5.3). 

\medskip\noindent
{\bf 1.4.2. Definition.} If the isogeny property holds for $y:\Spec\dbF\to\scrN$, let the {\it isogeny set} $o(y)$ of $y$ be the set of $\dbF$-valued points of $\scrN$ that are $G(\dbA_f^{(p)})$-translate of points in the set $\{y(h)|h\in \grI(y)\}$. Let $o^{\text{big}}(y)$ be the isogeny set of the $\dbF$-valued of $\scrM$ defined by $y$. 

\medskip
Obviously, for each point $\tilde y\in o(y)$ we have $o(\tilde y)=o(y)$. Therefore, if one knows that the isogeny property holds for all $\dbF$-valued points of $\scrN$, then one gets a natural disjoint decomposition into $G(\dbA_f^{(p)})$-equivariant sets
$$\scrN(\dbF)=\bigsqcup_{o(y)\in \grK} o(y),\leqno (2)$$
where $\grK:=\{o(y)|y\in\scrN(\dbF)\}$. Let $\Theta$ be the Frobenius automorphism of $\dbF$ whose fixed field is the residue field $k(v)$ of $v$. It acts naturally on the set $\scrN(\dbF)$ in such a way that: (i) each isogeny set $o(y)$ is left invariant by $\Theta$ (see Corollary 2.3.5) and (ii) it commutes with the right action of $G(\dbA_f^{(p)})$ on the set $\scrN(\dbF)$. Thus, Formula (2) is in fact a decomposition into $\dbZ\Theta\times G(\dbA_f^{(p)})$-equivariant sets. 

\medskip\noindent
{\bf 1.4.3. Definition.} Let $\scrE$ be the schematic closure in the reductive group $\pmb{Aut}(A)_{\dbQ}$ over $\dbQ$ of invertible elements of $\End(A)\otimes_{\dbZ} \dbQ$ of the set of those $\dbQ$--automorphisms $a$ of $A$ that satisfy the following two properties:

\medskip
{\bf (i)} The crystalline realization of $a$ is a $B(\dbF)$-valued point of $\scrG$.

\smallskip
{\bf (ii)} For each prime $l\neq p$, the $l$-adic realization of $a$ is a $\dbQ_l$-valued point of $G$ (this makes sense, cf. second paragraph of Subsubsection 1.1.4). 

\medskip
Due to the property (ii), the $l$-adic realizations allow us to view naturally $\scrE(\dbQ)$ as a subgroup of $G(\dbA_f^{(p)})$. 

\medskip\noindent
{\bf 1.4.4. Definition.} We say that the {\it endomorphism property} holds for the point $y:\Spec\dbF\to\scrN$, if $\scrE$ is a reductive group over $\dbQ$ of the maximum possible dimension i.e., the crystalline realizations of elements of $\scrE(\dbQ)$ allow us to view $\scrE_{B(\dbF)}$ as a reductive subgroup of $\scrG_{B(\dbF)}$ with the property that we have an identity of Lie algebras over $\dbQ_p$
$$\Lie(\scrE)\otimes_{\dbQ} \dbQ_p=\{x\in\Lie(\scrG)[{1\over p}]|\phi(x)=x\}\cap (\End(A)^{\text{opp}}\otimes_{\dbZ} \dbQ_p).$$ 
\indent
Here $\End(A)^{\text{opp}}$ is the opposite $\dbZ$-algebra of $\End(A)$. Moreover, here and in all that follows, the action of $\scrE$ on $M[{1\over p}]$ is the composite of the standard isomorphism $\scrE\arrowsim \scrE^{\text{opp}}$ with the natural functorial action of the opposite group $\scrE^{\text{opp}}$ of $\scrE$ on $M[{1\over p}]$ induced by the functorial action of $\End(A)^{\text{opp}}$ on $M$. 

If $y:\Spec\dbF\to\scrN$ is a basic point, then the endomorphism property holds for it if and only if $\scrE$ is a $\dbQ$--form of $\scrG_{B(\dbF)}$ (cf. Subsection 2.2).

\medskip\noindent
{\bf 1.4.5. Definition.} We say the {\it CM lift property} holds for the point $y:\Spec\dbF\to\scrN$, if there exists an element $h\in\grI(y)$ for which the morphism $y(h):\Spec\dbF\to\scrM$ factors through $\scrN$ as in Definition 1.4.1 and for which there exist a finite, discrete valuation ring extension $V$ of $W(\dbF)$ and a lift $\tilde z(h):\Spec V\to\scrN$ of $y(h):\Spec\dbF\to\scrN$ whose generic fibre is a special point (i.e., the abelian scheme $\tilde z(h)^*(\grA)$ over $V$ has complex multiplication). If moreover:

\medskip
{\bf (a)} the Tate and Hodge conjectures hold for the generic fibre of $\tilde z(h)^*(\grA)$, then we say that the {\it Tate--Hodge property} holds for the point $y:\Spec\dbF\to\scrN$;

\smallskip
{\bf (b)} we can take $V=W(\dbF)$, then we say that the {\it unramified CM lift property} holds for the point $y:\Spec\dbF\to\scrN$.

\medskip\noindent
{\bf 1.4.6. On the Langlands--Rapoport conjecture.} We define
$$\grX(y):=\{h(M)|h\in\grI(y)\}=\text{Im}(\grI(y)\to \scrG(B(\dbF))/\scrG(W(\dbF))).$$ 
Due to the property 1.4.3 (i), the group $\scrE(\dbQ)$ acts via crystalline realizations on $\grX(y)$. Due to the isomorphism (1), one can naturally identify $\grX(y)$ with a certain subset of $G_{\dbZ_{(p)}}(B(\dbF))/G_{\dbZ_{(p)}}(W(\dbF))$ which is uniquely determined up to $G_{\dbZ_{(p)}}(W(\dbF))$-conjugation. Using such an identification, one gets a natural action of $\Theta$ on the set $\grX(y)$ (see Subsection 2.3). Thus $\Theta$ also acts on the product set $\grX(y)\times G(\dbA_f^{(p)})$ and on the quotient set $\scrE(\dbQ)\backslash [\grX(y)\times G(\dbA_f^{(p)})]$.

If the isogeny property holds for $y:\Spec\dbF\to\scrN$, then Faltings deformation theory implies that we have a natural bijection of $\dbZ\Theta\times G(\dbA_f^{(p)})$-equivariant sets
$$o(y)\arrowsim\scrE(\dbQ)\backslash [\grX(y)\times G(\dbA_f^{(p)})]\leqno (3)$$
with the property that $y$ is mapped to the equivalence class $[(1_M,1_W)]$ defined by $(1_M,1_W)\in \grX(y)\times G(\dbA_f^{(p)})$ (see Corollary 2.5.4 (a)). 
The $\dbZ\Theta$-equivariant part of the isomorphism (3) follows from its well known equivalent for the set $o^{\text{big}}(y)$ (see [Mi2]) and from the fact that the restriction of the natural map $o(y)\to o^{\text{big}}(y)$ to the set of equivalence classes of the form $[(h(M),1_W)]$ with $h\in\grI_p(y)$, is injective. If the isogeny property, the endomorphism property, and (in the case when $y$ is not basic) the CM lift property hold for $y$, then the natural map $o(y)\to o^{\text{big}}(y)$ is injective (cf. Corollary 2.5.4 (a)). 

There exists several forms of the Langlands--Rapoport conjecture in the literature with different levels of mistakes, impreciseness (due to some choices involved in constructing groupoids, categories of motives, and fibre functors), and generality (see [LR], [Mi3-5,7], [Pf], and [Re]). For applications to the zeta functions of Shimura varieties of Hodge type that involve a reductive group over $\dbQ$ whose derived group is simply connected, in essence all that is needed is an integral formula conjectured by Langlands and Kottwitz (see [La] and [Ko]) and known to follow from essentially any form of the Langlands--Rapoport conjecture (see [Mi3]). 

Based on the previous paragraph and on the fact that [Mi7] is already available, our approach to the Langlands--Rapoport conjecture will be a very practical one which is meant to provide necessary ingredients for the machinery of [Mi7] and which is automatically functorial in nature. Accordingly, we mention here that the Langlands--Rapoport conjecture predicts that the following three things hold.

\medskip
{\bf (a)} The isogeny property holds for every point $y:\Spec\dbF\to\scrN$ (thus Formulas (2) and (3) also hold). 

\smallskip
{\bf (b)} The endomorphism and the CM properties hold for every point $y:\Spec\dbF\to\scrN$.

\smallskip
{\bf (c)} The set $\grK$ is in natural bijection to the set of $G(\overline{\dbQ})$-conjugates (i.e., of isomorphism classes) $[\kappa]$ of {\it admissible homomorphisms} $\kappa:\grP\to \grG$ of $\overline{\dbQ}/\dbQ$--groupoids, where $\grP$ is the pseudo-motivic groupoid of $\dbF$ as in [LR], [Mi1-5], [Pf], [Re], and [Mi7, Sect. 6]  and where $\grG$ is the $\overline{\dbQ}/\dbQ$--groupoid defined by $G$. If $[\kappa]$ maps to $o(y)\in\grK$, then $\scrE(\dbQ)$ is naturally identified with the group of those elements of $G(\overline{\dbQ})$ which act on $\psi$ via a rational multiple and which fix $\kappa$ under inner conjugation.   

\medskip
More recently, the manuscript [Mi5] speaks about the unconditional form of the Langlands--Rapoport conjecture: briefly speaking, this is the same as above but in (c) one requires the admissibility to be in the strongest possible way (that determines the isomorphism classes uniquely). In this introduction  we will not be very precise in differentiating between the conditional form and the unconditional form of the Langlands--Rapoport conjecture. We will only say that the conditional form of the  Langlands--Rapoport conjecture and suitable forms of the Tate--Hodge property imply the unconditional form of the  Langlands--Rapoport conjecture.

We denote by $\star^{\ad}$ and $\star^{\der}$ the adjoint and the derived group schemes (respectively) of a reductive group scheme $\star$. For Shimura types we refer to [Mi4] and [Va3]. Let $(\tilde G^{\ad},\tilde X^{\ad})$ be the adjoint Shimura pair of a Shimura pair $(\tilde G,\tilde X)$, cf. [Va1, Subsubsect. 2.4.1]. The Shimura pair $(\tilde G,\tilde X)$ or its associated Shimura variety  $\Sh(\tilde G,\tilde X)$ is called compact if the $\dbQ$--rank of the adjoint group $\tilde G^{\ad}$ is $0$ (i.e., if and only if the canonical model $\Sh(\tilde G,\tilde X)$ is a pro-\'etale cover of a projective $E(\tilde G,\tilde X)$-scheme, cf. [BHC, Tm. 12.3 and Cor. 12.4]). 

We recall that the abelian variety $A$ is called {\it supersingular} if all Newton polygon slopes of $(M,\phi)$ are ${1\over 2}$ (equivalently, if all Newton polygon slopes of $(\End(M)[{1\over p}],\phi)$ are $0$). 

We have the following basic result which is the very essence of all our main results.

\bigskip\noindent
{\bf 1.5. Basic Theorem.} {\it Let $(G_0,X_0)$ be an arbitrary simple, adjoint Shimura pair of $A_n$, $B_n$, $C_n$, or $D_n^{\dbR}$ Shimura type. We assume that the group $G_{0,\dbQ_p}$ is unramified (i.e., it extends to a reductive group scheme over $\dbZ_p$). If $p=2$, then we also assume that the Shimura pair $\Sh(G_0,X_0)$ is compact (in fact this assumption can be removed based on [Va14], cf. footnote 1). Then there exists an injective map $(G,X)\hookrightarrow (\pmb{GSp}(W,\psi),S)$ of Shimura pairs such that the following eight properties hold (for each prime $v$ of $E(G,X)$ that divides $p$):

\medskip
{\bf (a)} We have $(G_0,X_0)=(G^{\ad},X^{\ad})$ and the derived group $G^{\der}$ is simply connected.

\smallskip
{\bf (b)} There exists a $\dbZ$-lattice $L$ of $W$ such that the property (*) holds (thus we can use the notations introduced above, with $v$ an arbitrary prime of $E(G,X)$ that divides $p$). If $p=2$, then the property (**) holds as well. 

\smallskip
{\bf (c)} Let $y:\Spec\dbF\to\scrN$ be a basic point. If $(G_0,X_0)$ is of $B_n$, $C_n$, or $D_n^{\dbR}$ type, then the abelian variety $A$ is supersingular and thus the Tate conjecture holds for $A$.

\smallskip
{\bf (d)} If $(G_0,X_0)$ is of $B_n$, $C_n$, or $D_n^{\dbR}$ type, then the Conjecture 1.2 holds for each basic point $y:\Spec\dbF\to\scrN$. If $(G_0,X_0)$ if of $A_n$ type and if a certain technical condition (***) (defined in Subsection 5.3 (f)) holds (resp. and the technical condition (***) does not hold), then the Conjecture 1.2 holds for each basic point $y:\Spec\dbF\to\scrN$ (resp. we can choose the family of tensors $(v_{\alpha})_{\alpha\in\scrJ}$ such that each $v_{\alpha}$ is an endomorphism of $W$).

\smallskip
{\bf (e)} The isogeny property, the endomorphism property, and the unramified CM lift property hold for all basic points $y:\Spec\dbF\to\scrN$ and in fact all $\dbF$-valued basic points of $\scrN$ form one isogeny set. If moreover $(G_0,X_0)$ is compact, then the isogeny property holds for all $\dbF$-valued points of $\scrN$.

\smallskip
{\bf (f)} There exists a $G(\dbA_f^{(p)})$-invariant, open subscheme $\scrL$ of $\scrN$ which contains all $\dbF$-valued basic points and which is as well an open subscheme of $\scrN^{\text{cl}}$. If moreover $(G_0,X_0)$ is compact, then in fact we have $\scrL=\scrN=\scrN^{\text{cl}}$.

\smallskip
{\bf (g) (due in essence to Zink and Milne)} We assume that $(G_0,X_0)$ is compact. If the endomorphism property and the CM lift property hold for all points $y:\Spec\dbF\to\scrN$, then the (conditional) Langlands--Rapoport conjecture holds for $\scrN$.

\smallskip
{\bf (h)} If $(G_0,X_0)$ if of $A_n$ type, then we assume that the technical condition (***) holds. Then the unconditional form of the Langlands--Rapoport conjecture holds for the isogeny set of basic points of $\scrN$.}

\medskip
Part (d) represents progress towards the proof of Conjecture 1.2. The guiding idea behind (d) is to check that, modulo ``relative PEL tricks'' as in [Va1, Subsubsect. 4.3.16], the Tate--Hodge property holds for one basic point of $\scrN$ (cf. Lemma 2.7 and Subsection 5.3 (b) and (c)). Part (f) is a completely new result (before it was known only when $(G_0,X_0)$ is of either $C_n$ type with $G_{0,\dbR}$ having no simple, compact factor or $A_n$ type; see [Zi], [Va3, Cor. 4.10], and [Va5, Prop. 3.2]). Part (g) follows from (e) and Formulas (2) and (3), cf. [Zi] and [Mi3-5,7]. In Section 6 we check that if $(G_0,X_0)$ is compact, then the endomorphism property and the unramified CM lift property hold as well for all points $y:\Spec\dbF\to\scrN$ (see Corollary 6.4 (a) whose essence is contained in the proof of Corollary 6.2). The proof of Corollary 6.2 involves all Shimura pairs $(G_1,X_1)$ of Hodge type that are compact and goes by {\it induction} on the dimension of the hermitian symmetric domain $X_1$ (this explains why (g) is stated as it is). 

For the sake of completeness, in this paragraph we assume that $(G_0,X_0)$ is not compact and we recall what is known in this case. If $(G_0,X_0)$ is of $A_n$ or $C_n$ type, then one can assume that the injective map $(G,X)\hookrightarrow (\pmb{GSp}(W,\psi),S)$ is a PEL type embedding. More precisely, one can assume that $G_{\dbZ_{(p)}}$ is the closed subgroup scheme of $\pmb{GSp}(L\otimes_{\dbZ} \dbZ_{(p)},\psi)$ that fixes each element of a semisimple $\dbZ_{(p)}$-subalgebra of $\End(L\otimes_{\dbZ} \dbZ_{(p)})$ (see Subsection 5.1 (d) and (f)). Thus if $(G_0,X_0)$ is of $A_n$ or $C_n$ type, then it is well known that the endomorphism property holds for all points $y:\Spec\dbF\to\scrN$ and in [Zi, Thm. 4.4] it is proved that a slightly weaker version of the isogeny and the CM lift properties hold for all points $y:\Spec\dbF\to\scrN$ (if the Hasse principle holds for $G$ --like if $(G_0,X_0)$ is of either $A_n$ type with $n$ odd or $C_n$ type--, then loc. cit. implies that the isogeny and the CM lift properties hold for all points $y:\Spec\dbF\to\scrN$). 

\bigskip\noindent
{\bf 1.6. On contents.} Section 2 contains preliminary material on tori, Faltings deformation theory, and on the Tate--Hodge property. In Section 3 we recall basic things on three stratifications of $\scrN_{k(v)}$. In Section 4 we present reduction steps towards the proof of the isogeny property. In Section 5 we prove Theorem 1.5. Section 6 contains our main applications. More precisely, in Section 6 we combine Theorem 1.5 with results of Milne and Pfau (see [Mi3-5,7] and [Pf]) to get general forms of Theorem 1.5 (e), (f), and (g) that do not depend on a particular choice of the injective map $(G,X)\hookrightarrow (\pmb{GSp}(W,\psi),S)$ and that do not require any extra assumption besides ``is compact''. In particular, we prove the conditional Langlands--Rapoport conjecture for all special fibres of integral canonical models of Shimura varieties of $A_n$, $B_n$, $C_n$, and $D_n^{\dbR}$ type that are compact (see Corollary 6.2). Many proofs (like those of Theorem 1.5 (c), (d), and (f)) also apply to large classes of Shimura varieties of $D_n^{\dbH}$ type (the so called unmixed $D_n^{\dbH}$ type).${}^2$ $\vfootnote{2}{In fact, due to the relative PEL trick of [Va3, Rm. 4.8.2 (b)] (to be compared with [Va13, Thm. 8.4]), all our results extend to the $D_n^{\dbH}$ type. The mentioned trick allows us to pass directly the endomorphism property from the $C_{2n}$ type to the $D_n^{\dbH}$ type (regardless of the fact that we are or we are not in a compact context).}$

The Appendix groups together different properties of Shimura $F$-crystals. The reader ought to refer to its Subsections A1 to A6 only when quoted. 

\bigskip\noindent
{\bf 1.7. Leitfaden.} For reader's convenience, we include the following diagram of implications between the different properties we have introduced in Subsection 1.4. We will abbreviate ``property'' by ``P'' and ``property for basic points'' by ``Pb''. Moreover the relevant subsections are listed on top of arrows. First, we have: 

$$
\spreadmatrixlines{1\jot}
\CD
Conjecture\;\; 1.2\;\; or\;\; Tate-Hodge\;\; P @>{4.2}>> Endomorphism\;\; P\\
@V{}VV @VV{}V\\
Conjecture\;\; 1.2\;\; or\;\;Tate-Hodge\;\; Pb @>{4.2}>> Endomorphism\;\; Pb
\endCD
$$
Second, if the derived group $G^{\der}$ of $G$ is simply connected, then we have:
$$Endomorphism\;\; Pb @>{4.3, 4.4, \grC_1}>> Isogeny\;\;Pb @>{4.1, \grC_2}>> Isogeny\;\;P @>{5.5,\;proof\;\;6.2}>> CM\;\;lift\;\; P.$$
Here $\grC_1$ refers to some group theoretical assumptions on the injective map $(G,X)\hookrightarrow (\pmb{GSp}(W,\psi),S)$ of Shimura pairs and on the existence of a basic point of $\scrN$ for which Conjecture 1.2 holds, while $\grC_2$ refers to the assumption that $(G,X)$ is compact. Third, the starting point for applying these implications is to show that:

\medskip\noindent
{\it Under good arrangements on the injective map $(G,X)\hookrightarrow (\pmb{GSp}(W,\psi),S)$, we can check that the Conjecture 1.2 (and even the Tate--Hodge property in many cases) holds for one basic point of $\scrN$ (cf. Theorem 1.5 (d) proved in Subsection 5.3).}

\medskip
See [Va13] for an overview of extra previous works on the Langlands--Rapoport conjecture performed by Ihara, Kottwitz, Reimann, and Milne.

\bigskip\smallskip
\noindent
{\boldsectionfont 2. Preliminaries}
\bigskip 

In this section we list our basic notations and recall basic material on tori, Faltings deformation theory, and the Tate--Hodge property. 

\bigskip\noindent
{\bf 2.1. Notations and conventions.} In all that follows, we will use the following notations $p$, $\dbF$, $W(\dbF)$, $B(\dbF)$, $(G,X)\hookrightarrow (\pmb{GSp}(W,\psi),S)$,  $d$, $(v_{\alpha})_{\alpha\in\scrJ}$, $L$, $G_{\dbZ_{(p)}}$, $E(G,X)$, $\Sh_{\flat}(G,X)$, $K_p=\pmb{GSp}(L,\psi)(\dbZ_p)$, $H=G_{\dbZ_{(p)}}(\dbZ_p)$, $v$, $O_{(v)}$, $\scrN^{\text{cl}}$, $\scrN$, $(\grA,\lambda_{\grA})$, $\eta_N:(L/NL)_{\scrN}\arrowsim \grA[N]$ with $N\in\dbN\setminus p\dbN$, $y:\Spec\dbF\to\scrN$, $(A,\lambda_A,M,\phi,(t_{\alpha})_{\alpha\in\scrJ},\psi_M,(\eta_{N,y})_{N\in\dbN\setminus p\dbN})$, $\eta_{\dbQ_l,y}$, $\scrG$, $\grI(y)$, $y(h):\Spec\dbF\to\scrM$ with $h\in\grI(y)$, $o(y)$, $\Theta$, $o^{\text{big}}(y)$, $k(v)$, $\scrE$, $\grP$, and $\grG$ we have introduced before Theorem 1.5. Always $l\in\dbN^*$ will be a prime different from $p$ and $\scrT(\sharp)$ will be the tensor algebra of some vector space $\sharp$. We will view $\pmb{Aut}(A)$ as an affine group scheme over $\dbZ$. Let $\scrE^+$ be the schematic closure in $\pmb{Aut}(A)_{\dbQ}$ of the set of those $\dbQ$--automorphisms $a$ of $A$ whose $l$-adic realizations are $\dbQ_l$-valued points of $G$ for all primes $l\neq p$; thus $\scrE$ is a subgroup of $\scrE^+$. Let $L_{(p)}:=L\otimes_{\dbZ} \dbZ_{(p)}$. 

If $k$ is a perfect field that contains $\dbF_p$, let $\sigma_k$ be the Frobenius automorphism of $k$, $W(k)$, and $B(k):=W(k)[{1\over p}]$. Let $\sigma:=\sigma_{\dbF}$. Let $Z(\boxdot)$ be the center of  a reductive group scheme $\boxdot$. Let $Z^0(\boxdot)$ be the largest torus of $Z(\boxdot)$; the finite group scheme $Z(\boxdot)/Z^0(\boxdot)$ is of multiplicative type. If $\boxdot$ is a reductive group scheme over an affine scheme $\Spec R$ and if $Q\to R$ is a finite, flat monomorphism, let $\Res_{R/Q} \boxdot$ be the group scheme over $\Spec Q$ obtained from $\boxdot$ through the Weil restriction of scalars (see [BT, Subsect. 1.5] and [BLR, Ch. 7, Sect. 7.6]). If moreover $R$ is an \'etale $Q$-algebra, then  $\Res_{R/Q} \boxdot$ is a reductive group scheme over $\Spec Q$.

For a point $y_{\star}:\Spec k\to\scrN$, let $(A_{\star},\lambda_{A_{\star}},M_{\star},\phi_{\star},(t_{\star,\alpha})_{\alpha\in\scrJ},\psi_{M_{\star}},(\eta_{N,y_{\star}})_{N\in\dbN\setminus p\dbN})$ be the septuple that is analogous to $(A,\lambda_A,M,\phi,(t_{\alpha})_{\alpha\in\scrJ},\psi_M,(\eta_{N,y})_{N\in\dbN\setminus p\dbN})$ but obtained working with $y_{\star}$ instead of with $y$. Let $\scrG_{\star}$ be the reductive subgroup scheme of $\pmb{\GL}_{M_{\star}}$ whose generic fibre is the subgroup of $\pmb{\GL}_{M_{\star}[{1\over p}]}$ that fixes $t_{\star,\alpha}$ for all $\alpha\in\scrJ$.

As the elements of $Z^0(G)(\dbQ)$ define $\dbQ$--endomorphisms of the universal abelian scheme $\grA$ over $\scrN$, we can identify canonically $Z^0(G)_{B(\dbF)}=Z^0(\scrG_{B(\dbF)})$. 

Let $(A_0,\lambda_{A_0})$ be a principally polarized  abelian variety model of $A$ over a finite field $\dbF_{p^q}$. Let $\pi_0$ be the Frobenius endomorphism of $A_0$. Let $(M_0,\phi_0,\psi_{M_0})$ be the principally quasi-polarized Dieudonn\'e module of $(A_0,\lambda_{A_0})$. By taking $q>>0$, we can assume that the following two properties hold:

\medskip
\item  {\bf (i)} For each $\alpha\in\scrJ$, we have $t_{\alpha}\in\scrT(M_0[{1\over p}])$ and therefore $\scrG$ is the pull-back of a reductive, closed subgrou scheme $\scrG_0$ of $\pmb{\GL}_{M_0}$. 

\smallskip
\item {\bf (ii)} The schematic closure of $\{\pi_0^s|s\in\dbZ\}$ in the $\dbQ$--group of invertible elements of the \'etale $\dbQ$--algebra $\dbQ[\pi_0]$ is connected i.e., is a torus (called the {\it Frobenius torus} of $\pi_0$).  

\medskip\noindent
The element $\pi_0:=\phi_0^q\in\pmb{\GL}_{M_0}(B(\dbF_{p^q}))$ is the crystalline realization of $\pi_0$ (thus the notations match) and fixes $t_{\alpha}$ for all $\alpha\in\scrJ$. Thus we have $\pi_0\in\scrG(B(\dbF))$. 

\bigskip\noindent
{\bf 2.2. On basic points.} The point $y:\dbF\to\scrN$ is basic if and only if $\pi_0\in Z^0(\scrG_{B(\dbF)})$ (cf. proof of [Va13, Thm. 8.3 (a)]) and thus if and only if $\pi_0\in Z^0(G)(\dbQ)$. Therefore, if $y$ is a basic point, then $\{x\in\Lie(\scrG)[{1\over p}]|\phi(x)=x\}\subset\End(A)^{\text{opp}}\otimes_{\dbZ} \dbQ_p$ and thus the endomorphism property holds for $y:\dbF\to\scrN$ if and only if we have $\Lie(\scrE)\otimes_{\dbQ} \dbQ_p=\{x\in\Lie(\scrG)[{1\over p}]|\phi(x)=x\}$ i.e., if and only if $\scrE$ is a $\dbQ$--form of $\scrG_{B(\dbF)}$. If $\scrE$ is a $\dbQ$--form of $\scrG_{B(\dbF)}$, then $\scrE_{\dbQ_p}$ splits over an unramified extension of $\dbQ_p$ and therefore $\scrE(\dbQ)$ is dense in $\scrE(\dbQ_p)$ (cf. [Mi4, Lem. 4.10]).  

\bigskip\noindent
{\bf 2.3. Action of $\Theta$.} Let $L_p:=L_{(p)}\otimes_{\dbZ_{(p)}} \dbZ_p$, $G_{\dbZ_p}:=G_{\dbZ_{(p)}}\times_{\Spec\dbZ_{(p)}} \Spec\dbZ_p$, and $G_{W(\dbF)}:=G_{\dbZ_{(p)}}\times_{\Spec\dbZ_{(p)}} \Spec W(\dbF)$. Based on Formula (1) we can identify 
$$(M,\phi,(t_{\alpha})_{\alpha\in\scrG},\psi_M)=(L_p^{\vee}\otimes_{\dbZ_p} W(\dbF),b_y(1\otimes\sigma),(v_{\alpha})_{\alpha\in\scrG},\psi),\leqno (4)$$
where $b_y\in G(B(\dbF))$. Let $\mu_y:\dbG_m\to G_{W(\dbF)}$ be a Hodge cocharacter of the Shimura $F$-crystal $(L_p^{\vee}\otimes_{\dbZ_p} W(\dbF),b_y(1\otimes\sigma),G_{W(\dbF)})$ over $\dbF$ (to be compared with end of Subsubsection 1.1.4). Under the identification (4), we can identify $\grX(y)$ with the subset of $G(B(\dbF))/G_{\dbZ_p}(W(\dbF))$ defined by those elements $g\in G(B(\dbF))$ with the properties that $g$ acts on $\psi$ via a $\dbG_m(\dbQ_p)$-multiple and that $g^{-1}b_y\sigma(g)\in G_{\dbZ_p}(W(\dbF))\mu_y({1\over p})G_{\dbZ_p}(W(\dbF))$.

Let $e_v\in\dbN^*$ be such that $k(v)$ has $p^{e_v}$ elements. 
Let $\Theta$ act on $G(B(\dbF))/G_{\dbZ_p}(W(\dbF))$ such that for $g\in G(B(\dbF))$, the element $[g]\in G(B(\dbF))/G_{\dbZ_p}(W(\dbF))$ is mapped to the element $[b_y\sigma(b_y)\cdots\sigma^{e_v-1}(b_y)\sigma^{e_v}(g)]\in G(B(\dbF))/G_{\dbZ_p}(W(\dbF))$. It is well known that $\Theta$ leaves invariant the subset $\grX(y)$ of $G(B(\dbF))/G_{\dbZ_p}(W(\dbF))$ (for instance, see [Mi3, p. 188]). Thus identification (4) defines naturally an action of $\Theta$ on $\grX(y)$.

\medskip\noindent
{\bf 2.3.1. On $o^{\text{big}}(y)$.} Let $\grI^{\text{big}}(y)$, $\scrE^{\text{big}}$, and $\grX^{\text{big}}(y)$ be the analogues of $\grI(y)$, $\scrE$, and $\grX(y)$ but obtained working with the composite point $y:\Spec\dbF\to\scrN\to\scrM$ instead of $y:\Spec\dbF\to\scrN$. Identification (4) also allows us to identify $\grX^{\text{big}}(y)$ with a subset of $\pmb{GSp}(W,\psi)(B(\dbF))/\pmb{GSp}(L_p,\psi)(W(\dbF))$ and to define in a similar way an action of $\Theta$ on $\break\pmb{GSp}(W,\psi)(B(\dbF))/\pmb{GSp}(L_p,\psi)(W(\dbF))$ which extends the action of $\Theta$ on $G(B(\dbF))/G_{\dbZ_p}(W(\dbF))$ and which (cf. loc. cit.) leaves invariant $\grX^{\text{big}}(y)$. The Langlands--Rapoport conjecture is known to hold for Siegel modular varieties,  cf. [Mi2]. In particular, we have a natural $\dbZ\Theta\times \pmb{GSp}(W,\psi)(\dbA_f^{(p)})$-equivariant bijection 
$$o^{\text{big}}(y)\arrowsim\scrE^{\text{big}}(\dbQ)\backslash [\grX^{\text{big}}(y)\times \pmb{GSp}(W,\psi)(\dbA_f^{(p)})]\leqno (5)$$
which for $h^{\text{big}}\in\grI^{\text{big}}(y)$ maps $y(h^{\text{big}})$ to the equivalence class $[(h^{\text{big}}(M),1_W)]$.

\bigskip\noindent
{\bf 2.4. Complements on tori.} We recall that a torus over a field is called {\it anisotropic}, if it has no subtorus isomorphic to $\dbG_m$. Let $T$ be a torus over $\dbQ$ such that we have a short exact sequence $0\to\dbG_m\to T\to T_c\to 0$, where the torus $T_{c,\dbR}$ is compact. Thus $T$ has a unique subtorus isomorphic to $\dbG_m$. We assume that we have a homomorphism $h:\Res_{\dbC/\dbR} \dbG_m\to T_{\dbR}$ such that the subtorus $\dbG_m$ of $\Res_{\dbC/\dbR} \dbG_m$ maps onto the subtorus $\dbG_m$ of $T$. In this subsection we are interested in the smallest subtorus $\tilde T$ of $T$ such that $h$ factors through $\tilde T_{\dbR}$. Let $T_1$ be a subtorus of $T$ that contains the $\dbG_m$ subtorus of $T$. Let $T_2$ be a subtorus of $T$ which is naturally isogenous to $T/T_1$ (the below assumptions will imply that such a subtorus is unique). We assume that there exists a number field $F_0$ and a torus $T_2^{F_0}$ over $F_0$ such that we have an isogeny 
$$T_2\to\Res_{F_0/\dbQ} T_2^{F_0}.$$
If $l$ is a prime that splits in $F_0$, then we have $(\Res_{F_0/\dbQ} T_2^{F_0})_{\dbQ_l}=\prod_{i\in\Hom(F_0,\dbQ_l)}T_2^{F_0}\times_{F_0,i} \dbQ_l$. 

\medskip\noindent
{\bf 2.4.1. Lemma.} {\it We assume that the following three conditions hold:

\medskip
{\bf (i)} the torus $\tilde T$ surjects onto $T/T_2$ and its image in $T/T_1$ is non-trivial;

\smallskip
{\bf (ii)} there exists a prime $l$ which splits in $F_0$, for which the torus $T_{1,\dbQ_l}$ is split, and which has the property that for each $i\in\Hom(F_0,\dbQ_l)$ the torus $T_2^{F_0}\times_{F_0,i} \dbQ_l$ is anisotropic and has no proper subtorus;

\smallskip
{\bf (iii)} either $F_0=\dbQ$ or $F_0\neq \dbQ$ and there exists an element $i_0\in\Hom(F_0,\dbQ_l)$ such that the torus $T_2^{F_0}\times_{F_0,i_0} \dbQ_l$ is not isogenous to any other torus of the form $T_2^{F_0}\times_{F_0,i} \dbQ_l$ with $i\in\Hom(F_0,\dbQ_l)\setminus\{i_0\}$.

\medskip
Then we have $\tilde T=T$.}

\medskip
\proof
We can assume that $T=T_1\times_{\dbQ} T_2$ and that the isogeny $T_2\to\Res_{F_0/\dbQ} T_2^{F_0}$ is an isomorphism. As $T=T_1\times_{\dbQ} T_2$, we have natural projection homomorphisms $T\to T_1$, $T\to T_2$, $\tilde T\to T_1$, and $\tilde T\to T_2$. Let $T_0$ be the identity component of $\Im(\Ker(\tilde T\to T_1)\to T_2)$. 

We first assume that $F_0=\dbQ$. From the condition (i) and the last part of the condition (ii) we get that $\tilde T$ surjects onto both $T_1$ and $T_2$. From the condition (ii) we also get that $T_0$ is non-trivial and in fact that it is $T_2$ itself. Thus $\tilde T=T_1\times_{\dbQ} T_2=T$. 

We assume now that $F_0\neq \dbQ$. From conditions (i) to (iii) we get that $T_0$ is non-trivial and moreover $T_{0,\dbQ_l}$ contains the torus $T_2^{F_0}\times_{F_0,i_0} \dbQ_l$ (viewed as a direct factor of $T_{2,\dbQ_l}=(\Res_{F_0/\dbQ} T_2^{F_0})_{\dbQ_l}=\prod_{i\in\Hom(F_0,\dbQ_l)}T_2^{F_0}\times_{F_0,i} \dbQ_l$ and thus of $T_{\dbQ_l}$). As $T_0$ is a subtorus of $T_2=\Res_{F_0/\dbQ} T_2^{F_0}$, we get that $T_0=T_2$. Thus $\tilde T$ contains $T_2$. From this and the condition (i) we get that $\tilde T=T$.\endproof

\medskip\noindent
{\bf 2.4.2. Lemma.} {\it Let $\scrT_{\dbZ_p}$ be a torus over $\dbZ_p$ whose generic fibre $\scrT_{\dbQ_p}$ is a maximal torus of $G_{\dbQ_p}$. Then there exists $g\in G(\dbQ_p)$ such that $\scrT_{\dbZ_p}$ is a maximal torus of $gG_{\dbZ_p} g^{-1}$.}

\medskip
\proof
This is a particular case of [Va13, Lem. 2.3].\endproof

\bigskip\noindent
{\bf 2.5. Faltings deformation theory.} Let $y:\Spec\dbF\to\scrN$. Let $z:\Spec W(\dbF)\to\scrN$ be a lift of $y$. Let $F^1\subset M$ and $\mu:\dbG_m\to\scrG$ be as in Subsubsection 1.1.4. For $c\in\dbN$, let $\scrR_c:=W(\dbF)[[x_1,\ldots,x_c]]$. Let $\Phi_c$ be the Frobenius endomorphism of $\scrR_c$ which is compatible with $\sigma$ and which takes $x_i$ to $x_i^p$ for all $i\in\{1,\ldots,c\}$. Let $d\Phi_c:\oplus_{i=1}^c \scrR_cdx_i\to \oplus_{i=1}^c \scrR_cdx_i$ be the $(x_1,\ldots,x_c)$-adic completion of the differential of $\Phi_c$.

Let $b:=\dim(X)$ be the dimension of $X$ viewed as a complex manifold. Let $\scrV$ be a smooth, closed subscheme of $\scrG$ which is isomorphic to $\Spec W(\dbF)[x_1,\ldots,x_b]$, which contains the identity section of $\scrG$, and for which we have an identity $\Lie(\scrG_{\dbF})=\Lie(\scrV_{\dbF})\oplus\Lie(\scrP_{\dbF})$ of $\dbF$-vector spaces, where $\scrP_{\dbF}$ is the normalizer of $F^1/pF^1$ in $\scrG_{\dbF}$ and where by $\Lie(\scrV_{\dbF})$ we mean the tangent space of $\scrV_{\dbF}$ at the identity element (in practice, this tangent space is an abelian Lie subalgebra of $\Lie(\scrG_{\dbF})$ and this is why we denote it as a Lie algebra). For instance, we can take $\scrV$ to be the $\dbG_a^b$ subgroup scheme of $\scrG$ whose Lie algebra is the largest direct summand of $\Lie(\scrG)$ on which $\dbG_m$ acts through $\mu$ via the identical character of $\dbG_m$ (see [Va14, Part I, Subsects. 3.2 to 3.4]). Let $\Spf \scrR$ be the completion of $\scrV$ along its identity section. Let $u_{\text{univ}}:\Spec \scrR\to \scrV$ be the universal (natural) morphism. We fix an identification $\scrR=\scrR_b=W(\dbF)[[x_1,\ldots,x_b]]$. Let $\Phi:=\Phi_b$; it is a Frobenius endomorphism of $\scrR$. The completion of the local ring of $y$ in $\scrN$ can be identified with $\scrR$. Let $\scrY:\Spec\scrR\to\scrN$ be the natural formally-\'etale morphism. Faltings deformation theory asserts three basic things (cf. [Fa, Thm. 10 and rmks. after], [Va1, Subsect. 5.4], [Va14, Part I, Subsects. 3.3 and 3.4]).

\medskip
{\bf (a)} First it assets that the principally quasi-polarized filtered $F$-crystal with tensors  over $\Spec \scrR/p\scrR$ associated naturally to $\scrY:\Spec \scrR\to\scrN$ is (up to $W(\dbF)$-automorphisms of $\scrR$ that leave invariant the ideal $(x_1,\ldots,x_b)$) isomorphic to
$$\scrC_z:=(M\otimes_{W(k)} \scrR,F^1\otimes_{W(k)} \scrR,u_{\text{univ}}(\phi\otimes\Phi),\nabla,(t_{\alpha})_{\alpha\in\scrJ},\psi_M),$$
where $\nabla$ is an integrable, topologically nilpotent, versal connection on $M\otimes_{W(k)} \scrR$ which is uniquely determined by the identity $\nabla\circ u_{\text{univ}}(\phi\otimes\Phi)=(u_{\text{univ}}(\phi\otimes\Phi)\otimes d\Phi_b)\circ\nabla$ and which annihilates $t_{\alpha}$ for all $\alpha\in\scrJ$ (see [Va1, Subsect. 5.4] and [Va14, Part I, Subsects. 3.3 and 3.4]). 
\smallskip
{\bf (b)} Second it asserts that each quintuple 
$$\tilde \scrC:=(\tilde M,\tilde F^1,\tilde\phi,\tilde\nabla,(\tilde t_{\alpha})_{\alpha\in\scrJ},\psi_{\tilde M})$$ 
over $\scrR_c/p\scrR_c$ which modulo $(x_1,\ldots,x_c)$ is $(M,F^1,\phi,(t_{\alpha})_{\alpha\in\scrJ},\psi_M)$ and which satisfies the same properties as $\scrC_z$, is the pull-back of $\scrC_z$ via a unique morphism $\Spec \scrR_c\to\Spec \scrR$ of $W(\dbF)$-schemes which at the level of ring homomorphisms takes the ideal $(x_1,\ldots,x_b)$ of $\scrR$ to the ideal $(x_1,\ldots,x_c)$ of $\scrR_c$. Here by $\tilde\scrC$ being over $\scrR_c/p\scrR_c$ and by satisfying the same properties as $\scrC_z$, we mean that $\tilde M$ is a free $\scrR_c$-module of rank $2d$, that $\psi_{\tilde M}$ is a perfect alternating form on $\tilde M$, that $\tilde F^1$ is a direct summand of $\tilde M$ which has rank $d$ and which is anisotropic with respect to $\psi_{\tilde M}$, that $\tilde\phi:\tilde M\to\tilde M$ is a $\Phi_c$-linear endomorphism such that $\tilde M$ is $\scrR_c$-generated by $\tilde\phi({1\over p}\tilde F^1+\tilde M)$ and for which we have $\psi_{\tilde M}(\tilde\phi(u),\tilde\phi(v))=p\Phi(\psi_{\tilde M}(u,v))$ for all $u,v\in\tilde M$, that $\tilde\nabla$ is an integrable, topologically nilpotent connection on $\tilde M$ which is uniquely determined by the identity $\tilde\nabla\circ \tilde\phi=(\tilde\phi\otimes d\Phi_c)\circ\tilde\nabla$, and that for each $\alpha\in\scrJ$ the tensor $\tilde t_{\alpha}\in\scrT(\End(\tilde M[{1\over p}]))$ is fixed by $\tilde\phi$ and belongs to the $F^0$-filtration of $\scrT(\End(\tilde M[{1\over p}]))$ defined by $\tilde F^1$ (see [Va14, Part I, Appendix, Thm. B6.4]). Each $\tilde t_{\alpha}$ is annihilated by $\tilde\nabla$, cf. also [Va14, Part I, Appendix, property B6.3 (i)].
\smallskip
{\bf (c)} Let $k$ be a perfect field of characteristic $p$. Third it asserts that both (a) and (b) continue to hold for an arbitrary $k$-valued point of $\scrN$. 

\medskip
In what follows we will need the following six applications of Faltings deformation theory. From the versality of $\nabla$ (see (a)) we get directly the following consequence which is a particular case of [Va14, Part I, Thm. 1.5 (b)].

\medskip\noindent
{\bf 2.5.1. Corollary.} {\it The morphism $\scrN_{W(\dbF)}\to\scrM_{W(\dbF)}$ induces formally closed embeddings at the level of complete local rings of residue $\dbF$.} 

\medskip\noindent
{\bf 2.5.2. Simple properties.} The morphism $\bar\scrY:\Spec \scrR/p\scrR\to\scrN\to\scrM$ induced by $\scrY$, depends only on $y$ and not on $z$. This holds even if $p=2$, due to Serre--Tate deformation theory and the fact that the quadruple $(M\otimes_{W(k)} \scrR,u_{\text{univ}}(\phi\otimes\Phi),\nabla,\psi_M)$ is the principally quasi-polarized $F$-crystal of a unique (up to isomorphism) principally quasi-polarized $p$-divisible group over $\Spec \scrR/p\scrR$ (cf. the fully faithfulness part of [BM, Thm. 4.1]). Similarly, regardless of what $p$ is, there exists a unique morphism $\Spec \scrR_c/p\scrR_c\to\scrN\to\scrM$ which lifts $y:\Spec\dbF\to\scrN\to\scrM$ and which has the property that the principally quasi-polarized $F$-crystal over $\scrR_c/p\scrR_c$ of the pull-back to $\Spec \scrR_c/p\scrR_c$ of $(\grA,\lambda_{\grA})$ is $(\tilde M,\tilde\phi,\tilde\nabla,\psi_{\tilde M})$. The same properties hold for $k$-valued points of $\scrN$.

\medskip\noindent
{\bf 2.5.3. Lemma.} {\it Let $y,y_1:\Spec\dbF\to\scrN$ be two points that define the same $\dbF$-valued point of $\scrN^{\text{cl}}$. The quadruple which is attached to $y_1$ and which is the analogue of the quadruple $(M,\phi,(t_{\alpha})_{\alpha\in\scrJ},\psi_{M})$, is of the form $(M_1,\phi_1,(t_{1,\alpha})_{\alpha\in\scrJ},\psi_{M_1})=(M,\phi,(t_{1,\alpha})_{\alpha\in\scrJ},\psi_{M})$. We assume that $\scrG_1=\scrG$ i.e., the subgroup of $\pmb{\GL}_{M[{1\over p}]}$ that fixes $t_{\alpha}$ for all $\alpha\in\scrJ$ is the same as the subgroup of $\pmb{\GL}_{M[{1\over p}]}$ that fixes $t_{1,\alpha}$ for all $\alpha\in\scrJ$ (for instance, this holds if we have $t_{\alpha}=t_{1,\alpha}$ for all $\alpha\in\scrJ$. Then we have $y=y_1$.}

\medskip
\proof
Let $z_1:\Spec W(\dbF)\to\scrN$ be a lift of $y_1$ such that the Hodge filtration of $M$ defined by the abelian scheme $z_1^*(\grA)$ over $W(\dbF)$ which lifts $A$, is $F^1$ (cf. [Va14, Part I, Lem. 3.5.2]).  

We first assume that either $p>2$ or $p=2$ and the $2$-rank of $A$ is $0$. This assumption implies that the $p$-divisible groups of $z^*(\grA)$ and $z_1^*(\grA)$ are canonically identified as lifts of the $p$-divisible group of $A$ (cf. [Va11, Prop. 2.2.6]) and thus $z$ and $z_1$ induce the same $W(\dbF)$-valued point of $\scrN^{\text{cl}}$. As $\scrN$ is the normalization of $\scrN^{\text{cl}}$ and as $W(\dbF)$ is a normal ring, we get that the two $W(\dbF)$-valued points $z$ and $z_1$ of $\scrN$ coincide. This implies that $y=y_1$. 

We assume that $p=2$ and that the $2$-rank of $A$ is positive. Let $a$ be the multiplicity of the Newton polygon slope $-1$ for $(\Lie(\scrG)[{1\over p}],\phi)=(\Lie(\scrG_1)[{1\over p}],\phi)$. Let $b$ be the number of points $y_1:\Spec\dbF\to\scrN$ with the property that $\scrG_1=\scrG$. We can count the number $N_y(F^1)$ of $W(\dbF)$-valued points of $\scrN^{\text{cl}}$ which lift the $\dbF$-valued point of $\scrN^{\text{cl}}$ defined by either $y$ or $y_1$ and whose filtered $F$-crystal is $(M,F^1,\phi)$. On one hand, due to the assumption 1.1.4 (**) we have $N_y(F^1)=2^a$  (cf. [Va14, Part I, Appendix, Thm. B7 (c)]) and on the other hand it is at least $b2^a$ (cf. loc. cit. and property 2.5 (a); more precisely, each one of the $b$ possibilities of $y_1$'s produces $2^a$ distinct such $W(\dbF)$-valued points of $\scrN^{\text{cl}}$). Thus if $y\neq y_1$, then $b\ge 2$ and we reached a contradiction. Therefore we have $b=1$ and $y=y_1$.\endproof

\medskip\noindent
{\bf 2.5.4. Corollary.} {\it We assume that the isogeny property holds for $y:\Spec\dbF\to\scrN$.

\medskip
{\bf (a)} Then the subset $o(y)$ of $\scrN(\dbF)$ is invariant under the action of $\Theta$ and the Formula (3) holds.

\smallskip
{\bf (b)} Then $\scrE$ is an open subgroup of the group $\scrE^+$ of Subsection 2.1.}

\medskip
\proof
To check that $\Theta(o(y))=o(y)$, it suffices to show that for each $h\in \grI(y)$ we have $\Theta(y(h))\in o(y)$. Let $\Theta(h)\in \grI(y)$ be such that $\Theta([h(M),1_W])=[\Theta(h)(M),1_W]\in o^{\text{big}}(y)$, cf. Formula (5) (the inclusion $G(B(\dbF))\backslash G_{\dbZ_p}(W(\dbF))\subset \pmb{GSp}(W,\psi)(B(\dbF))/\pmb{GSp}(L_p,\psi)(W(\dbF))$ is invariant under the actions of $\Theta$ of Subsection 2.3). Thus $\Theta(y(h))=y(\Theta(h))$ (cf. Lemma 2.5.3) and therefore $\Theta(y(h))\in o(y)$ (cf. hypothesis).

Let $y_1,y_2\in o(y)$. We can identify $M[{1\over p}]=M_1[{1\over p}]=M_2[{1\over p}]$ and $t_{\alpha}=t_{1,\alpha}=t_{2,\alpha}$ for all $\alpha\in\scrJ$. Let $h_1,h_2\in \grI(y)$ and $t_1,t_2\in G(\dbA_f^{(p)})$ be such that $y_1=y(h_1)t_1$ and $y_2=y(h_2)t_2$. The images of $y_1$ and $y_2$ in $\scrM(\dbF)$ (equivalently in $o^{\text{big}}(y)$) coincide if and only if there exists an element $a\in\scrE^{\text{big}}(\dbQ)$ such that we have $ay(h_1)=ay(h_2)$ and $at_1=t_2$, cf. Subsubsection 2.3.1. Based on this and Lemma 2.5.3, we get that we have $y_1=y_2$ if and only if we can choose $a\in\scrE^{\text{big}}(\dbQ)$ such that moreover the crystalline realization of $a$ fixes each $t_{\alpha}$ with $\alpha\in\scrJ$. If $a\in\scrE^{\text{big}}(\dbQ)$ is such that $at_1=t_2$, then the $l$-adic realizations of $a$ allow us to view $a$ as an element of $G(\dbA_f^{(p)})$ and therefore we have $a\in\scrE(\dbQ)$ if and only if the crystalline realization of $a$ is an element of $\scrG(B(\dbF))$. Based on the last two sentences we get that the identity $y_1=y_2$ holds if and only if there exist an element $a\in\scrE(\dbQ)$ such that we have $ay(h_1)=ay(h_2)$ and $at_1=t_2$. From this the Formula (3) follows. Thus (a) holds. 

To check (b), we first remark that obviously $\scrE$ is a subgroup of $\scrE^+$. Let $\Gamma^+$ (resp. $\Gamma$) be the subgroup of $\scrE^+(\dbQ)$ (resp. $\scrE(\dbQ)$) which normalizes $M$; it is the group of $\dbZ_{(p)}$-valued points of the schematic closure of $\scrE^+$ (resp. of $\scrE$) in $\pmb{Aut}(A)_{\dbZ_{(p)}}$. As the morphism $\scrN\to\scrN^{\text{cl}}$ is finite, the fibres of the map $\scrN(\dbF)\to\scrN^{\text{cl}}(\dbF)$ are finite. Thus the natural map $o(y)\to o^{\text{big}}(y)$ has finite fibres. From this and the fact that Formulas (3) and (5) hold (cf. (a) and Subsubsection 2.3.1), we get that the quotient set $\Gamma^+/\Gamma$ is finite. As $\Gamma^+$ is Zariski dense in an open subgroup of $\scrE^+$, we conclude that $\scrE$ is an open subgroup of $\scrE^+$.\endproof

\medskip\noindent
{\bf 2.5.5. Corollary.} {\it We assume that the isogeny property and the endomorphism property hold for $y:\Spec\dbF\to\scrN$. We also assume that either $y$ is a basic point or the CM lift property holds for $y$ as well. Then the natural map $o(y)\to o^{\text{big}}(y)$ is injective.}

\medskip
\proof
Based on the proof of Corollary 2.5.4 (a), it suffices to show that if $a\in\scrE^{\text{big}}(\dbQ)$ is such that $a\in G(\dbA_f^{(p)})$, then $a\in\scrE(\dbQ)$. Let $A_0$, $\dbF_{p^q}$, and $\pi_0$ be as in Subsection 2.1. If the CM lift property holds for $y:\Spec\dbF\to\scrN$ (resp. if $y$ is a basic point), then we can view naturally $\pi_0$ as an element of $G(\dbQ)$ (resp. of $Z^0(G)(\dbQ)$) in such a way that its centralizers in either $G_{\dbQ_l}$ or $\scrG_{B(\dbF)}$ are reductive groups of equal dimension (the reductiveness part follows from the property 2.1 (ii) and the fact that the centralizer of a torus in a reductive group is reductive). But as the endomorphism property holds for $y$, this dimension is precisely $\dim(\scrE)$. From this we get for each prime $l\neq p$, the group of those $\dbQ$--automorphisms of $A$ whose $l$-adic realizations belongs to $G(\dbQ_l)$, is $\scrE(\dbQ)$. Thus $a\in\scrE(\dbQ)$.\endproof

\medskip\noindent
{\bf 2.5.6. Lemma.} {\it We assume that $p=2$ and that $(G,X)$ is compact. Then the condition 1.1.4 (**) holds.} 

\medskip
\proof
This is a particular case of [Va14, Part I, Thm. 1.7 (c)].\endproof 

\bigskip\noindent
{\bf 2.6. Deformations of isogenies.} Let $h\in\grI(y)$. Let $\scrY(h):\Spec \scrR\to\scrM$ be a morphism that is constructed similarly to the composite morphism $\scrY:\Spec \scrR\to\scrN\to\scrM$ which lifts the point $y(h):\Spec\dbF\to\scrM$ and which modulo the ideal $(x_1,\ldots,x_b)$ of $\scrR$ produces a point $z(h):\Spec W(\dbF)\to\scrM$ with the property that the direct summand $F^1(h)$ of $h(M)$ defined by the Hodge filtration of the abelian scheme over $W(\dbF)$ that lifts $A(h)$ and that is naturally associated to $z(h)$, is such that $(h(M),F^1(h),\phi,h\scrG h^{-1})$ is a filtered Shimura $F$-crystal over $\dbF$. Thus the principally quasi-polarized filtered $F$-crystal with tensors over $\Spec \scrR/p\scrR$ associated to $\scrY(h)$ and to $z(h)$, is (up to $W(\dbF)$-automorphisms of $\scrR$ that leave invariant the ideal $(x_1,\ldots,x_b)$) of the form 
$$\scrC_z(h):=(h(M)\otimes_{W(k)} \scrR,F^1(h)\otimes_{W(k)} \scrR,u_{\text{univ}}(h)(\phi\otimes\Phi),\nabla(h),(t_{\alpha})_{\alpha\in\scrJ},\psi_M),$$ where the direct summand $F^1(h)$ of $h(M)$ is of the form $hh_1(F^1)$ for some element $h_1\in \scrG(W(k))$ and where $u_{\text{univ}}(h)\in h\scrG h^{-1}(\scrR)$ is a suitable universal element associated to the closed subscheme $\scrV(h):=hh_1\scrV h_1^{-1}h^{-1}$ of $h\scrG h^{-1}$ (or of $\pmb{\GL}_{h(M)}$). 

We assume that there exists a geometric point $\tilde y_1:\Spec k\to\Spec \scrR/p\scrR\to\Spec\scrR$ such that $y_1:=\scrY(h)\circ \tilde y_1:\Spec k\to\scrM$ factors through $\scrN$ in such a way that the pull-back of $\scrC_z(h)$ via it is the quadruple $(M_1,\phi_1,(t_{1,\alpha})_{\alpha\in\scrJ},\psi_{M_1})$ attached to the resulting factorization $y_1:\Spec k\to\scrN$. Let $\scrG_1$ be the reductive subgroup scheme of $\pmb{\GL}_{M_1}$ which is the schematic closure of the subgroup of $\pmb{\GL}_{M_1[{1\over p}]}$ that fixes $t_{1,\alpha}$ for all $\alpha\in\scrJ$. 
Let $\bar\scrY(h):\Spec \scrR/p\scrR\to\scrM$ be the morphism defined naturally by $\scrY(h)$. 

\medskip\noindent
{\bf 2.6.1. Fact.} {\it The morphism $\bar\scrY(h)$ factors uniquely through $\scrN$ in such a way that the resulting factorization $y(h):\Spec\dbF\to\scrN$ is as in Definition 1.4.1.}

\medskip
\proof
Let $c:=b-\tilde c$, where $\tilde c$ is the codimension in $\Spec \scrR/p\scrR$ of the point of $\Spec \scrR/p\scrR$ through which $\tilde y_1$ factors. Let $\scrZ:\Spec W(k)[[x_1,\ldots,x_c]]\to\Spec \scrR$ be a dominant morphism that lifts $\tilde y_1$. Let $\bar\scrZ:\Spec k[[x_1,\ldots,x_c]]\to\Spec \scrR$ be induced by $\scrZ$. From properties 2.5 (b) and (c) and Subsection 2.5.2 applied over $k$, we get that the morphism $\bar\scrY(h)\circ\bar\scrZ$ factors uniquely through $\scrN$ in such a way that the resulting morphism $\Spec k[[x_1,\ldots,x_c]]\to\scrN$ lifts $y_1$. This implies that the morphism $\bar\scrY(h)$ factors through $\scrN$ in such a way that $y_1:\Spec k\to\scrN$ is indeed the resulting factorization of $\tilde y_1$. The fact that the resulting factorization $y(h):\Spec\dbF\to\scrN$ is as in Definition 1.4.1, follows via specializations from the very smoothness of $\scrN$ and from the properties enjoyed by $y_1$. The uniqueness part follows from the uniqueness of the factorization $y(h):\Spec\dbF\to\scrN$ (cf. Lemma 2.5.3) and from the uniqueness part of Subsection 2.5 (b).\endproof 

\medskip\noindent
{\bf 2.6.2. Remark.} If either $p>2$ or $p=2$ and the $2$-rank of $y^*(\grA)$ is $0$, then in fact $\scrY(h)$ itself factors uniquely through $\scrN$. To check this it suffices to show that each lift $z(h):\Spec W(\dbF)\to\Spec\scrR\to\scrM$ of $y(h)$ factors uniquely through $\scrN$. But this follows from [Va11, Prop. 2.2.6] once we remark that the Hodge filtration of $h(M)$ associated naturally to the lift $z(h)$ of $y(h)$ determines uniquely $z(h)$ (to be compared with [Va14, Part I, Lem. 3.5.2 (b) and (c)]).
 
\bigskip\noindent
{\bf 2.7. Lemma (the Tate--Hodge trick).} {\it Let $y:\Spec\dbF\to\scrN$. We assume that there exists a finite, discrete valuation ring extension $V$ of $W(\dbF)$ and a lift $\tilde z:\Spec V\to\scrN$ of $y$ such that the Mumford--Tate group of the generic fibre $\scrA$ of the abelian scheme $\tilde z^*(\grA)$ over $V$ is naturally identified with a torus $T$ of $G$ for which the following property holds:

\medskip
{\bf ($\natural$)} the torus $T$ is the subgroup of $\pmb{GSp}(W,\psi)$ that fixes a semisimple $\dbQ$--subalgebra of $\End(W)$ (this semisimple $\dbQ$--subalgebra can be taken to be the centralizer of $T$ in $\End(W)$). 

\medskip
Then the Tate--Hodge property holds for $y:\Spec\dbF\to\scrN$. Therefore the Conjecture 1.2 holds for the point $y:\Spec\dbF\to\scrN$ (i.e., regardless of the choice of the family of tensors $(v_{\alpha})_{\alpha\in\scrJ}$, each $t_{\alpha}$ is the crystalline realization of a Hodge cycle $c_{\alpha,A}$ on $A$).}

\medskip
\proof
It is well known that the Mumford--Tate conjecture holds for abelian varieties over number fields that have complex multiplication. From this and ($\natural$) we get that the Tate conjecture holds for $\scrA$, cf. [Mi6, Thm. 2.6]. From this and [Mi6, Thm. 6.2] we get that the Hodge conjecture holds for $\scrA$. Thus the Hodge--Tate property holds for $y:\Spec\dbF\to\scrN$ and moreover each $z^*(w_{\alpha}^{\grA})$ is the Hodge cycle associated to an algebraic cycle $c_{\alpha,\scrA}$ on $\scrA$. Let $c_{\alpha,A}$ be the algebraic cycle on $A$ which is the specialization of $c_{\alpha,\scrA}$. The crystalline realization of $c_{\alpha,\scrA}$ is $t_{\alpha}$ (cf. the very definition of $t_{\alpha}$) and therefore the crystalline realization of $c_{\alpha,A}$ is as well $t_{\alpha}$.\endproof 

\bigskip\smallskip
\noindent
{\boldsectionfont 3. Stratifications}
\bigskip 

We recall that $k(v)$ is the residue field of $v$. In this section we study three stratifications of $\scrN_{k(v)}$ and present few basic properties of them that are needed in what follows. In particular, we prove a refined general analogue of Oort's result that the supersingular locus of the special fibre of an integral canonical model of a Siegel modular variety is the only closed Newton polygon stratum of the special fibre (see [Oo1] and Example 3.7). We consider an algebraically closed field $k$ of characteristic $p$.

\bigskip\noindent
{\bf 3.1. Rational stratification.} It is known that there exists a $G(\dbA_f^{(p)})$-invariant stratification of $\scrN_{k(v)}$ in reduced, locally closed subschemes defined by the following property:

\medskip\noindent
{\it Two points $y_1,y_2\in\scrN(k)$ belong to the same stratum if and only if there exists an isomorphism
$$(M_1[{1\over p}],\phi_1,(t_{1,\alpha})_{\alpha\in\scrJ})\arrowsim (M_2[{1\over p}],\phi_2,(t_{2,\alpha})_{\alpha\in\scrJ}).$$}

\noindent 
We recall (cf. Subsection 2.1), that the quadruple $(M_i,\phi_i,(t_{i,\alpha})_{\alpha\in\scrJ},\psi_{M_i})$ is the analogue of the quadruple $(M,\phi,(t_{\alpha})_{\alpha\in\scrJ},\psi_{M})$ but obtained for the $k$-valued point $y_i$ instead of the $\dbF$-valued point $y$. 
This type of stratification of $\scrN_{k(v)}$ was first introduced in [RR] and our shorter definition follows [Va10, Subsect. 5.3]; we call it the rational stratification of $\scrN_{k(v)}$. 

By the {\it basic locus} of $\scrN_{k(v)}$ we mean the closed stratum $\grs_0$ of the rational stratification  of $\scrN_{k(v)}$ whose $\dbF$-valued points are the basic points, cf. [Va13, Thm. 8.3 (a)] and Subsection 2.2. We have $y\in\grs_0(\dbF)$ if and only if all Newton polygon slopes of $(\Lie(\scrG_{B(\dbF)}),\phi)$ are $0$.

\medskip\noindent
{\bf 3.1.1. Lemma.} {\it We assume that the adjoint group $G^{\ad}$ is simple. Then the basic locus $\grs_0$ is a pro-\'etale cover of a projective $k(v)$-scheme.}

\medskip
\proof
To prove this, we can assume that $\scrN_{k(v)}$ is not a pro-\'etale cover of a projective $k(v)$-scheme. Thus $(G,X)$ is not compact, cf. [Va14, Part I, Lem. 2.2.4]. Let $H_0$ be a compact, open subgroup of $G(\dbA_f^{(p)})$. For $H_0$ small enough, the pull-back to $\grs_0$ of $\grA$ is defined over $\grs_0/H_0$. Thus also the pull-back  to $\grs_0$ of the $p$-divisible group of $\grA$ is defined over $\grs_0/H_0$. 

For $H_0$ small enough, the resulting $p$-divisible group $\scrW$ over $\grs_0/H_0$ is a direct sum of an \'etale $p$-divisible group $\scrW_0$, of the Cartier dual $\scrW_1$ of $\scrW_0$, and of a $p$-divisible group $\scrW_{(0,1)}$ which does not have integral Newton polygon slopes at any geometric point of $\grs_0/H_0$. More precisely, $\scrW_0$, $\scrW_{(0,1)}$, and $\scrW_1$ correspond to suitable subrepresentations of the representation  on $L_p$ of the $\dbG_m$ subtorus of $Z^0(G_{\dbZ_p})$ which is naturally associated to the Newton quasi-cocharacters of the Shimura $F$-crystals $(M,\phi,\scrG)$ over $\dbF$ attached to basic points $y:\Spec\dbF\to \grs_0\to\grs_0/H_0$. 

Moreover, $\scrW_0$ becomes constant after pull-back to $(\grs_0/H_0)_{\dbF}$ (this is so as the mentioned pull-back of $\scrW_0$ lifts to a constant \'etale $p$-divisible group over $\scrN_{W(\dbF)}$). 

Let $\bar\grs_0/H_0$ be a compactification of $\grs_0/H_0$ coming from a toroidal compactification of $\scrN/H_0$; thus the natural abelian scheme over $\grs_0/H_0$ extends to a semiabelian scheme over $\bar\grs_0/H_0$ whose fibres over geometric points of $(\bar\grs_0/H_0)\setminus  (\grs_0/H_0)$ are not abelian varieties. Let $\Spec \bar V\to\bar\grs_0/H_0$ be a morphism, where $\bar V$ is an $\dbF$-algebra which is a discrete valuation ring. We assume that the generic point $\Spec \bar K$ of $\Spec \bar V$ maps to a point of $\grs_0/H_0$. Then the pull-back to $\Spec K$ of $\scrW$ extends to a $p$-divisible group over $\Spec \bar V$, cf. the very description of $\scrW$ in the last two paragraphs. This implies that the morphism $\Spec \bar V\to\bar\grs_0/H_0$ factors through $\grs_0/H_0$, cf. [dJ, Criterion 2.5]. Thus $\grs_0/H_0$ is a projective $k(v)$-scheme.\endproof

\medskip\noindent
{\bf 3.1.2. Definition.} A point $y:\Spec\dbF\to\scrN$ is called {\it pivotal} (resp. {\it Levi}) if its attached Shimura $F$-crystal is pivotal (resp. Levi) in the sense of Definition A3 (c) (resp. A3 (b)) of Appendix. 

\medskip\noindent
{\bf 3.1.3. Example.} We assume that $\scrN=\scrM$. A point $y:\Spec\dbF\to\scrN$ is basic (resp. pivotal) if and only if $A$ is a supersingular abelian variety (resp. if and only if the $p$-divisible group of $A$ is a product of  $p$-divisible groups of supersingular elliptic curves). 

In general, each basic point $y:\Spec\dbF\to\scrN$ is a Levi point.

\medskip\noindent
{\bf 3.1.4. Lemma.} {\it There exist special points $z\in\scrN(W(\dbF))$ that lift pivotal points $y\in\grs_0(\dbF)$.}

\medskip
\proof
Let $T_{\dbZ_p}$ be a maximal torus of $G_{\dbZ_p}$. We consider the family of Shimura $F$-crystals over $\dbF$ of the form $\{\scrC_g:=(L_p^{\vee}\otimes_{\dbZ_p} W(\dbF),g(1_{L_p^{\vee}}\otimes\sigma)\mu_0({1\over p}),G_{W(\dbF)})|g\in G_{W(\dbF)}(W(\dbF))\}$, where $\mu_0:\dbG_m\to T_{\dbF}$ is a cocharacter whose extension to $\dbC$ via an $O_{(v)}$-monomorphism $W(\dbF)\hookrightarrow\dbC$ is $G(\dbC)$-conjugate to the Hodge cocharacters $\mu_x:\dbG_m\to G_{\dbC}$ that are associated naturally to each $x\in X$. We know that there exists $w\in G_{W(\dbF)}(W(\dbF))$ which normalizes  $T_{W(\dbF)}$ and such that $\scrC_w$ is pivotal in the sense of [Va8, Def. 8.2], cf. [Va8, Cor. 11.1 (c) and Thm. 8.3]. From [Va10, proof of Thm. 5.2.3] we get that there exists a special point $z\in\scrN(W(\dbF))$ that lifts a point $y\in\scrN(\dbF)$ with the property that $(M,\phi,\scrG)$ is isomorphic to $\scrC_w$. Thus $y\in\grs_0(\dbF)\subseteq \scrN(\dbF)$ is a pivotal point.\endproof

\medskip\noindent
{\bf 3.1.5. Definition.} Let $y:\Spec\dbF\to\scrN$ be a Levi point. Let $\scrL$ be the Levi subgroup scheme of $(M,\phi,\scrG)$ in the sense of Subsection A2. We say the {\it weak isogeny property holds} for the Levi point $y:\Spec\dbF\to\scrN$, if for each element $h\in\grI(y)\cap \scrL(B(k))$ the point $y(h):\Spec\dbF\to\scrM$ factors naturally  through a morphism $y(h):\Spec\dbF\to\scrN$ in such a way that for a (any) lift $z(h):\Spec W(\dbF)\to\scrN$ of it, every tensor $t_{\alpha}\in\scrT(\End(M[{1\over p}]))$ is the crystalline realization of the Hodge cycle $z(h)_{B(\dbF)}^*(w_{\alpha}^{\grA})$ on $z(h)^*(\grA)_{B(\dbF)}$.

\medskip
Note that if $y:\Spec\dbF\to\scrN$ is a basic point, then the isogeny property for it is the same as the weak isogeny property for it.

The following exercise is not used in this paper; it is included only to point out a general direct (but more computational) way to check that the isogeny property holds for all points $y:\Spec\dbF\to\scrN$.  

\medskip\noindent
{\bf 3.1.6. Exercise.} The isogeny property holds for all points $y:\Spec\dbF\to\scrN$ if and only if the weak isogeny property holds for all Levi points $y:\Spec\dbF\to\scrN$. Hint: use the non-positive standard form of $(M,\phi,\scrG)$ introduced in [Va10, Subsect. 3.2]. 

\bigskip\noindent
{\bf 3.2. Proposition.} {\it Each stratum $\grs$ of the rational stratification of $\scrN_{k(v)}$ is an open closed subscheme of the Newton polygon stratification of $\scrN_{k(v)}$ defined by the $F$-crystal over $\scrN_{k(v)}$ of the $p$-divisible group of $\grA_{k(v)}$. Moreover, $\grs_0$ is the reduced scheme of the pull-back via the morphism $\scrN_{k(v)}\to\scrN^{\text{cl}}_{k(v)}$ of a reduced, closed subscheme $\grs_0^{\text{cl}}$ of $\scrN^{\text{cl}}_{k(v)}$.}

\medskip
\proof
The first part is a particular case of [Va10, Thm. 5.3.1 (b)]. To prove the second part, let $\grs_0^{\text{cl}}$ be the image of $\grs_0$ in $\scrN^{\text{cl}}$. As the morphism $\scrN\to\scrM_{O_{(v)}}$ is finite, $\grs_0^{\text{cl}}$ is a reduced, closed subscheme of $\scrN^{\text{cl}}_{k(v)}$. Thus to end the proof, it suffices to show that if $y\in\grs_0(\dbF)$ and $y_1\in\scrN(\dbF)$ are two points that map to the same $\dbF$-valued point of $\scrN^{\text{cl}}(\dbF)$, then we have $y_1\in\grs_0(\dbF)$. We can identify $(A,\lambda_A)=(A_1,\lambda_{A_1})$ and $(M,\phi,\psi_M)=(M_1,\phi_1,\psi_{M_1})$ and thus we can view both $\scrG$ and $\scrG_1$ as subgroup schemes of $\pmb{\GL}_M$. Let $(A_0,\lambda_{A_0})$, $\dbF_{p^q}$, $(M_0,\phi_0,(t_{\alpha})_{\alpha\in\scrJ},\psi_{M_0})$, and $\pi_0$ be as in Subsection 2.1. For $q>>0$, we can assume that each $t_{1,\alpha}$ is also a tensor of the tensor algebra of $\End(M_0)[{1\over p}]$. Thus $\scrG$ and $\scrG_1$ are pull-backs of reductive, closed subgroup schemes $\scrG_0$ and $\scrG_{1,0}$ of $\pmb{\GL}_{M_0}$. As $y$ is basic we have $\pi_0\in Z^0(G)(\dbQ)\leqslant Z^0(\scrG_0)(B(\dbF_{p^q}))=Z^0(\scrG_{1,0})(B(\dbF_{p^q}))$ (here we are using the canonical identifications $Z^0(\scrG_{B(\dbF)})=Z^0(G)_{B(\dbF)}=Z^0(\scrG_{1,B(\dbF)})$). This implies that all Newton polygon slopes of $(\Lie(\scrG_{1,0})[{1\over p}],\phi_0)$ are $0$. Thus we have $y_1\in\grs_0(\dbF)$.\endproof

\bigskip\noindent
{\bf 3.3. Lemma.} {\it Let $y, y_1:\Spec\dbF\to\scrN$ be two points that map to the same $\dbF$-valued point of $\scrN^{\text{cl}}$. We assume that the isogeny property holds for both $y$ and $y_1$. Then $\scrE$ and its analogue $\scrE_1$ for $y_1$ have the same identity component. If moreover $y$ is a basic point and the endomorphism property holds for it, then we have $y=y_1$.}

\medskip
\proof
From the very definitions of $\scrE^+$ and of its analogue $\scrE_1^+$ for $y_1$, we can identify $\scrE_1^+=\scrE^+$. From this and Corollary 2.5.4 (b) applied to both $y$ and $y_1$ we get that $\scrE_1$, $\scrE$, and $\scrE_1^+=\scrE^+$ have the same identity components. 

We assume that $y$ is basic and that the endomorphism property holds for it. As $\Lie(\scrG)[{1\over p}]=\Lie(\scrE)\otimes_{\dbQ} B(\dbF)=\Lie(\scrE_1)\otimes_{\dbQ} B(\dbF)\subseteq \Lie(\scrG_1)[{1\over p}]$ (cf. Subsection 2.2 for the first equality), by reasons of dimensions we get that $\Lie(\scrG)[{1\over p}]=\Lie(\scrG_1)[{1\over p}]$. This implies that $\scrG_{B(\dbF)}=\scrG_{1,B(\dbF)}$ and thus $\scrG=\scrG_1$. From this and Lemma 2.5.3 we get that $y=y_1$.\endproof

\bigskip\noindent
{\bf 3.4. Traverso and level $1$ stratifications.} It is known that [Va14, Part I, property 3.5.1 (iii)] implies (cf. [Va6, Sect. 8.7] and [Va7, Subsubsect. 4.2.4, Cor. 4.3, and Ex. 4.6]) that there exists a $G(\dbA_f^{(p)})$-invariant stratification of $\scrN_{k(v)}$ (in the sense of [Va6, Subsect. 8.2]) in regular, equidimensional, locally closed subschemes (of either $\scrN_{k(v)}$ or of pull-backs of $\scrN_{k(v)}$ to algebraically closed fields that contain $k(v)$ and that have countable transcendental degrees) that is defined by the following property:

\medskip\noindent
{\it Two points $y_1,y_2\in\scrN(k)$ factor through the same stratum if and only if there exists an isomorphism $(M_1,\phi_1,(t_{1,\alpha})_{\alpha\in\scrJ})\arrowsim (M_2,\phi_2,(t_{2,\alpha})_{\alpha\in\scrJ})$.}

\medskip
We call such a stratification the Traverso stratification of $\scrN_{k(v)}$.
In what follows we also need the level $1$ form of the Traverso stratification.  It is known (see [Va6, Sect. 8.6] and [Va7, Subsubsect. 4.2.3, Cor. 4.3, and Ex. 4.6]; see also [Va8, Basic Thm. D and Rm. 12.4 (a)]) that there exists a $G(\dbA_f^{(p)})$-invariant stratification of $\scrN_{k(v)}$ in regular, equidimensional, locally closed subschemes that is defined by the following property:

\medskip\noindent
{\it Two points $y_1,y_2\in\scrN(k)$ belong to the same stratum if and only if there exists an isomorphism $(M_1,\phi_1,(t_{1,\alpha})_{\alpha\in\scrJ})\arrowsim (M_2,g_{12}\phi_2,(t_{2,\alpha})_{\alpha\in\scrJ})$ for some $g_{12}\in \Ker(\scrG_2(W(k))$$\to\scrG_2(k))$.}

\medskip
We call such a stratification the level $1$ stratification of $\scrN_{k(v)}$.

\medskip\noindent
{\bf 3.4.1. Lemma.} {\it Each stratum $\grs$ of either the level $1$ or the Traverso stratification of $\scrN_{k(v)}$ is a quasi-affine scheme.}  

\medskip
\proof
The isomorphism class of the principally quasi-polarized Barsotti--Tate group of level $1$ associated to a point $y_i:\Spec k\to\grs$ (i.e., of the reduction modulo $p$ of $(M_i,\phi_i,p\phi_i^{-1},\psi_i)$), does not depend on the choice of the $k$-valued point $y_i$ of $\grs$. This implies that the image of $\grs$ in $\scrN_{k(v)}^{\text{cl}}$ is a reduced, locally closed subscheme of a stratum of the Ekedahl--Oort (i.e., of level $1$) stratification of $\scrM_{k(v)}$ and thus (cf. [Oo2, Thm. 1.2]) it is a quasi-affine scheme. As the morphism $\scrN\to\scrN^{\text{cl}}$ is finite, we conclude that $\grs$ is a quasi--affine scheme.\endproof

\bigskip\noindent
{\bf 3.5. Lemma.} {\it Let $\grz$ be the schematic closure in $\scrN_{k(v)}$ of the $G(\dbA_f^{(p)})$-orbit of some point $y:\Spec\dbF\to\scrN$. We assume that $\grz$ is $0$ dimensional (for instance, this holds if $y$ is a pivotal point, cf. [Va8, Cor. 11.1 (c) and Basic Thm. D]).

\medskip
{\bf (a)} Then $y$ is a basic point.

\smallskip
{\bf (b)} If the isogeny property holds for $y$, then the endomorphism property holds for $y$.} 

\medskip
\proof
Let $H_0$ be a compact, open subgroup of $G(\dbA_f^{(p)})$. As $\grz$ is $0$ dimensional, the map
$$[\pmb{Aut}(A)(\dbZ_{(p)})\cap\scrE^+(\dbQ)]\backslash [\{M\}\times G(\dbA_f^{(p)})/H_0]\to o^{\text{big}}(y)/H_0$$
has a finite domain. This implies that the set $\scrE^+(\dbQ)\backslash G(\dbA_f^{(p)})/H_0$ is finite. 

We will show that the assumption that $y$ is not a basic point leads to a contradiction. Let $A_0$, $\pi_0$, and $\dbF_{p^q}$ be as in Subsection 2.1. As $y$ is not basic, the crystalline realization of $\pi_0$ is an element $\pi_0\in \scrG(B(\dbF))$ which does not belong to $Z(\scrG)(B(\dbF))$. As the elements of $Z(G)(\dbQ)$ define naturally $\dbQ$--endomorphisms of the universal abelian scheme $\grA$ over $\scrN$, we get that if $l$ is a prime distinct from $p$, then the $l$-adic realization of $\pi_0$ is an element of $G(\dbQ_l)$ that does not belong to $Z(G)(\dbQ_l)$. Therefore $\scrE^+_{\dbQ_l}$ is contained in a reductive subgroup $\scrF$ of $G_{\dbQ_l}$ distinct from $G_{\dbQ_l}$ (we can take $\scrF$ to be the centralizer in $G_{\dbQ_l}$ of the mentioned $l$-adic realization). We now choose $l$ such that both groups $G_{\dbQ_l}$ and $\scrF$ are split (equivalently, that both $G_{\dbQ_l}$ and the Frobenius torus of $\pi_0$ are split). As the quotient variety $Q:=\scrF\backslash G_{\dbQ_l}$ is a smooth, connected, affine variety of positive dimension, it is easy to see that the analytic $l$-adic subvariety $\scrF(\dbQ_l)\backslash G(\dbQ_l)$ of $Q(\dbQ_l)$ is not compact. But as the set $\scrE^+(\dbQ)\backslash G(\dbA_f^{(p)})/H_0$ is finite and as $\scrE^+(\dbQ_l)\leqslant \scrF(\dbQ_l)$, one gets that the analytic $l$-adic subvariety $\scrF(\dbQ_l)\backslash G(\dbQ_l)$ is compact. Contradiction. Thus $y$ is a basic point.

To check (b), we first remark that as in the previous paragraph we argue that if $\scrE^+$ is not a $\dbQ$--form of $\scrG_{B(\dbF)}$ (equivalently, of $G$), then we similarly reach a contradiction. From this and Corollary 2.5.4 (b) we get that $\scrE$ is a $\dbQ$--form of $\scrG_{B(\dbF)}$ and thus the endomorphism property holds for $y$ (cf. Subsection 2.2).\endproof

\bigskip\noindent
{\bf 3.6. Lemma.} {\it We assume that the Shimura pair $(G,X)$ is compact (in the sense before Theorem 1.5).}${}^3$ $\vfootnote{3}{The proof of Lemma 3.6 is the only place in the paper where we require to work with compact Shimura varieties. Therefore, all the results of the paper continue to hold for all those connected components of strata of the Traverso stratification of $\scrN_{k(v)}$ which are locally closed subschemes of $\scrN_{\dbF}$ and whose schematic closures in $\scrN_{\dbF}$ contain points of $\grs_0$ (like for all those strata that involve abelian varieties of $p$-rank $0$); one can check that this applies to all cases of non-compact Shimura pairs. Also, all the results of the paper continue to hold provided they do not rely on Lemma 3.6.}$ {\it Let $\grz$ be a closed, reduced, $G(\dbA_f^{(p)})$-invariant subscheme of $\scrN_{k(v)}$. Then the intersection $\grz\cap\grs_0$ is non-empty.} 

\medskip
\proof
We consider the stratification of $\grz$ which is the pull-back of the Traverso stratification of $\scrN_{k(v)}$. Let $\grz_0$ be a stratum of this stratification of $\grz$ which is a closed, reduced subscheme of $\grz$. The scheme $\grz_0$ is $G(\dbA_f^{(p)})$-invariant and therefore (cf. Lemma 3.5 (a)) it has either a non-empty intersection with the basic locus or positive dimension. It is known that $\scrN$ is a pro-\'etale cover of a projective, smooth $O_{(v)}$-scheme, cf. [Va4, Cor. 4.3]. Thus $\grz_0$ is a pro-\'etale cover of a projective, smooth $k(v)$-scheme (as this holds for $\scrN_{k(v)}$) and it is a quasi-affine scheme (cf. Lemma 3.4.1). Therefore $\grz_0$ can not be of positive dimension. We conclude that $\grz\cap\grs_0\supseteq\grz_0\cap\grs_0\neq\emptyset$.\endproof

\bigskip\noindent
{\bf 3.7. Example.} Let $y:\Spec\dbF\to\scrN$ be a point. Let $\grl$ be the stratum of the Traverso stratification of $\scrN_{k(v)}$ that contains $y$. Let $\grc$ be the connected component of $\grl$ that contains $y$. Let $\grt$ be the smallest closed, reduced subscheme of $\grl$ that contains $\grc$ and is $G(\dbA_f^{(p)})$-invariant; thus $\grc$ is a connected component of $\grt$. If $y$ is not a basic point, then the dimension of $\grt$ is positive and this implies that $\grc$ has positive dimension.

We assume now that the Shimura pair $(G,X)$ is compact. The schematic closure of $\grt$ in $\scrN_{k(v)}$ contains points of $\grs_0$, cf. Lemma 3.6. From this and the very definition of $\grt$, we get that the schematic closure of $\grc$ in $\scrN_{k(v)}$ contains also points of $\grs_0$.

\bigskip\noindent
{\bf 3.8. Theorem.} {\it Let $y:\Spec\dbF\to\scrN$ be a point. Let $\grc$ be the connected component of the Traverso stratum which is a locally closed subscheme of $\scrN_{k(v)}$ and which contains the point $y$. Let $\grr$ (resp. $\grl$) be the connected component of the rational (resp. level $1$) stratum of $\scrN_{k(v)}$ which contains the point $y$. 

\medskip
{\bf (a)} If $y$ is not a basic point, then we have $\dim(\grc)>0$.

\smallskip
{\bf (b)} If $y$ is not a basic point and $\grl$ has a non-trivial intersection with $\grs_0$, then we have $\dim_y(\grr\cap\grl)\ge 2$ (here by $\dim_y$ we mean the dimension at the point $y$).

\smallskip
{\bf (c)} If $y$ is a basic point which is not pivotal, then we have $\dim_y(\grr\cap\grl)\ge 1$.}

\medskip
\proof 
Part (a) follows from either Example 3.7 or Thorem A6 (a) of Appendix (see also below for a third proof). Also (b) and (c) are properties of Shimura $F$-crystals and they are a particular case of Theorem A6 (b) and (c) of Appendix.\endproof

\bigskip\smallskip
\noindent
{\boldsectionfont 4. Reduction steps for the isogeny property}
\bigskip 

In this section we get reduction steps towards the proof of the isogeny property.

\bigskip\noindent
{\bf 4.1. Theorem.} {\it Let $y:\Spec\dbF\to\scrN$ be a point. Let $\grl$ be the stratum of the Traverso stratification of $\scrN_{k(v)}$ that contains $y$. Let $\grc$ be the connected component of $\grl$ that contains $y$. We assume that the following two conditions hold:

\medskip
{\bf (i)} the isogeny property holds for each basic point $\Spec\dbF\to\scrN$;
\smallskip
{\bf (ii)} the schematic closure of $\grc$ in $\scrN_{k(v)}$ intersects non-trivially the basic locus $\grs_0$.

\medskip
Then the isogeny property holds for the point $y:\Spec\dbF\to\scrN$.}

\medskip
\proof
We consider a morphism $y_{\ell}:\Spec \dbF[[x]]\to\scrN$ which lifts $y$ and which has the property that the resulting morphism $\Spec \dbF((x))\to\scrN$ factors through the generic point $y_g$ of $\grc$ (we recall that $\grl$ is regular). We know that $y_g$ specializes to a point $y_1\in\grs_0(\dbF)$, cf. condition (ii).  Let $\scrQ$ be an integral, normal, noetherian, local ring  of residue field $\dbF$ such that we have a morphism $w:\Spec \scrQ\to \scrN$ with the properties that the closed point of $\Spec\scrQ$ maps to the point $y_1$ of $\grs_0$ and the generic point of $\Spec\scrQ$ maps to the point $y_g$ of $\grc$. Let $\grl^{\text{big}}$ be the reduced, locally closed subscheme of $\scrM_{\dbF}$ which is a stratum of the Traverso stratification of $\scrM_{\dbF_p}$ and which contains the composite point $y:\Spec\dbF\to\scrN_{\dbF}\to\scrM_{\dbF}$. It is known  (cf. [Va2, Subsubsect. 5.3.2 (c)]) that there exists an infinite sequence 
$$\cdots\to\grl^{\text{big}}_s\to\cdots\to\grl^{\text{big}}_3\to\grl^{\text{big}}_2\to\grl^{\text{big}}_1\to\grl^{\text{big}}_0=\grl^{\text{big}}$$
of faithfully flat, finite morphisms between regular schemes with the property that:

\medskip
{\bf (iii)} the pull-back $(\scrD_s,\lambda_{\scrD_s})$ to $\grl^{\text{big}}_s$ of the principally quasi-polarized Barsotti--Tate group of level $p^s$ of the universal principally polarized abelian scheme over $\scrM$ is constant i.e., there exists an isomorphism $\nu^{\text{big}}_s:(A[p^s],\lambda_A[p^s])_{\grl^{\text{big}}_s}\arrowsim (\scrD_s,\lambda_{\scrD_s})$.

\medskip
We can assume that:  

\medskip
{\bf (iv)} for each $s\in\dbN^*$ we have $\nu^{\text{big}}_{s,{\grl^{\text{big}}_{s+1}}}=\nu^{\text{big}}_{s+1}[p^s]$.

\medskip
Let $\scrR_1:=W(\dbF)[[x]]$ and let $\Phi_1$ be the Frobenius endomorphism of $\scrR_1$ that is compatible with $\sigma$ and that takes $x$ to $x^p$. Let $R$ be the normalization of $\dbF[[x]]=\scrR_1/p\scrR_1$ in an algebraic closure $k$ of $\dbF((x))$; it is an integral, normal, perfect $\dbF$-algebra. 
Let $(\scrI,\lambda_{\scrI})$ be the principally quasi-polarized $p$-divisible group over $\dbF[[x]]$ of $y_{\ell}^*(\grA,\lambda_{\grA})$. Let 
$$(M_{\ell},\phi_{M_{\ell}},(t_{\ell,\alpha})_{\alpha\in\scrJ},\nabla_{M_{\ell}},\psi_{M_{\ell}})$$ be the principally quasi-polarized $F$-crystal with tensors over $\dbF[[x]]$ associated naturally to the morphism $y_{\ell}:\Spec \dbF[[x]]\to\scrN$. Thus $M_{\ell}$ is a free $\scrR_1$-module of rank $2d$, $\nabla_{M_{\ell}}$ is an integrable, topologically nilpotent connection on $M_{\ell}$, and each $t_{\ell,\alpha}$ is a tensor of the tensor algebra of $\End(M_{\ell})[{1\over p}]$ that is fixed by $\phi_{M_{\ell}}$, that is annihilated by $\nabla_{\ell}$,  and that is the crystalline realization of  the Hodge cycle $z_{\ell,E(G,X)}^*(w_{\alpha}^\grA)$ on $z_{\ell,E(G,X)}^*(\grA)$, where $z_{\ell}:\Spec\scrR_1\to\scrN$ is an arbitrary morphism that lifts $y_{\ell}$. 

Due to properties (iii) and (iv), there exists an isomorphism 
$$\iota:(A[p^{\infty}],\lambda_A[p^{\infty}])_R\arrowsim (\scrI,\lambda_{\scrI})_R.$$ 
We can assume that this isomorphism $\iota$ lifts the identity automorphism of $(A[p^{\infty}],\lambda_A[p^{\infty}])$. The crystalline realization of $\iota$ is an isomorphism (denoted in the same way)
$$\iota:(M_{\ell},\phi_{M_{\ell}},(t_{\ell,\alpha})_{\alpha\in\scrJ},\nabla_{M_{\ell}},\psi_{M_{\ell}})_R\arrowsim (M\otimes_{W(\dbF)} \scrR_1,\phi\otimes\Phi_1,(t_{\alpha})_{\alpha\in\scrJ},\delta,\psi_M)_R$$
of principally quasi-polarized $F$-crystals with tensors over $R$ (here $\delta$ is the flat connection on $M\otimes_{W(\dbF)} \scrR_1$ that annihilates $M\otimes 1$ and the lower left index $R$ means a pull-back to an $F$-crystal over $R/pR$). The fact that $\iota$ maps each $t_{\ell,\alpha}$ to $t_{\alpha}$ follows from the fact that there exists no element of the tensor algebra of $\End(M)\otimes_{W(\dbF)} \Ker(W(R)\twoheadrightarrow R)[{1\over p}]$ which is fixed by $\phi\otimes \Phi_R$, where $\Phi_R$ is the Frobenius endomorphism of $W(R)[{1\over p}]$.  

Let $h\in\grI(y)$. By composing $\iota$ with the natural $\dbZ[{1\over p}]$-isogeny 
$$(M\otimes_{W(\dbF)} \scrR_1,\phi\otimes\Phi_1,(t_{\alpha})_{\alpha\in\scrJ},\delta,\psi_M)_R \rightsquigarrow (h(M)\otimes_{W(\dbF)} \scrR_1,\phi\otimes\Phi_1,(t_{\alpha})_{\alpha\in\scrJ},\delta,\psi_M)_R,$$ 
we get a $\dbZ[{1\over p}]$-isogeny
$$\iota(h):(M_{\ell},\phi_{M_{\ell}},(t_{\ell,\alpha})_{\alpha\in\scrJ},\nabla_{M_{\ell}},\psi_{M_{\ell}})_R \rightsquigarrow (h(M)\otimes_{W(\dbF)} \scrR_1,\phi\otimes\Phi_1,(t_{\alpha})_{\alpha\in\scrJ},\delta,\psi)_R$$
which is induced naturally by an analogous isogeny between principally quasi-polarized $p$-divisible groups over $R$.
Let $i_{\ell}^{\text{large}}:\Spec R\to\scrM$ be the morphism defined naturally by $y_{\ell}$ and let $j_{\ell}^{\text{large}}(h):\Spec R\to\scrM$  be the morphism obtained naturally from $i_{\ell}^{\text{large}}$ via the $\dbZ[{1\over p}]$-isogeny $\iota(h)$ (in the same way we constructed $y(h):\Spec\dbF\to\scrM$ from $y:\Spec\dbF\to\scrN$). Such a $\dbZ[{1\over p}]$-isogeny is defined over a discrete valuation ring which is a finite, normal $\dbF[[x]]$-subalgebra of $R$. As such a discrete valuation ring is isomorphic to $\dbF[[x]]$, not to introduce extra notations we can assume that $j_{\ell}^{\text{large}}(h):\Spec R\to\scrM$ factors naturally through a morphism $j_{\ell}(h):\Spec \dbF[[x]]\to\scrM$. The morphism $j_{\ell}(h)$ lifts $y(h)$, cf. constructions (here it is critical that $(\scrI,\lambda_{\scrI})$ is geometrically constant and not only of constant Newton polygon). From Subsection 2.5 (b) and Fact 2.6.1 we get that $y(h)$ factors through $\scrN$ as in Definition 1.4.1 if and only the point $j_{\ell,g}(h):\Spec k\to\scrM$ induced naturally by $j_{\ell}(h)$ factors through $\scrN$ in a way analogous to Definition 1.4.1. But we can choose (enlarge) $\scrQ$ such that the morphism $j_{\ell,g}(h):\Spec k\to\scrM$ is the composite of the natural morphism $\Spec k\to\Spec \scrQ$ with a morphism $w(h):\Spec \scrQ\to\scrM$ obtained from $w$ via a natural $\dbZ[{1\over p}]$-isogeny. As above, from Subsection 2.5 (b) and Fact 2.6.1 we get that $w(h)$ factors through $\scrN$ in a way analogous to Definition 1.4.1 if and only if the point $y_1(h):\Spec\dbF\to\scrM$ induced naturally by $w(h)$ factors through $\scrN$ as in Definition 1.4.1. 

We have $y_1(h)=y_1(h_1)$ for a suitable element $h_1\in\grI^{\text{big}}(y_1)$, where $\grI^{\text{big}}(y_1)$ is as in Subsubsection 2.3.1. We assume that we can choose $h_1$ to be an element of $\grI(y_1)$. From this assumption and the condition (i), we get that $y_1(h)$ factors through $\scrN$ as in Definition 1.4.1. We conclude that $y(h)$ factors through $\scrN$ as in Definition 1.4.1. As $h\in\grI(y)$ was arbitrary, we get that the isogeny property holds for $y$.

Thus to end the proof of the theorem, it suffices to show that we can choose $h_1$ to be an element of $\grI(y_1)$. The proof of this is essentially the same as of [Va14, Part I, Subsubsects. 5.1 to 5.3]. We recall the details. We consider a morphism $y_{s}:\Spec \dbF[[x]]\to\Spec \scrQ$ which is dominant and for which the composite morphism $w(h)\circ y_s:\Spec \dbF[[x]]\to\Spec \scrM$ lifts $y_1(h)$. Let $(M_s,\phi_{M_s},(t_{s,\alpha})_{\alpha\in\scrJ},\nabla_s,\psi_{M_s})$ be the principally quasi-polarized $F$-crystal with tensors over $\dbF[[x]]$ associated naturally to $w(h)\circ y_s:\Spec \dbF[[x]]\to\Spec \scrM$ (the part on tensors makes sense here due to de Jong's extension theorem, to be compared with [Va14, Part I, Subsect. 5.1]).

Let $\scrK_1$ be the field of fractions of $\scrR_1$ and let $\scrO_1$ be the local ring of $\scrR_1$ which is a discrete valuation ring of mixed characteristic $(0,p)$. We have $\scrO_1/p\scrO_1=\dbF((x))$. Let $\scrG_s$ be the schematic closure in $\pmb{\GL}_{M_s}$ of the subgroup of $\pmb{\GL}_{M_s\otimes_{\scrR_1} \scrK_1}$ that fixes each $t_{s,\alpha}$ with $\alpha\in\scrJ$. If $\scrR_1\to\scrO_1\to W(k)$ is the Teichm\"uller lift that lifts the natural embedding $\dbF[[x]]\hookrightarrow k$, then there exists an isomorphism 
$$(M_s,(t_{s,\alpha})_{\alpha\in\scrJ},\psi_{M_s})\otimes_{\scrR_1} W(k)\arrowsim (h(M),(t_{\alpha})_{\alpha\in\scrJ},\psi_{M})\otimes_{W(\dbF)} W(k),\leqno (6)$$ cf. constructions. This implies that $\scrG_{s,\scrO_1}$ is a reductive subgroup scheme of $\pmb{\GL}_{M_s\otimes_{\scrR_1} \scrO_1}$.

From this, as in [Va14, Part I, Thm. 5.2] we argue that $\scrG_s$ itself is a reductive subgroup scheme of $\pmb{\GL}_{M_s}$. As in [Va14, Part I, Subsect. 5.3] we argue that we have a cocharacter $\mu_s:\dbG_m\to\scrG_s$ such that there exists a direct sum decomposition $M_s=F^1_s\oplus F^0_s$ with the properties that for $i\in\{0,1\}$ the group scheme $\dbG_m$ acts on $F^i_s$ through $\mu_s$ via the weight $-i$ and that $F^1_s/pF^1_s$ is the kernel of $\phi_s$ modulo $p$. 

We consider the morphism $\scrR_1\to W(\dbF)$ that maps $x$ to $0$. We have
$$(h(M_1),\phi_1,(t_{1,\alpha})_{\alpha\in\scrJ},\psi_{h(M_1)})=(M_s,\phi_s,(t_{s,\alpha})_{\alpha\in\scrJ},\psi_{M_s})\otimes_{\scrR_1} W(\dbF)\leqno (7)$$
(cf. the fact that $y_s$ specializes to $y_1(h)$). Let $F^1_1:=F^1_s\otimes_{W(k)} W(\dbF)\subset h(M_1)$.

As in Subsection 2.5 (a) and (b), Faltings deformation theory implies that the quintuple $(M_s,\phi_s,F^1_s,(t_{s,\alpha})_{\alpha\in\scrJ},\psi_{M_s})$ is induced from a deformation of $(h(M_1),F^1_1,\phi_1,(t_{1,\alpha})_{\alpha\in\scrJ},\psi_{h(M_1)})$ to $\scrR_b/p\scrR_b$ constructed as in Subsection 2.5 (a) (see [Va14, Part I, Appendix, Thm. B6.4]). In particular, we have an isomorphism
$$(h(M_1),(t_{1,\alpha})_{\alpha\in\scrJ},\psi_{h(M_1)})\otimes_{W(\dbF)} W(\dbF)[[x]]\arrowsim (M_s,(t_{s,\alpha})_{\alpha\in\scrJ},\psi_{M_s}).\leqno (8)$$
From (6) to (8) and (1) we get that there exists an isomorphism $(M_1,(t_{1,\alpha})_{\alpha\in\scrJ},\psi_{M_1})\otimes_{W(\dbF)} W(k)\arrowsim (h(M_1),(t_{1,\alpha})_{\alpha\in\scrJ},\psi_{h(M_1)})\otimes_{W(\dbF)} W(k)$. From this and [Va2, Part I, Lem. B4] we get that there exist isomorphisms
$$(M_1,(t_{1,\alpha})_{\alpha\in\scrJ},\psi_{M_1})\arrowsim (h(M_1),(t_{1,\alpha})_{\alpha\in\scrJ},\psi_{h(M_1)}).$$
This implies that we can assume that $h\in\scrG_1(B(\dbF))$. The fact that we have $h\in\grI(y_1)$ follows from the existence of the cocharacter $\mu_s$.  \endproof

\bigskip\noindent
{\bf 4.2. Proposition.} {\it We assume that the Conjecture 1.2 holds for the point $y:\Spec\dbF\to\scrN$. Then the endomorphism property holds for $y:\Spec\dbF\to\scrN$.}

\medskip
\proof
We consider the direct sum decomposition
$$\End(W)=\Lie(G)\oplus\Lie(G)^\perp,\leqno (9)$$
where $\Lie(G)^\perp$ is the perpendicular on $\Lie(G)$ with respect to the trace form on $\End(W)$ (cf. [Va14, Part I, Appendix, Lem. A2 (b)]). Let $\Pi$ be the projector of $\End(W)$ on $\Lie(G)$ along  $\Lie(G)^\perp$; we identify it as well with an element of $\End(W)\otimes_{\dbQ} \End(W)^{\vee}=\End(W)\otimes_{\dbQ} \End(W)\subseteq\scrT(\End(W))$ fixed by $G$. Thus we can assume that there exists an element $\alpha_0\in\scrJ$ such that we have $v_{\alpha_0}=\Pi$. 

Let $A_0$, $\dbF_{p^q}$, $(M_0,\phi_0,(t_{\alpha})_{\alpha\in\scrJ},\psi_{M_0})$, and $\pi_0$ be as in Subsection 2.1. The centralizer $\scrC$ of $\pi_0$ in $\End(M[{1\over p}])$ is the $B(\dbF)$-subalgebra generated by crystalline realizations of $\End(A)^{\text{opp}}\otimes_{\dbZ} \dbQ$. Therefore we can identify naturally
$$\scrC=\End(A)^{\text{opp}}\otimes_{\dbZ} B(\dbF).\leqno (10)$$
As $\pi_0$ fixes $t_{\alpha_0}$ and due to (9),
we have a direct sum decomposition of $B(\dbF)$-vector spaces
$$\scrC=(\scrC\cap \Lie(\scrG_{B(\dbF)}))\oplus (\scrC\cap \Lie(\scrG_{B(\dbF)})^\perp),\leqno (11)$$
where $\Lie(\scrG_{B(\dbF)})^\perp$ is the perpendicular on $\Lie(\scrG_{B(\dbF)})$ with respect to the trace form on $\End(M[{1\over p}])$. But $t_{\alpha_0}$ is the crystalline realization of an algebraic cycle $c_{\alpha_0,A}$ on $A$, cf. the assumption that the Conjecture 1.2 holds for $y$. From this, (10), and (11), we get that we have a direct sum decomposition of $\dbQ$--vector spaces
$$\End(A)^{\text{opp}}\otimes_{\dbZ} \dbQ=[(\End(A)^{\text{opp}}\otimes_{\dbZ} \dbQ)\cap \Lie(\scrG_{B(\dbF)})]\oplus [(\End(A)^{\text{opp}}\otimes_{\dbZ} \dbQ)\cap \Lie(\scrG_{B(\dbF)})^\perp].$$
In other words, there exists a Lie subalgebra $(\End(A)^{\text{opp}}\otimes_{\dbZ} \dbQ)\cap \Lie(\scrG_{B(\dbF)})$ of $\End(A)^{\text{opp}}\otimes_{\dbZ} \dbQ$ such that a $\dbQ$--basis for it is also a $\dbQ_p$-basis for $\{x\in\Lie(\scrG_{B(\dbF)})|\phi(x)=x\}\cap (\End(A)^{\text{opp}}\otimes_{\dbZ} \dbQ_p)$ and thus a $B(\dbF)$-basis for $\scrC\cap\Lie(\scrG_{B(\dbF)})$. This implies that there exists a reductive subgroup $\scrE^{\c}$ of the group of invertible elements of $\End(A)^{\text{opp}}\otimes_{\dbZ} \dbQ$ whose Lie algebra is $(\End(A)^{\text{opp}}\otimes_{\dbZ} \dbQ)\cap \Lie(\scrG_{B(\dbF)})$, cf. [Va14, Part I, Appendix, Lem. A2 (a)] applied to the pair of fields $(\dbQ,B(\dbF))$. Due to the fact that the cycle $\Pi$ is algebraic, for each prime $l\neq p$ the $l$-adic realization of an arbitrary element of $\scrE^{\c}(\dbQ)$ is an element of $G(\dbQ_l)$. Thus $\scrE^{\c}$ is a subgroup of the group $\scrE$ of Definition 1.4.3. By reasons of dimensions we get that $\scrE^{\c}$ is the identity component of $\scrE$ and that $\Lie(\scrE)\otimes_{\dbQ} \dbQ_p=\{x\in\Lie(\scrG_{B(\dbF)})|\phi(x)=x\}\cap (\End(A)^{\text{opp}}\otimes_{\dbZ} \dbQ_p)$. Thus the endomorphism property holds for $y:\Spec\dbF\to\scrN$.\endproof   

\bigskip\noindent
{\bf 4.3. Proposition.} {\it Let $y:\Spec\dbF\to\scrN$ be a basic point for which the following six properties hold:

\medskip
{\bf (i)} the derived group $G^{\der}$ is simply connected;

\smallskip
{\bf (ii)} the isogeny property holds for $y$;

\smallskip
{\bf (iii)} the endomorphism property holds for $y$;

\smallskip
{\bf (iv)} there exists a semisimple $\dbZ_{(p)}$-subalgebra of $\End(L_{(p)})$ such that its centralizer in $\pmb{GSp}(L_{(p)},\psi)$ is a reductive group scheme $G_{1,\dbZ_{(p)}}$ with the properties that $Z^0(G)=Z^0(G_1)$ with $G_1$ as the generic fibre of $G_{1,\dbZ_{(p)}}$ and that the monomorphism $G^{\der}_{W(\dbF)}\to G^{\der}_{1,W(\dbF)}$ is either an isomorphism or a product of monomorphisms of one of the following three forms: $\pmb{Sp}_{2n}\hookrightarrow \pmb{SL}_{2n}$, $\pmb{Spin}_{2n-1}\hookrightarrow \pmb{SL}_{2^{n-1}}$, $\pmb{Spin}_{2n}\hookrightarrow \pmb{SL}_{2^{n-2}}\times \pmb{SL}_{2^{n-2}}$;

\smallskip
{\bf (v)} the Hasse principle holds for the tori $G^0/G^{\der}$ and  $G^0_1/G^{\der}_1$, where $G^0=G\cap\pmb{Sp}(W,\psi)$ and $G^0_1:=G_1\cap\pmb{Sp}(W,\psi)$;

\smallskip
{\bf (vi)} either the abelian variety $A$ is supersingular or we have $G=G_1$.

\medskip
Then we have $\grs_0(\dbF)=o(y)$ (i.e., the basic locus consists of one isogeny set).}

\medskip
\proof
Let $G_1$ be the generic fibre of $G_{1,\dbZ_{(p)}}$. We have injective maps $(G,X)\hookrightarrow (G_1,X_1)\hookrightarrow (\pmb{GSp}(W,\psi),S)$ of Shimura pairs, cf. property (iv). Let $y_{'}:\Spec\dbF\to\scrN$ be a basic point. To prove the proposition it suffices to show that $y_{'}\in o(y)$. To prove this we can replace $y_{'}$ by a point in the $G(\dbA_f^{(p)})$-orbit of $y_{'}$. Thus we can assume that both points $y,y_{'}:\Spec\dbF\to\scrN$ factor through the same connected component of $\scrN$, cf. [Va1, Lem. 3.3]. For the sake of clarity, we will divide the proof into five steps as follows.

\smallskip
{\bf Step 1.} 
We know that there exist elements $h\in \grI(y)$ and $t\in G_1(\dbA_f^{(p)})$ such that we have an identity $y_{'}=y(h)t$ of $\dbF$-valued points of $\scrM$, cf. [Va13, Thm. 8.3 (b)] (here we are using properties (i) and (iv); property (v) for the case of $G_1^0$ follows in fact from the property (iv) and it is inserted here only to emphasize that it did play a role in the proof of loc. cit.). As both points $y,y_{'}:\Spec\dbF\to\scrN$ factor through the same connected component of $\scrN$, we have $t\in G_1^0(\dbA_f^{(p)})$ (cf. also the proof of Loc. cit.). If $G=G_1$, then (due to property (iv)) it is well known that $\scrN=\scrN^{\text{cl}}$ and thus $y_{'}\in o(y)$. 

Thus from now on we can assume that $G\neq G_1$; therefore $A$ is a supersingular abelian variety (cf. property (vi)). Let $t_{'}=t^{-1}$. 
Due to the property (ii), we can assume that $h(M)=M$. 
Thus we have $y=y_{'}t_{'}$ and we can identify $(M_{'},\phi_{'},\psi_{M_{'}})=(M,\phi,\psi_M)$ via a canonical $\dbZ_{(p)}$-isogeny $(A_{'},\lambda_{A_{'}})\rightsquigarrow(A,\lambda_A)$. 

Let $b=\dim(X)$ and $\scrR_c$ with $c\in\dbN$ be as in Subsection 2.5. Let $\scrD_y$ be the formal deformation space of $(A,\lambda_A)$ corresponding to $y:\Spec\dbF\to\scrN$ and let $\scrD_y^{\text{big}}$ be the formal deformation space of $(A,\lambda_A)$ corresponding to the composite point $y:\Spec\dbF\to\scrN\to\scrM$. We fix identifications $\scrD_y=\Spf \scrR_b$ and $\scrD_y^{\text{big}}=\Spf \scrR_{{d(d+1)}\over 2}$. Similarly we define $\scrD_{y_{'}}$ and its translation $\scrD_{y_{'},t_{'}}$ via $t_{'}$.  Both $\scrD_y$ and $\scrD_{y_{'},t_{'}}$ are formal deformation subspaces of $\scrD_y^{\text{big}}$ (the first one is defined by a $W(\dbF)$-epimorphism $\scrR_{{d(d+1)}\over 2}\twoheadrightarrow\scrR_b$). For $n\in\dbN^*$ let $\scrD_{y,n}:=\Spf \scrR_b/(p,x_1,\ldots,x_b)^n$. Similarly we define $\scrD_{y,n}^{\text{big}}$ and $\scrD_{y_{'},t_{'},n}$.

Let $\scrE^0$ be the subgroup of $\scrE$ that fixes $\lambda_A$ (equivalently $\psi_M$). Let $\scrE_1$ and $\scrE_1^0$ be the analogues of $\scrE$ and $\scrE^0$ (respectively) but associated to the point $y_{'}:\Spec\dbF\to\scrN\to\scrN_1$ defined by $y$. Here $\scrN_1$ is the integral canonical model of $(G_1,X_1,G_{1,\dbZ_{(p)}}(\dbZ_p),v_1)$ (with $v_1$ the prime of $E(G_1,X_1)$ divided by $v$) and $\scrN\to\scrN_1$ is the natural factorization of $\scrN\to\scrM$. 

Let $g\in\scrE_1^0(\dbQ_p)$ be such that it takes $t_{',\alpha}$ to $t_{\alpha}$ for all $\alpha\in\scrJ$ and we have $g(M)=M$, cf. proof of [Va13, Thm. 8.3 (b)]. We have $g\scrG_{'}g^{-1}=\scrG$. For each $n\in\dbN^*$ there exists a $\dbZ_{(p)}$-automorphism $g_n$ of $(A,\lambda_A)$ such that $g$ and $g_n$ are congruent modulo $p^n$. Let $t_n$ be the image of $g_n$ in $G_1^0(\dbA_f^{(p)})$.

\smallskip
{\bf Step 2.} We will show that the two formal subschemes of $\scrD_{y}^{\text{big}}$ corresponding to the basic (supersinsingular) loci of $\scrD_{y_{'},t_{'}}$ and $\scrD_y$ are permuted under the natural action of a $\dbZ_{(p)}$-automorphism $g_0$ of $(A,\lambda_A)$ on $\scrD_{y}^{\text{big}}$. Let $n_0\in\dbN^*$ be such that each $\dbZ_{(p)}$-automorphism of $(A,\lambda_A)$ acts trivially on the supersingular locus of $\scrD_y^{\text{big}}$. Briefly speaking, $n_0$ is the largest number such that we have ${1\over {p^{n_0}}}\dbZ$--isogenies between suitable extensions of $(A,\lambda_A)$ and the pull-backs of the universal abelian scheme over $\scrD_y^{\text{big}}$ via arbitrary generic points of the supersingular locus of $\scrD_y^{\text{big}}$. Let $g_0:=g_{n_0}$. Let $n\in\dbN^*$ with $n\ge n_0$. By replacing $t_{'}$ with $t_nt_{'}$ and each $t_{',\alpha}$ with $g_n(t_{',\alpha})$, we can assume that $g_n$ is the identity automorphism. Such a replacement is allowed as $g_n$ and $g_0$ act in the same way on the supersingular locus of $\scrD_{y}^{\text{big}}$. As $g_n$ is the identity element and as $g$ is congruent to it modulo $p^n$, the reductive group schemes $\scrG_{',W_n(\dbF)}$ and $\scrG_{W_n(\dbF)}$ coincide. From this and Faltings deformation theory (see Subsection 2.5 (a)) we easily get that $\scrD_{y,p^n}=\scrD_{y_{'},t_{'},p^n}$. As $n\ge n_0$ is arbitrary, by passing to limit $n\to\infty$ we conclude that to show that $g_0$ has the desired property we can assume that $\scrD_y=\scrD_{y_{'},t_{'}}$ and this case is obvious.

\smallskip
{\bf Step 3.} Let $\grc_{'}^{\text{cl}}$ be the image in $\scrN^{\text{cl}}_{k(v)}$ of the connected component $\grc_{'}$ of the basic (supersingular) locus $\grs_0$ which contains the point $y_{'}$. Similarly we define $\grc$ and $\grc^{\text{cl}}$ but working with $y$. We view both $\grc_{'}^{\text{cl}}$ and $\grc^{\text{cl}}$ as reduced, closed subschemes of $\scrN^{\text{cl}}_{k(v)}$.

Due to Step 2, $\grc_{'}^{\text{cl}}t_{n_0}t_{'}$ is isomorphic to $\grc^{\text{cl}}$ under the isomorphism of $\scrD_y^{\text{big}}$ defined by $g_0$. As $A$ is supersingular, it is easy to see that each $\dbF$-valued point of $\grc_{'}^{\text{cl}}$ is of the form $y_{'}(h_{'})$ for some element $h_{'}\in \grI(y_{'})$. Similarly, each $\dbF$-valued point of $\grc$ belongs to $o(y)$. Thus $\grc(\dbF)G(\dbA_f^{(p)})$ is a $G(\dbA_f^{(p)})$-invariant subset of $o(y)$ (and thus also of $o^{\text{big}}(y)$, cf. Corollary 2.5.5 and the property (iii)) equal to the subset $\grc_{'}^{\text{cl}}(\dbF)t_{n_0}t_{'}G(\dbA_f^{(p)})$ of $o^{\text{big}}(y)$.

Therefore, a simple count of $\dbF$-valued points (cf. Formula (3)) shows (as in the proof of Corollary 2.5.4 (b)) that the endomorphism property holds for $y_{'}$ and in fact that for each prime $l\neq p$, $\scrE_{',\dbQ_l}$ is the conjugate of $\scrE_{\dbQ_l}$ under the element $t_{n_0}t_{'}\in G_1^0(\dbA_f^{(p)})$. 

\smallskip
{\bf Step 4.} Due to the property (iv), the torus $Z^0(G)=Z^0(G_1^0)$ contains $Z(G^{\der})$ and $G^{\ad}$ is its own normalizer in $G_1^{\der}$. From this and the property (vi) we get that $G^0$ is its own normalizer in $G_1^0$. We consider the torsor of $G^0$ that parametrizes the existence of an element $j\in \scrE_1^0(\dbQ)$ that takes under inner conjugation $\scrE_0$ onto $\scrE$. This makes sense as $G^0$ is its own normalizer in $G_1^0$. This torsor defines a class $\gamma\in H^1(\dbQ,G^0)$ which is trivial at all primes $l\neq p$, cf. end of Step 3. As $\scrE^{\der}_{\dbR}$ and $\scrE_{',\dbR}^{\der}$ are compact subgroups of the compact group $\scrE_{1,\dbR}^0$, from the property (iv) we easily get that they are $\scrE_1^0(\dbR)$-conjugate. Thus $\gamma$ is also trivial at the infinite place of $\dbQ$. The existence of the element $g$ implies that the class $\gamma$ is also trivial at the prime $p$. The Hasse principle holds for $G^{\der}$ (cf. property (i)) and for $G^0/G^{\der}$ (cf. property (iv)) and thus also for $G^0$. Therefore the class $\gamma$ (i.e., the mentioned torsor) is trivial. Thus there exists an element $j\in \scrE_1^0(\dbQ)$ such that we have $j\scrE_{'}j^{-1}=\scrE$.

\smallskip
{\bf Step 5.} The existence of $j$ allows us to redo Step 1 as follows. It defines a $\dbQ$--isogeny from $(A,\lambda_A)$ to $(A_{'},\lambda_{A_{'}})$ under which we can identify $(M[{1\over p}],\phi,\psi_M)=(M_{'}[{1\over p}],\phi_{'},\psi_{M_{'}})$ in such a way that we have $\scrE_{'}=\scrE$ and there exists an element $\tilde g\in \scrE_1(\dbQ_p)$ which fixes $\psi_M$ and which takes each $t_{\alpha}$ to $t_{',\alpha}$. As $\scrE_{'}=\scrE$ we also have $\scrG_{B(\dbF)}=\scrG_{',B(\dbF)}$. Thus $\tilde g$ normalizes $\scrG_{B(\dbF)}$ and therefore we have $\tilde g\in\scrG(B(\dbF))$. This implies that $\tilde g\in\scrE(\dbQ_p)$. Thus we have $t_{\alpha}=t_{',\alpha}$ for all $\alpha\in\scrJ$. As above we had $h(M)=M=g(M)$, we can choose $\tilde g$ such that $\tilde g(M)=M_{'}$. By replacing $y$ with $y(\tilde g)$, we can assume that $\tilde g=1_M$. From this and the fact that $t_{\alpha}=t_{',\alpha}$ for all $\alpha\in\scrJ$, we get that there exists $t\in G(\dbA_f^{(p)})$ such that we have an identity $y_{'}=yt$ of $\dbF$-valued points of $\scrM$ (as $\scrE_{'}=\scrE$ such an element $t\in G_1(\dbA_f^{(p)})$ normalizes $G(\dbA_f^{(p)})$ and thus it belongs to $G(\dbA_f^{(p)})$). As we have $t_{\alpha}=t_{',\alpha}$ for all $\alpha\in\scrJ$, we have an identity $y_{'}=yt$ of $\dbF$-valued points of $\scrN$ (cf. Lemma 2.5.3). Thus $y_{'}\in o(y)$.\endproof

\bigskip\noindent
{\bf 4.4. Proposition.} {\it Let $y:\Spec\dbF\to\scrN$ be a basic point. We assume that one of the following two properties hold:

\medskip
{\bf (i)} either $A$ is supersingular and the Conjecture 1.2 holds for $y$ 

\smallskip
{\bf (ii)} or the adjoint group $G^{\ad}$ is simple and there exists a semisimple $\dbZ_{(p)}$-subalgebra of $\End(L_{(p)})$ such that its centralizer in $\pmb{GSp}(L_{(p)},\psi)$ is $G_{\dbZ_{(p)}}$.

\medskip
Then the isogeny property holds for $y$.}

\medskip
\proof 
We first assume that the property (i) holds. We consider the moduli space $\grs$ over $\dbF$ of quadruples $(\tilde A,\lambda_{\tilde A},(c_{\alpha,\tilde A})_{\alpha\in\scrJ},(\eta_{N,\tilde A})_{N\in \dbN\setminus p\dbN})$, where $(\tilde A,\lambda_{\tilde A})$ is a principally polarized supersingular abelian variety endowed with a family of algebraic cycles (up to numerical equivalence) $(c_{\alpha,\tilde A})_{\alpha\in\scrJ}$ and with compatible level-$N$ symplectic similitude structures $\eta_{N,\tilde A}$ for all $N\in\dbN\setminus p\dbN$, which satisfies the following two axioms modeled on $\scrN$ and $\grs_0$:

\medskip
{\bf (a)} For each prime $l$, the $l$-adic realization of $c_{\alpha,\tilde A}$ corresponds to $v_{\alpha}$ via the isomorphism $W\otimes_{\dbQ} \dbQ_l\arrowsim (\text{proj.}\text{lim.}_{n\in\dbN^*} \tilde A[l^n])\otimes_{\dbZ_l} \dbQ_l$ induced by the $\eta_{l^n,\tilde A}$'s with $n\in\dbN^*$. 

\smallskip
{\bf (b)} There exists a $\dbQ$--isogeny $(\tilde A,\lambda_{\tilde A},(c_{\alpha,\tilde A})_{\alpha\in\scrJ})\rightsquigarrow (A,\lambda_A,(c_{\alpha,A})_{\alpha\in\scrJ})$ under which the principally quasi-polarized Dieudonn\'e module of $(\tilde A,\lambda_{\tilde A})$ endowed with the family of tensors that are the crystalline realizations of $(c_{\alpha,\tilde A})_{\alpha\in\scrJ}$, can be identified with a quadruple of the form $(h(M),\phi,(t_{\alpha})_{\alpha\in\scrJ},\psi_{h(M)})$ for some element $h\in\grI(y)$.

\medskip
The moduli space $\grs$ is a reduced, closed, $G(\dbA_f^{(p)})$-invariant subscheme of $\scrM_{\dbF}$ and it a pro-\'etale cover of a projective $\dbF$-scheme. The set $\grs(\dbF)$ is exactly the set of $\dbF$-valued points of $\scrM$ defined by the points $y(h)G(\dbA_f^{(p)})$ with $h\in\grI(y)$. 

Let $\grc$ be the connected component of $\grs$ to which $y$ belongs. Due to the property (i), it is easy to see that $\grc$ is a closed subscheme of both $\scrN_{\dbF}$ and $\scrN_{\dbF}^{\text{cl}}$ with the property that each point of $\grc(\dbF)$ is of the from $y(h)$ for some element $h\in\grI(y)$ (cf. also Lemma 2.5.3 for the case of $\scrN_{\dbF}$). Thus to prove the proposition in the case when (i) holds, it suffices to show that the connected components of $\grs$ are permuted transitively by $G(\dbA_f^{(p)})$. 

But this follows from the following two facts:

\medskip
{\bf (c)} The $\dbF$-valued pivotal points of $\grs$ form one orbit under the action of $G(\dbA_f^{(p)})$.

\smallskip
{\bf (d)} Each connected component of $\grs$ has $\dbF$-valued pivotal points (see Definition 3.1.2).

\medskip
Property (c) is a direct consequence of the fact that the set of lattices $h(M)$ of $M[{1\over p}]$ with the property that $h\in\grI(y)$ and that the triple $(h(M),\phi,h\scrG h^{-1})$ is pivotal (i.e., is $\scrE(\dbQ_p)$ isomorphic to the one associated to a fixed pivotal point of $\grs_0(\dbF)$), are conjugate under $\scrE(\dbQ_p)$ and thus also under $\scrE(\dbQ)$ (cf. Subsection 2.2). 

To check (d), we first recall that each connected component $\grc_0$ of $\grs$ is a pro-\'etale cover of a projective $\dbF$-scheme. Based on this and Lemma 3.4.1, we get that there exists a stratum $\grl$ of the level $1$ stratification of $\grc_0$ which is non-empty of dimension $0$ (it is [Va8, Basic Thm. D and Rm. 12.4 (a)] that allows us to speak naturally about the level $1$ stratification of $\grc_0$; it is induced by the level $1$ stratification of $\scrN_{\dbF}$ in case $\grc_0\subset\scrN_{\dbF}$). But from Theorem A6 (c) of Appendix we get that the points of $\grl\cap\grc_0$ are pivotal. We conclude that $\grl$ is the stratum corresponding to pivotal points; thus the property (d) holds.

Thus if (i) holds we conclude that $\grs=\grs_0=\grc G(\dbA_f^{(p)})$.

We now assume that (ii) holds. The same proof applies, with the only modification that instead of algebraic cycles we have to use endomorphisms. We only have to add that the resulting moduli space $\grs$ is a disjoint union of a finite number of basic loci and thus it is a pro-\'etale cover of a projective $\dbF$-scheme, cf. Lemma 3.1.1.\endproof  

\bigskip\smallskip
\noindent
{\boldsectionfont 5. Proof of the Basic Theorem 1.5}
\bigskip 

In this section we prove the basic Theorem 1.5. If $E$ is a number field unramified over $p$, let $E_{(p)}$ be the normalization of $\dbZ_{(p)}$ in $E$; it is a finite, \'etale $\dbZ_{(p)}$-algebra. Let the simple, adjoint Shimura pair $(G_0,X_0)$ and the prime $p$ be as in Theorem 1.5. 

\bigskip\noindent
{\bf 5.1. Constructions.} We begin the proof of Theorem 1.5 by constructing the injective map $(G,X)\hookrightarrow (\pmb{GSp}(W,\psi),S)$ and the $\dbZ$-lattice $L$ of $W$. For the sake of clarity, the steps of the construction will be itemized from (a) to (h). We will construct only $(G,X)\hookrightarrow (\pmb{GSp}(W,\psi),S)$ and the $\dbZ_{(p)}$-lattice $L_{(p)}$ of $W$ (implicitly $L$ will be an arbitrary $\dbZ$-lattice of $W$ which extends $L_{(p)}$ and which is self dual with respect to $\psi$). 

\smallskip
{\bf (a)} Let $F_0$ be a totally real number field such that we have $G_0=\Res_{F_0/\dbQ} J_0$, with $J_0$ an absolutely simple, adjoint  group over $F_0$ (cf. [De2, Subsubsect. 2.3.4 (a)]). The field $F_0$ is unique up to isomorphisms (i.e., up to $\Gal(\dbQ)$-conjugation). As the group $G_{0,\dbQ_p}=\Res_{F_0\otimes_{\dbQ} \dbQ_p/\dbQ_p} J_{0,F_0\otimes_{\dbQ}\dbQ_p}$ is unramified, the $\dbQ_p$-algebra $F_0\otimes_{\dbQ} \dbQ_p$ is unramified. Thus the field $F_0$ is unramified above $p$.

\smallskip
{\bf (b)} If $(G_0,X_0)$ is of $A_1$, $B_n$, $C_n$, or $D_n^{\dbR}$ type, let $E_0$ be an arbitrary totally imaginary quadratic extension of $F_0$ for which the following two properties hold:

\medskip
\item
{\bf (b.i)} There exists a prime $l\neq p$ which splits in $F_0$ and $E_0$ has exactly $[F_0:\dbQ]$ primes that divide $l$, all but one of them being unramified over $l$. 

\item
{\bf (b.ii)} It is unramified above $p$ and moreover there exists no prime of $F_0$ that divides $p$ and that splits in $E_0$.

\medskip
If $(G_0,X_0)$ is of $A_n$ type with $n\ge 2$, let $E_0$ be the unique totally imaginary quadratic extension of $F_0$ such that the $\Gal(\dbQ)$-set $\Hom(E_0,\overline{\dbQ})$ is naturally identified with the $\Gal(\dbQ)$-set defined by the ending notes of the Dynkin diagram of $G_{0,\overline{\dbQ}}$ (see either [De2] or [Va3, Subsect. 2.2] for the action of $\Gal(\dbQ)$ on the Dynkin diagram of $G_{0,\overline{\dbQ}}$); as $G_{0,\dbQ_p}$ splits over an unramified finite field extension of $\dbQ_p$, the field $E_0$ is unramified over $p$. 

\smallskip
{\bf (c)} Let $G_{0,\dbZ_{(p)}}$ be the unique adjoint group scheme over $\dbZ_{(p)}$ such that $G_{0,\dbZ_{(p)}}(\dbZ_p)$ is an a priori fixed hyperspecial subgroup of $G_0(\dbQ_p)$, cf. [Va5, Lem. 2.3 (a)]. Let $G^{\sc}_{\dbZ_{(p)}}$ be the simply connected semisimple group scheme cover of $G_{0,\dbZ_{(p)}}$.

\smallskip
{\bf (d)} We assume $(G_0,X_0)$ is of $A_n$ type. We take $(G,X)\hookrightarrow (\pmb{GSp}(W,\psi),S)$ and $L$ to be as in [Va5, Prop. 3.2]. Thus the property (*) holds, $G^{\sc}_{\dbZ_{(p)}}$ is the derived group scheme of $G_{\dbZ_{(p)}}$, the centralizer $\scrB_{(p)}$ of $G_{\dbZ_{(p)}}$ in $\pmb{\GL}_{L_{(p)}}$ is a semisimple $\dbZ_{(p)}$-algebra, and $G_{\dbZ_{(p)}}$ is the centralizer of $\scrB_{(p)}$ in $\pmb{GSp}(L_{(p)},\psi)$. In particular, $(G,X)$ is a Shimura variety of PEL type and the injective map $(G,X)\hookrightarrow (\pmb{GSp}(W,\psi),S)$ is a PEL type embedding. We recall from [Va5, Prop. 3.2] and its proof, that $L_{(p)}$ has a natural structure of an $E_{0,(p)}$-module and $Z(G_{\dbZ_{(p)}})$ is a subtorus of $\Res_{E_{0,(p)}/\dbZ_{(p)}} \dbG_m$. More precisely, $Z(G_{\dbZ_{(p)}})$ is the largest subtorus $T_{\dbZ_{(p)}}$ of $\Res_{E_{0,(p)}/\dbZ_{(p)}} \dbG_m$ that is the extension of $Z(\pmb{\GL}_{L_{(p)}})$ by a torus $T_{c,\dbZ_{(p)}}$ which over $\dbR$ is compact and which is isogenous to $\Res_{F_{0,(p)}/\dbZ_{(p)}} [\Res_{E_{0,(p)}/F_{0,(p)}} \dbG_m]/\dbG_m$. Let 
$$W\otimes_{\dbQ} \dbC=\oplus_{j\in\Hom(E_0,\dbC)} W_j$$
be the direct sum decomposition such that $W_j$ is the largest complex vector subspace of $W\otimes_{\dbQ} \dbC$ on which $E_0$ acts as scalar multiplication via the embedding $j:E_0\hookrightarrow\dbC$. If $i:F_0\hookrightarrow\dbR$ is such that the group $J_0\times_{F_0,i} \dbR$ is compact and if $j_{i,1},j_{i,2}:E_0\hookrightarrow \dbC$ are the two embeddings that extend $i$, then one has an independent choice to make on the Hodge types of $W_{j_{i,1}}$ and $W_{j_{i,2}}$ defined by a (any) element $x\in X$ (cf. proof of [Va5, Prop. 3.2]):

\medskip\noindent
{\bf (d.i)} If $W\otimes_{\dbQ} \dbC=F^{-1,0}_x\oplus F^{0,-1}_x$ is the Hodge decomposition defined by $x$, then either we have $W_{j_{i,1}}\subseteq F^{-1,0}_x$ and $W_{j_{i,2}}\subseteq F^{0,-1}_x$ or we have $W_{j_{i,1}}\subseteq F^{0,-1}_x$ and $W_{j_{i,2}}\subseteq F^{-1,0}_x$.

\medskip
For most applications, it is irrelevant which choices we make. However, in Subsection 5.3 (f) below we will impose a technical condition (***) which might or might not hold, depending on the choices we make. Thus, if possible, the choices will be such that the technical condition (***) holds. 

\medskip
{\bf (e)} We assume that $(G_0,X_0)$ is of $B_n$ or $D_n^{\dbR}$ (resp. of $C_n$) type and that is compact. We will take $(G,X)\hookrightarrow (\pmb{GSp}(W,\psi),S)$ as in [De2, Prop. 2.3.10] and [Va1, Subsects. 6.5 and 6.6]. To make this very explicit, we will introduce an extra number field $E$ which contains $F_0$, which is unramified over $p$, and for which the group scheme $G^{\sc}_{E_{(p)}}:=G^{\sc}_{\dbZ_{(p)}}\times_{\Spec\dbZ_{(p)}} \Spec E_{(p)}$ is split. Let $G^{\sc}_{E_{(p)}}\hookrightarrow \pmb{\GL}_{\tilde W_{(p)}}$ be the standard spin (resp. symplectic) representation. Thus $\tilde W_{(p)}$ is a free $E_{(p)}$-module of rank $2^{n-1}$ (resp. $2n$). Let $L_{(p)}:=\tilde W_{(p)}\otimes_{F_{0,(p)}} E_{0,(p)}$ but viewed as a free $\dbZ_{(p)}$-module. The semisimple group scheme $G^{\sc}_{\dbZ_{(p)}}$ is naturally a closed subgroup scheme of $\pmb{\GL}_{L_{(p)}}$. Based on [Va5, Lem. 2.3 (b)], we can speak about the reductive group scheme $\tilde G_{\dbZ_{(p)}}$ which is the closed subgroup scheme of $\pmb{\GL}_{L_{(p)}}$ generated by $G^{\sc}_{\dbZ_{(p)}}$ and by the torus $T_{\dbZ_{(p)}}$ obtained as in (d). If $(G_0,X_0)$ is not of $D_n^{\dbR}$ type, then we take $G_{\dbZ_{(p)}}:=\tilde G_{\dbZ_{(p)}}$. If $(G_0,X_0)$ is of $D_n^{\dbR}$ type, then we take $G_{\dbZ_{(p)}}$ to be generated by $\tilde G_{\dbZ_{(p)}}$ and by the largest torus of $\pmb{Sp}(L_{(p)},\psi)$ which centralizes $G_{\dbZ_{(p)}}$. Let $G$ and $\tilde G$ be the generic fibres of $G_{\dbZ_{(p)}}$ and $\tilde G_{\dbZ_{(p)}}$ (respectively). Let $Z^{00}(G)$ be the generic fibre of $Z^{00}(G_{\dbZ_{(p)}}):=Z^0(\tilde G_{\dbZ_{(p)}})=T_{\dbZ_{(p)}}$. Similarly to (d), $L_{(p)}$ has a natural structure of an $E_{0,(p)}$-module and $Z^{00}(G_{\dbZ_{(p)}})$ is a subtorus of $\Res_{E_{0,(p)}/\dbZ_{(p)}} \dbG_m$. As in [Va1, Subsects. 6.5 and 6.6], [Va5, Prop. 3.2, end of proof], or [Va14, Part I, proof of Lem. 4.2.1] we argue that we can choose $\psi$ such that it restricts to a perfect alternating form on $L_{(p)}$. As $(G_0,X_0)$ is compact, it is easy to see that the property (b.i) implies (cf. Lemma 2.4.1 applied to $(T,T_1):=(\tilde G/G^{\der},\dbG_m)$) that:

\medskip\noindent
{\bf (e.i)} The smallest subtorus $\diamond$ of $Z^0(G)$ such that each element of $X$ factors through the extension to $\dbR$ of the subgroup of $G$ generated by $G^{\der}$ and $\diamond$, is $Z^{00}(G)$ (thus it is $Z^0(G)$ itself if and only if $(G_0,X_0)$ is not of $D_n^{\dbR}$ type).

\medskip
As $Z^{00}(G)$ is the extension of a torus isogenous to $\Res_{F/\dbQ} [(\Res_{E/F} \dbG_m)/\dbG_m]$ by $Z(\pmb{\GL}_W)$, from the property (b.ii) we get that:

\medskip\noindent
{\bf (e.ii)} The torus $Z^{00}(G)_{\dbQ_p}$ is the extension of an anisotropic torus over $\dbQ_p$ by $Z(\pmb{\GL}_W)_{\dbQ_p}$.

\medskip
{\bf (f)} We assume that $(G_0,X_0)$ is of $A_n$, $B_n$, $C_n$, or $D_n^{\dbR}$ type and that all simple factors of $G_{0,\dbR}$ are non-compact. If $(G_0,X_0)$ is of $A_n$ type, we also assume that $n$ is odd and that $G^{\der}_{\dbR}$ is a product of $\pmb{SU}({{n+1}\over 2},{{n+1}\over 2})$ groups. Then from [De2, Rms. 2.3.12 and 2.3.13] we get that in (d) and (e), we could take $G_{\dbZ_{(p)}}$ to be the reductive subgroup scheme of $\pmb{\GL}_{L_{(p)}}$ generated by $G^{\sc}_{\dbZ_{(p)}}$ and $Z(\pmb{\GL}_{L_{(p)}})$ (this makes sense, cf. [Va5, Lem. 2.3 (b)]). But for the sake of uniformity, we will take $G_{\dbZ_{(p)}}$ to be the reductive subgroup scheme of $\pmb{\GL}_{L_{(p)}}$ generated by $G^{\sc}_{\dbZ_{(p)}}$ and by $Z(\pmb{\GL}_{L_{(p)}})$ only in the case when $(G_0,X_0)$ is of $C_n$ type. The reflex field is $E(G,X)=\dbQ$, cf. loc. cit. Let $Z^{00}(G):=Z(\pmb{\GL}_W)$; the analogue of the property (e.i) obviously holds. 

If $(G_0,X_0)$ is of $C_n$ type, then the monomorphism $G^{\der}_{\dbZ_{(p)}}\times_{\Spec\dbZ_{(p)}} \Spec W(\dbF)\hookrightarrow \pmb{Sp}(L_{(p)}\otimes_{\dbZ_{(p)}} W(\dbF),\psi)$ is a product of standard composite monomorphisms $\pmb{Sp}_{2n}\hookrightarrow \pmb{Sp}_{2n}^s\hookrightarrow\pmb{Sp}_{2ns}$ with $s:={d\over {n[F_0:\dbQ]}}$. From this and the identity $Z^0(G)=\dbG_m$ we get that the centralizer $\scrB_{(p)}$ of $G_{\dbZ_{(p)}}$ in $\pmb{\GL}_{L_{(p)}}$ is a semisimple $\dbZ_{(p)}$-algebra and that $G_{\dbZ_{(p)}}$ is the centralizer of $\scrB_{(p)}$ in $\pmb{GSp}(L_{(p)},\psi)$; thus $(G,X)\to (\pmb{GSp}(W,\psi),S)$ is a PEL type embedding.

\medskip
{\bf (g)} We check that the properties 4.3 (i), (iv), and (v) hold. As $G^{\der}$ is simply connected, the property 4.3 (i) holds. Let $\scrB_{(p)}$ be the centralizer of $G_{\dbZ_{(p)}}$ in $\End(L_{(p)})$. It is a semisimple $\dbZ_{(p)}$-algebra. Let $G_{1,\dbZ_{(p)}}$ be the reductive subgroup scheme of $\pmb{GSp}(L_{(p)},\psi)$ that is the centralizer of $\scrB_{(p)}$ in $\End(L_{(p)})$. It is easy to see that the condition 4.3 (iv) holds. Let $G^0:=G\cap \pmb{Sp}(W,\psi)$ and $G_1^0:=G_1\cap\pmb{Sp}(W,\psi)$. We will check that the Hasse principle holds for $G^0/G^{\der}$ and $G^0_1/G_1^{\der}$. It is easy to see that $G^0/G^{\der}$ and $G^0_1/G_1^{\der}$ are isomorphic and thus we can work only with $G^0/G^{\der}$. We will only consider the $A_n$, $B_n$, and $C_n$ types, as the $D_n^{\dbR}$ type is very much the same. The torus $Z^0(G^0)$ is the maximal compact subtorus of $\Res_{E_0/\dbQ} \dbG_m$ and thus it sits in a short exact sequence $0\to \Res_{F_0/\dbQ} \dbG_m\to\Res_{E_0/\dbQ} \dbG_m\to Z^0(G^0)\to 0$. Therefore the Hasse principle holds for $Z^0(G^0)$. The torus $G^0/G^{\der}$ is isomorphic to $Z^0(G^0)$ (more precisely, it is the quotient of $Z^0(G^0)$ by $Z^0(G^0)[m]$, where $m$ is $2$ if $(G_0,X_0)$ is of $B_n$ or $C_n$ type and is $n+1$ if $(G_0,X_0)$ is of $A_n$ type).  From the last two sentences we get that the Hasse principle holds for $G^0/G^{\der}$. Thus the property 4.3 (v) holds as well.

\medskip
{\bf (h)} Regardless of what the type of $(G_0,X_0)$ is, properties (a) and (b) of Theorem 1.5 hold (for the last part of (b), cf. Lemma 2.5.6). In the next six subsections we will check one by one that all the properties (c) to (h) of Theorem 1.5 hold. 

\bigskip\noindent
{\bf 5.2. Proof of Theorem 1.5 (c).} In this subsection we assume that $(G_0,X_0)$ is not of $A_n$ type. Let $y:\Spec\dbF\to\scrN$ be a basic point and let $z:\Spec W(\dbF)\to\scrN$ be an arbitrary lift of it. We will use the canonical identification $Z^0(G)_{B(\dbF)}=Z^0(\scrG_{B(\dbF)})$ defined by crystalline realizations of $\dbQ$--endomorphisms of $\grA$ that correspond naturally to $\dbQ$--valued points of $Z^0(G)$.  Let $A_0$, $\dbF_{p^q}$, $(M_0,\phi_0,(t_{\alpha})_{\alpha\in\scrJ},\psi_{M_0})$, and $\pi_0$ be as in Subsection 2.1. As $y$ is basic, we can identify $\pi_0$ with a $\dbQ$--valued point of $Z^0(G)$ and thus also with an element of $Z^0(\scrG)(B(\dbF))$. Based on the property 5.1 (e.i), in fact we have $\pi_0\in Z^{00}(G)(\dbQ)$.

But $Z^{00}(G)$ is a subtorus of $\Res_{E_0/\dbQ} \dbG_m$ (see Subsection 5.1 (e)). As $\Res_{E_0/\dbQ} \dbG_m$ acts on $W$ via its natural structure of an $E_0$-vector space, we have canonical monomorphisms $\dbQ[\pi_0]\hookrightarrow E_0\hookrightarrow \End(W)$ of simple $\dbQ$--algebras. This implies that the double centralizer of $\dbQ[\pi_0]$ in $\End(W)$ is $\dbQ[\pi_0]$. Therefore the Lefschetz group of $A$ is equal to the center $T_0:=\Res_{\dbQ[\pi_0]/\dbQ} \dbG_m$ of the group scheme of invertible elements of $\End(A)\otimes_{\dbZ} \dbQ$. 

Thus to prove that the Tate conjecture holds for $A$, it suffices to show that the group $\{\pi_0^s|s\in\dbZ\}$ is dense in $T_0(\dbQ)$ (cf. [Mi6, Sect. 2]). 
The easiest way to prove this is to show that in fact we have $\dbQ[\pi_0]=\dbQ$ i.e., that the abelian variety $A$ is supersingular. If $Z^{00}(G)=\dbG_m$, then the fact that $\dbQ[\pi_0]=\dbQ$ is obvious. Thus we can assume that $Z^{00}(G)\neq \dbG_m$. As $(G_0,X_0)$ is of $B_n$, $C_n$, or $D_n^{\dbR}$ type and as $Z^{00}(\dbG_m)\neq \dbG_m$, from Subsection 5.1 (f) we get that $(G_0,X_0)$ is compact. Therefore $Z^{00}(G)_{\dbQ_p}$ is the extension of an anisotropic torus by $Z(\pmb{\GL}_W)_{\dbQ_p}=\dbG_m$, cf. property 5.1 (e.ii). Thus the only split subtorus of $Z^{00}(G)_{\dbQ_p}$ is $Z(\pmb{\GL}_W)_{\dbQ_p}$.

For Newton polygon slope quasi-cocharacters we refer to either [Pi] or [Va10, Subsect. 2.2]. The Newton polygon slope quasi-cocharacter of $(M,\phi)$ is a quasi-cocharacter of $Z^{00}(\scrG_{B(\dbF)}):=Z^{00}(G)_{B(\dbF)}$ which is the extension to $B(\dbF)$ of a quasi-cocharacter of $Z^0(G)_{\dbQ_p}$ and thus it factors through the extension to $B(\dbF)$ of the only split subtorus $Z(\pmb{\GL}_W)_{\dbQ_p}$ of $Z^{00}(G)_{\dbQ_p}$. Therefore the Newton polygon slope quasi-cocharacter of $(M,\phi)$ factors through $Z(\pmb{\GL}_{M[{1\over p}]})$. Thus all Newton polygon slopes of $(M,\phi)$ are equal and therefore are ${1\over 2}$. Thus $A$ is a supersingular abelian variety and therefore the property 4.3 (vi) holds as well. This ends the proof of Theorem 1.5 (c).${}^4$ $\vfootnote{4}{If $(G_0,X_0)$ is of $A_n$ type with $n>1$ and there exists no prime of $F_0$ which divides $p$ and which splits in $E_0$, then a similar argument shows that the abelian variety $A$ is supersingular and thus that the Tate conjecture holds for it.}$\endproof

\bigskip\noindent
{\bf 5.3. Proof of Theorem 1.5 (d).}
In this subsection we prove Theorem 1.5 (d). For the sake of clarity, we  itemize the six main tricks (ideas) we use  to prove Theorem 1.5 (d). In the first five tricks we assume that $(G_0,X_0)$ is not of $A_n$ type. The last trick explains the modifications required to be made in order to handle the $A_n$ type case as well. 

\smallskip
{\bf (a) The diagonal trick.} We assume that $(G_0,X_0)$ is of $B_n$, $C_n$, or $D_n^{\dbR}$ type. Let $m\in\dbN^*$. To prove Theorem 1.5 (d) we can replace $A$ by a power $A^m$ of it. In the crystalline context, this corresponds to a replacement of $M$ by $M^m$ and to a replacement of each $t_{\alpha}\in\scrT(\End(M[{1\over p}]))$ by a {\it diagonal} tensor $t_{\alpha}^{(m)}\in\scrT(\End(M^m[{1\over p}]))$. Like if $t_{\alpha}\in\End(M[{1\over p}])$, then $t_{\alpha}^{(m)}:=(t_{\alpha},\ldots,t_{\alpha})\in \End(M[{1\over p}])^m\subseteq \End(M^m[{1\over p}])$. At the level of Betti homologies with $\dbZ_{(p)}$ coefficients, this corresponds to a replacement of $L_{(p)}$ by $L^m_{(p)}$ and of $v_{\alpha}$ by $v_{\alpha}^{(m)}$. We recall from Subsection 5.1 (e), that we have $L_{(p)}=\tilde W_{(p)}\otimes_{F_{0,(p)}} E_{0,(p)}$. It is convenient to introduce a totally real number field $F_1$ which is unramified over $p$, which contains $F_0$, and for which we have $[F_1:F_0]=m$. Thus $F_{1,(p)}$ is a free $F_{0,(p)}$-module of rank $m$ and therefore we can identify
$$L^m_{(p)}=L_{(p)}\otimes_{F_{0,(p)}} F_{1,(p)}=\tilde W_{(p)}\otimes_{F_{0,(p)}} E_{0,(p)}\otimes_{E_{0,(p)}}  (E_{0,(p)}\otimes_{F_{0,(p)}} {F_{1,(p)}})\leqno (12)$$
as $G_{\dbZ_{(p)}}$-modules. Under these identifications we have 
$$v_{\alpha}^{(m)}=v_{\alpha}\otimes 1\in\scrT(\End(W))\otimes_{F_0} F_1=\scrT(\End(W\otimes_{F_0} F_1)).\leqno (13)$$
\indent
{\bf (b) The relative PEL trick.} Let $G_{1,\dbZ_{(p)}}^{\sc}:=\Res_{F_{1,(p)}/\dbZ_{(p)}} J_{0,F_{1,(p)}}^{\sc}$, where $J_{0,F_{0,(p)}}$ is the adjoint group scheme over $F_{0,(p)}$ whose generic fibre is $J_0$ and for which we have $G_{\dbZ_{(p)}}^{\sc}=\Res_{F_{0,(p)}/\dbZ_{(p)}} J_{0,F_{0,(p)}}^{\sc}$ and where $J_{0,F_{1,(p)}}^{\sc}$ is the simply connected semisimple group scheme cover of $J_{0,F_{1,(p)}}:=J_{0,F_{0,(p)}}\times_{\Spec F_{0,(p)}} \Spec F_{1,(p)}$. We view naturally $G_{1,\dbZ_{(p)}}^{\sc}$ as a simply connected semisimple subgroup scheme of $\pmb{\GL}_{L^m_{(p)}}$. Let $G_{1,\dbZ_{(p)}}$ be the reductive subgroup scheme of $\pmb{\GL}_{L^m_{(p)}}$ generated by $G_{1,\dbZ_{(p)}}^{\sc}$ and $Z^0(G_{\dbZ_{(p)}})$, cf. [Va5, Lem. 2.3 (b)]. Let $G_1$ be the generic fibre of $G_{1,\dbZ_{(p)}}$. We have $G^{\der}(\dbQ)=J_{0,F_{0,(p)}}^{\sc}(F_0)$ and $G_1^{\der}(\dbQ)=J_{0,F_{0,(p)}}^{\sc}(F_1)$. From this and (13), as $G^{\der}$ fixes each $v_{\alpha}$, we get that $G_1^{\der}(\dbQ)$ and thus also $G_1^{\der}$ fixes each $v_{\alpha}^{(m)}$. Moreover $Z^0(G)$ fixes each $v_{\alpha}$ and thus also each $v_{\alpha}^{(m)}$. From the last two sentences we get that:

\medskip
{\bf (b.i)} The reductive group $G_1$ fixes $v_{\alpha}^{(m)}$ for all $\alpha\in\scrJ$. 

Let $X_1$ be the $G_1(\dbR)$-conjugacy class of the composite of the natural monomorphism $G_{\dbR}\hookrightarrow G_{1,\dbR}$ with any monomorphism $\Res_{\dbC/\dbR} \dbG_m\hookrightarrow G_{\dbR}$ that defines an element of $X$. We have natural injective maps of Shimura pairs
$$(G,X)\hookrightarrow (G_1,X_1)\hookrightarrow (\pmb{GSp}(W^m,\psi_m),S_m),$$
where $\psi_m$ is a perfect, alternating form on $L^m$ such that we get naturally a Siegel pair $(\pmb{GSp}(W^m,\psi_m),S_m)$.${}^5$ $\vfootnote{5}{One can always modify the construction of the injective map $(G,X)\hookrightarrow (\pmb{GSp}(W,\psi),S)$ such that it factors through the natural injective map $(G,X)\hookrightarrow (G_1,X_1)$ (for this, one would only need to enlarge the field $E$ such that it also contains $F_1$). Thus for the proof of Theorem 1.5, one can always assume that $m=1$ (but in such a case one would either have to write $A$ as an $m$-th power of an abelian variety $B$ or work with the weaker form of the Conjecture 1.2 in which one replaces ``regardless of the choice" by ``we can choose").}$ The composite monomorphisms $G_{\dbZ_{(p)}}\hookrightarrow G_{1,\dbZ_{(p)}}\hookrightarrow \pmb{\GL}_{L^m_{(p)}}$ form a relative PEL situation (embedding) in the sense of [Va1, Subsubsect. 4.3.16] i.e., we have the following property (cf. loc. cit.):
\medskip
\item
{\bf (b.ii)} The group scheme  $G_{\dbZ_{(p)}}$ is the closed subgroup scheme of $G_{1,\,\dbZ_{(p)}}$ that fixes a semisimple $\dbZ_{(p)}$-subalgebra of $\End(L^m_{(p)})$.

\medskip
{\bf (c) The $A_1$ twisting trick.} By enlarging the finite field extension $F_1$ of $F$, we can assume that there exist injective maps
$$(G_2,X_2)\hookrightarrow (G_1,X_1)\hookrightarrow (\pmb{GSp}(W^m,\psi_m),S_m)$$
of Shimura pairs, such that the following two properties hold:
\medskip
\item
{\bf (c.i)} The schematic closure $G_{2,\dbZ_{(p)}}$ of $G_2$ in $G_{1,\dbZ_{(p)}}$ is a reductive group scheme such that we have $Z^0(G_{2,\dbZ_{(p)}})=Z^{00}(G_{\dbZ_{(p)}})$, the derived group $G_2^{\der}$ is simply connected, and $G_2^{\ad}$ is an adjoint group of the form $\Res_{F_0/\dbQ} J_2$, where $J_2$ is an absolutely simple, adjoint group of $A_1$ Dynkin type.

\item
{\bf (c.ii)} The reductive group scheme $G_{2,\dbZ_{(p)}}$ is the closed subgroup scheme of $\pmb{GSp}(L^m_{(p)},\psi_m)$ that fixes a semisimple $\dbZ_{(p)}$-subalgebra of $\End(L^m_{(p)})$. Moreover the derived group scheme $G_{2,\dbZ_p}^{\der}$ is a product of Weil restrictions of $\pmb{SL}_2$ group schemes. 
\medskip
The existence of the injective map $(G_2,X_2)\hookrightarrow (G_1,X_1)$ is argued using the same twisting arguments as in [Va3, Subsect. 4.5]. In fact, as $J_2$ we can take any absolutely simple, adjoint group over $F_0$ which has the following two properties:
\medskip
\item
{\bf (c.iii)} It is unramified over $\dbQ_p$.

\item
{\bf (c.iv)} For each embedding $i_0:F_0\hookrightarrow \dbR$ such that the group $J_0\times_{F_0,i_0} \dbR$ is non-compact (resp. it is compact), the group $J_2\times_{F_0,i_0} \dbR$ is a $\pmb{PGL}_2$ group (resp. it is compact).
\medskip
If $(G_0,X_0)$ is of $B_n$ type, let $s:=2n+1\ge 3$. If $(G_0,X_0)$ is of $D_n^{\dbR}$ type, let $s:=2n\ge 8$. Let $J_2^{\sc}$ be the simply connected semisimple group cover of $J_2$. The monomorphism $G_2^{\der}\to G_1^{\der}$ is a composite of the natural monomorphism $G_2^{\der}=\Res_{F_0/\dbQ} J_2^{\sc}\hookrightarrow \Res_{F_1/\dbQ} J_{2,F_1}^{\sc}$ with a monomorphism of the form $\Res_{F_1/\dbQ} J_{2,F_1}^{\sc}\hookrightarrow G_1^{\der}$ which is an inner twist of the Weil restriction from $F_1$ to $\dbQ$ of a composite monomorphism of the form:
\medskip
\item
{\bf (c.v)} $\pmb{Spin}_{2,1}\hookrightarrow \pmb{Spin}_{s-2,2}$, if $(G_0,X_0)$ is of $B_n$ or $D_n^{\dbR}$ type;

\item
{\bf (c.vi)} $\pmb{SL}_2\hookrightarrow \pmb{SL}_2^n\hookrightarrow \pmb{Sp}_{2n}$ (the first one being a diagonal monomorphism), if $(G_0,X_0)$ is of $C_n$ type.
\medskip
The fact that such an inner twist exists is guaranteed by [Va3, Lemmas 4.4 and 4.6].

Let $T_{2,\dbZ_{(p)}}$ be a maximal torus of $G_{2,\dbZ_{(p)}}$  which over $\dbR$ has $\dbR$-rank $1$ and whose extension to $\dbZ_p$ has an image in $G_{2,\dbZ_p}^{\ad}$ which is an anisotropic torus. Let $T_2$ be the generic fibre of $T_{2,\dbZ_{(p)}}$. We can choose $T_2$ such that there exists a prime $\tilde l$ distinct from $p$ and from the prime $l$ of the property 5.1 (b.i), for which the following property holds:

\medskip
\item
{\bf (c.vii)} The prime $\tilde l$ splits in $E_0$ and the image of $T_{2,\dbQ_{\tilde l}}$ in $G_{2,\dbQ_{\tilde l}}^{\der}$ is a product of anisotropic tori of dimension $1$, all but one of them being split after extension to the unramified quadratic extension of $\dbQ_{\tilde l}$. 

\medskip
Let $x_2:\Res_{\dbC/\dbR} \dbG_m\hookrightarrow T_{2,\dbR}$ be a homomorphism that defines an element of $X_2$, cf. [Ha, Lem. 5.5.3]. Let $\scrT_2$ be the integral canonical model of the Shimura quadruple $(T_2,\{x_2\},T_{2,\dbZ_{(p)}}(\dbZ_p),v_2^{T_2})$, where $v_2^{T_2}$ is a prime of $E(T_2,\{x_2\})$ that divides $v$. 

\medskip
{\bf (d) The Tate--Hodge trick.} We check that the smallest subtorus of $T_2$ such that $x_2$ factors through its extension to $\dbR$, is $T_2$ itself. If $(G_0,X_0)$ is not compact, then $Z^{00}(G)=\dbG_m$ (cf. Subsection 5.1 (f))) and the statement follows easily from the property (c.vii), cf. Lemma 2.4.1 applied to $(T,T_1):=(T_2,Z^{00}(G))$. If $(G_0,X_0)$ is compact, then the statement follows from properties 5.1 (e.i) and (c.vii) (cf. Lemma 2.4.1 applied again to $(T,T_1):=(T_2,Z^{00}(G))$). 

From the previous paragraph we get that $T_2$ is the Mumford--Tate group of every abelian variety $\scrA_2$ over $B(\dbF)$ which is associated naturally to a composite morphism
$z_2:\Spec W(\dbF)\to\scrT_2\to\scrN_m\to \scrM_m$, where $\scrM_m$ and $\scrN_m$ are the integral canonical models of $(\pmb{GSp}(W^m,\psi_m),S_m,\pmb{GSp}(L^m_{(p)},\psi_m)(\dbZ_p),p)$ and $(G_1,X_1,G_{1,\dbZ_{(p)}}(\dbZ_p),v_1)$ (respectively) and where $\scrT_2\to\scrN_m\to \scrM_m$ are the functorial morphisms (see [Va5, Subsubsect. 2.4.2], cf. [VZ, Cor. 5]). Here $v_1$ is the prime of $E(G_1,X_1)$ divided by $v_2^{T_2}$; it divides $v$. Due to the property (c.ii), the torus $T_2$ is the subgroup of $\pmb{GSp}(W^m,\psi_m)$ that fixes a semisimple $\dbQ$--subalgebra of $\End(W^m)$. From this and Lemma 2.7 we get that the Hodge--Tate property holds for the point $y_2:\Spec\dbF\to \scrT_2\to\scrN_m$ defined naturally by $z_2$.

\smallskip
{\bf (e) The basic point trick.} The point $y_2:\Spec\dbF\to \scrT_2\to\scrN_m$ is basic (cf. the anisotropic property of $\text{Im}(T_{2,\dbZ_p}\to G_{2,\dbZ_p}^{\ad})$ of (c); to be compared with [Va10, Subsubsect. 4.2.2]). To $y_2:\Spec\dbF\to\scrN_m$ corresponds a septuple 
$$(A_2,\lambda_{A_2},M_2,\phi_2,(t_{2,\alpha})_{\alpha\in\scrJ_2},\psi_{M_2},(\eta_{N,y_2})_{N\in \dbN\setminus p\dbN})$$ 
which is analogous to the septuple $(A,\lambda_A,M,\phi,(t_{\alpha})_{\alpha\in\scrJ},\psi_M,(\eta_{N,y})_{N\in \dbN\setminus p\dbN})$. Based on the property (b.i), we can assume that $\scrJ_2$ is a set which contains $\scrJ$ and that for each $\alpha\in\scrJ$ the tensor $t_{2,\alpha}$ is the crystalline realization of suitable Hodge cycles which whose Betti realizations correspond naturally to the tensor $v_{\alpha}^{(m)}$ introduced in (b) above. 

We know that there exists an algebraic cycle $c_{\alpha^{(m)},A_2}$ on $A_2$ whose crystalline realization is $t_{2,\alpha}^{(m)}$, cf. Lemma 2.7. Thus Conjecture 1.2 holds for the point $y_2:\Spec\dbF\to\scrN_m$. Therefore the endomorphism property holds for $y_2$ (cf. Proposition 4.2) and the isogeny property holds for $y_2$ (cf. Proposition 4.4). Thus all the six properties 4.3 (i) to (vi) hold in the context of $(G_2,X_2)\hookrightarrow (\pmb{GSp}(W^m,\psi_m),S_m)$, $L^m_{(p)}$, and $y_2:\Spec\dbF\to\scrN_m$ instead of $(G,X)\hookrightarrow (\pmb{GSp}(W,\psi),S)$, $L_{(p)}$, and $y:\Spec\dbF\to\scrN_m$ (above we have checked that properties 4.3 (i), (iv), (v), and (vi) hold in the last context of $(G,X)$ but the same arguments show that they hold as well in the context of $(G_1,X_1)$). Thus from Proposition 4.3 applied in the context of $(G_2,X_2)\hookrightarrow (\pmb{GSp}(W^m,\psi_m),S_m)$, $L^m_{(p)}$, and $y_2:\Spec\dbF\to\scrN_m$, we get that all the basic points $\Spec\dbF\to\scrN_m$ form one isogeny set $o(y_2)$. In particular, the point $y:\Spec\dbF\to\scrN_m$ defend naturally by $y:\Spec\dbF\to\scrN$ belongs to the isogeny set $o(y_2)$ and moreover the isogeny property holds for the point $y:\Spec\dbF\to\scrN_m$. 

The fact that $y:\Spec\dbF\to\scrN_m$ belongs to $o(y_2)$ implies that each $t_{\alpha}^{(m)}$ is the crystalline realization of an algebraic cycle $c_{\alpha^{(m)},A^m}$ on $A^m$. From this and (a) we get that Conjecture 1.2 holds for the basic point $y:\Spec\dbF\to\scrN$. This ends the proof of Theorem 1.5 (d) for the $B_n$, $C_n$, and $D_n^{\dbR}$ types. 

\smallskip
{\bf (f) The $A_n$ type variational trick.} We now assume that $(G_0,X_0)$ is of $A_n$ type and we explain the modifications needed to be made in (a) to (e) above in order to get that Theorem 1.5 (d) continues to holds for the $A_n$ type. If $n=1$, it is more practical to think of the $A_1$ type as a $C_1$ type and thus this case is handled by (a) to (e) above. Thus until Subsubsection 5.3.1 we can assume that $n>1$. The technical condition mentioned in Theorem 1.5 (d) is as follows.

\medskip
{\bf (***)} If neither $n$ is odd nor $G^{\der}_{\dbR}$ is a product of $\pmb{SU}({{n+1}\over 2},{{n+1}\over 2})$ groups, then we assume that there exists a subfield $F_2$ of $F_0$ and a totally imaginary quadratic extension $E_2$ of $F_2$ contained in $E_0$, such that the smallest subtorus $\diamond$ of $Z^0(G)$ with the property that each element of $X$ factors through the extension to $\dbR$ of the subgroup of $G$ generated by $G^{\der}$ and $\diamond$, is the subtorus $Z^{00}(G)$ of $\Res_{E_0/\dbQ} \dbG_m$ generated by $\dbG_m$ and by the maximal compact subtorus of $\Res_{E_2/\dbQ} \dbG_m$. 

\medskip
For instance, if $F_0=\dbQ$ then we can take $(F_2,E_2):=(F_0,E_0)$ and thus (***) holds. Let $m\in\dbN^*$. We will choose $F_1$ such that the properties listed in (a) hold. The only modification required to be made to (a), is provided by the single identity 
$$L^m_{(p)}=L_{(p)}\otimes_{F_{0,(p)}} F_{1,(p)}$$
(to be compared with  (12) where we had two identities) of $G_{\dbZ_{(p)}}$-modules. No modification is required to be made to (b). Referring to (c), one requires to take $(G_2,X_2)$ to be of $A_n$ type, to have $G_2^{\ad}=\Res_{F_2/\dbQ} J_2$, where $J_2$ is an absolutely simple group of $A_n$ Dynkin type, to have $E(G_2^{\ad},X_2^{\ad})=E(G_0,X_0)$, and to have $Z^0(G_2)=Z^{00}(G)$. In fact, as $J_2$ we can take any absolutely simple, adjoint group over $F_2$ which has the following three properties:
\medskip
\item
{\bf (f.i)} It is unramified over $\dbQ_p$.

\smallskip
\item
{\bf (f.ii)} The number field $E_2$ is the analogue of $E_0$ but obtained working with $(G_2^{\ad},X_2^{\ad})$ instead of with $(G_0,X_0)$ (i.e., the $\Gal(\dbQ)$-set $\Hom(E_2,\overline{\dbQ})$ is naturally identified with the $\Gal(\dbQ)$-set defined by the ending notes of the Dynkin diagram of $G^{\ad}_{2,\overline{\dbQ}}$). 

\smallskip
\item
{\bf (f.iii)} For each embedding $i_0:F_0\hookrightarrow \dbR$, by denoting $i_2:F_2\hookrightarrow\dbR$ the restriction of $i_0$ to $F_2$, the group $J_2\times_{F_2,i_2} \dbR$ is isomorphic to $J_0\times_{F_0,i_0} \dbR$.

\medskip
The main thing that defines $F_1$ is the requirement that the groups $J_0\times_{F_0} F_1$ and $J_2\times_{F_2} F_1$ are isomorphic. By assuming that the technical condition (***) holds, we are indirectly imposing the following condition on choices pertaining to the property 5.1 (d.i).

\medskip
{\bf (f.iv)} If $i_2:F_2\hookrightarrow \dbR$ is an embedding such that the group $J_2\times_{F_2,i_2} \dbR$ is compact and if $j_2:E_2\hookrightarrow \dbC$ is a fixed extension of $i_2$, then (due to (***)) for each embedding $j:E_0\hookrightarrow \dbR$ which is the restriction to $E_0$ of an embedding $F_1\otimes_{F_0} E_0\hookrightarrow \dbC$ that extends $j_2$, the Hodge type of $W_j$ is the same (e.g., for all such $j$ we have $W_j\subseteq F^{0,-1}_x$ for any $x\in X$).  

\medskip
Properties (c.i) and (c.ii) continue to hold but with $\pmb{SL}_2$ replaced by $\pmb{SL}_{n+1}$.

We will choose the maximal torus $T_2$ of $G_2$ such that there exists a prime $\tilde l>>0$ distinct from $p$ and such that the following property (analogous to (c.vii)) holds:

\medskip
\item
{\bf (f.v)} The prime $\tilde l$ splits in $E_2$ and the image of $T_{2,\dbQ_{\tilde l}}$ in $G_{2,\dbQ_{\tilde l}}^{\der}$ is a product of tori of dimension $n$ which are anisotropic and have no proper subtori, all but one of these tori being isomorphic. 

\medskip
Such tori of dimension over $\dbQ_{\tilde l}$ are constructed as follows. We choose a Galois extension $Q_{\tilde l}$ of $\dbQ_{\tilde l}$ whose Galois group is the symmetric group $S_{n+1}$. Let $\tilde Q_{\tilde l}$ be the fixed subfield of an $S_n$ subgroup of $S_{n+1}$. Then the torus $[\Res_{\tilde Q_{\tilde l}/\dbQ_{\tilde l}} \dbG_m]/\dbG_m$ over $\dbQ_{\tilde l}$ has dimension $n$, has no proper subtorus, and up to an isogeny is a maximal torus of $\pmb{SL}_{n+1}$. 

Only the first paragraph of (d) needs to be modified in the case when $(G_0,X_0)$ is not compact and $G_{0,\dbR}$ is not isomorphic to $[F_0:\dbQ]$ copies of $\pmb{SU}({{n+1}\over 2},{{n+1}\over 2})^{\ad}$. In such a case, from (***) and (f.v) one gets that the smallest subtorus of $T_2$ such that $x_2$ factors through its extension to $\dbR$, is $T_2$ itself (cf. Lemma 2.4.1 applied to $(T,T_1):=(T_2,Z^{00}(G))$). No modification is required to be made to (e). Thus if the technical condition (***) holds (resp. does not hold), then Theorem 1.5 (d) holds (resp. Theorem 1.5 (d) follows from the PEL type property of Subsection 5.1 (d)). This ends the proof of Theorem 1.5 (d) for the $A_n$ type with $n>1$. Thus Theorem 1.5 (d) holds.\endproof

\medskip\noindent
{\bf 5.3.1. Remark.} If $(G_0,X_0)$ is of $B_n$ or $C_n$ type, then $J_0$ has no outer automorphisms and based on this one can check that we can assume that $m=1$ (i.e., $F_1=F_0$); thus the Tate--Hodge property holds for basic points $\Spec\dbF\to\scrN$ (cf. Subsection 5.3 (d) and (e)).   

\bigskip\noindent
{\bf 5.4. Proof of Theorem 1.5 (e).}
We first assume that $y:\Spec\dbF\to\scrN$ is a basic point. The fact that the endomorphism property holds for $y$ follows from Theorem 1.5 (d) and Proposition 4.2. The fact that the isogeny property holds for $y$ and that we have an identity $\grs_0(\dbF)=o(y)$, is argued in the same way we argued in Subsection 5.3 (e) that the isogeny property holds for the point $y_2:\Spec\dbF\to\scrN_m$ and that all $\dbF$-valued basic points of $\scrN_m$ form one isogeny set $o(y_2)$. We now check that the unramified CM lift property holds for $y$. We know that there exists a special point $\tilde z:\Spec W(\dbF)\to\scrN$ whose closed point maps to a pivotal point $\tilde y\in\grs_0(\dbF)=o(y)$, cf. Lemma 3.1.4. As $o(\tilde y)=o(y)$ and as we have an identity $\grs_0(\dbF)=o(y)$, we get that the unramified CM lift property holds for $y$. Thus the first part of Theorem 1.5 (e) holds.

We now assume that $(G_0,X_0)$ is compact and that $y$ is not a basic point. Let $\grl$ and $\grc$ be as in Example 3.7. The schematic closure of $\grc$ in $\scrN_{k(v)}$ contains points of $\grs_0$, cf. Example 3.7. From this and  Theorem 4.1 we get that the isogeny property holds for $y$. Thus Theorem 1.5 (e) holds.\endproof

\bigskip\noindent
{\bf 5.5. Proof of Theorem 1.5 (f).} Let $H_0$ be a compact, open subgroup of $G(\dbA_f^{(p)})$ such that the natural morphism $f:\scrN\to \scrN^{\text{cl}}$ is a pro-\'etale cover of the morphism $f_{H_0}:\scrN/H_0\to \scrN^{\text{cl}}/H_0$. As $f_{H_0}$ is an isomorphism in characteristic $0$, there exists a smallest reduced, closed $k(v)$-subscheme $\grU^{\text{cl}}/H_0$ of $\scrN^{\text{cl}}/H_0$ such that $f_{H_0}$ is an isomorphism outside $\grU^{\text{cl}}/H_0$. Let $\grU$ and $\grU^{\text{cl}}$ be the reduced schemes of the inverse images of $\grU^{\text{cl}}/H_0$ to $\scrN$ and $\scrN^{\text{cl}}$ (respectively). As the morphism $\scrN\to \scrM_{O_{(v)}}$ induces formally closed embeddings at the level of completions of local rings of residue field $\dbF$ (cf. Corollary 2.5.1), we get that:

\medskip
{\bf (i)} The set $\grU(\dbF)$ is the union of those fibres of the map $f(\dbF):\scrN(\dbF)\to\scrN^{\text{cl}}(\dbF)$ which have at least two elements (points).

\medskip
Let $\grs_0$ be the basic locus of $\scrN_{k(v)}$. From the identity $\grs_0=o(y)$ and the fact that the endomorphism property holds for $y$ (cf. Theorem 1.5 (e)), we get (see Corollary 2.5.5 or Lemma 3.3): 

\medskip
{\bf (ii)} The intersection $\grU\cap\grs_0$ is the empty scheme.

\medskip
Obviously $\grU$ is $G(\dbA_f^{(p)})$-invariant and thus we can take $\scrL:=\scrN\setminus \grU$.

We now assume that $(G_0,X_0)$ is compact. If the scheme $\grU$ is non-empty, then the intersection $\grU\cap\grs_0$ is non-empty (cf. Lemma 3.6) and this contradicts the property (ii). Thus $\grU=\empty$ and therefore the morphism $f:\scrL=\scrN\to\scrN^{\text{cl}}$ is an isomorphism.\endproof

\bigskip\noindent
{\bf 5.6. Proof of Theorem 1.5 (g).}
Theorem 1.5 (g) follows from Subsection 1.4.1, cf. works of Zink and Milne (see [Zi] and [Mi1-5,7]). For reader's convenience, we recall the essence of [Zi] and [Mi1-5,7]. We know that the isogeny property holds for all points $y:\Spec\dbF\to\scrN$ (cf. Theorem 1.5 (e)) and thus we can speak about the isogeny sets $o(y)$ for which the Formulas (2) and (3) hold. We will use the notations of Subsections 1.4.2 to 1.4.6.  Let $T$ be a maximal torus of $G$ such that we have a short exact sequence $0\to Z(\pmb{\GL}_{W\otimes_{\dbQ} \dbR})\hookrightarrow T_{\dbR}\hookrightarrow T_c\to 0$, where the torus $T_c$ is compact. Let $x\in X$ be such that its image factors through $T_{\dbR}$. Thus we have an injective map $(T,\{x\})\hookrightarrow (G,X)$ of Shimura pairs. Let $\grT$ be the $\overline{\dbQ}/\dbQ$--groupoid associated to $T$. We have three basic steps.

\medskip
{\bf Step 1.} To the Shimura pair $(T,\{x\})$ one associates naturally a homomorphism $t_x:\grP\to\grT$ of $\overline{\dbQ}/\dbQ$--groupoids, cf. [Mi3, Cor. 3.34]. By composing $t_x$ with the monomorphism $\grT\hookrightarrow \grG$, we get a homomorphism $j_x:\grP\to\grG$. Each $G(\overline{\dbQ})$-conjugate of such a  $j_x$ is called an admissible homomorphism (as $G^{\der}$ is simply connected, one has a simple and general criterion [LR, Thm. 5.3] to characterize them more intrinsically). Let $\grJ$ be the set of isomorphism classes (up to $G(\overline{\dbQ})$-conjugation) of such admissible homomorphisms. To $j_x$ one associates a $\dbZ\Theta$-set $X_p(j_x)$, a $G(\dbA_f^{(p)})$-torsor $X^p(j_x)$ (under a right action), and the group $\scrE_{j_x}$ of elements of $G(\overline{\dbQ})$ which fix $j_x$ under inner conjugation. The Langlands--Rapoport conjecture predicts the existence of a $\dbZ\Theta\times G(\dbA_f^{(p)})$-equivariant bijection
$$\scrN(\dbF)\arrowsim\bigsqcup_{[j_x]\in\grJ} \scrE_{j_x}\backslash [X_p(j_x)\times X^p(j_x)].\leqno (14)$$
\indent
{\bf Step 2.} We take $y=y_x:\Spec\dbF\to\scrN$ to be the specialization of a point in the image of the functorial morphism $\Sh(T,\{x\})\to\Sh_{H}(G,X)$ of $E(G,X)$-schemes. The choice of $y$ allows us to identify:

\medskip
\item {\bf (i)}  $X^p(j_x)$ with $G(\dbA_f^{(p)})$ and $X_p(j_x)$ with the set $\grX(y)$ of those $W(\dbF)$-lattices of $M[{1\over p}]$ that are of the form $h(M)$ with $h\in \grI(y)$;

\smallskip
\item
{\bf (ii)} $\scrE_{j_x}$ with $\scrE(\dbQ)$ in such a way that the actions on $X_p(j_x)\times X^p(j_x)$ and $\grX(y)\times G(\dbA_f^{(p)})$ coincide (this is possible due to the assumption that the endomorphism property holds for $y:\Spec\dbF\to\scrN$). 

\medskip
These identifications are compatible with the $\dbZ\Theta\times G(\dbA_f^{(p)})$-actions. Therefore we have an identification of $\dbZ\Theta\times G(\dbA_f^{(p)})$-sets
$$\scrE_{j_x}\backslash [X_p(j_x)\times X^p(j_x)]=\scrE(\dbQ)\backslash [\grX(y)\times G(\dbA_f^{(p)})].\leqno (15)$$ 
\indent
{\bf Step 3.} Based on Formulas (2), (3), and (15), to prove that (14) holds it suffices to show that the natural map $\grJ\to\grK$
that takes $[j_x]$ to $o(y_x)=o(y)$ is a bijection. The fact that this map is surjective follows from the fact that the CM property holds for all $\dbF$-valued points of $\scrN$ (to be compared with [Mi7, Thm. 6.5]). The fact that this map is injective is well known (it follows from the fact that the Langlands--Rapoport conjecture holds for Siegel modular varieties, cf. [Mi2]). More precisely, if we have another special pair $(T_1,\{x_1\})$ of $(G,X)$ such that $o(y_{x_1})=o(y_x)$, then $o^{\text{big}}(y_{x_1})=o^{\text{big}}(y_x)$ and thus there exists $g\in \pmb{\text{Gsp}}(W,\psi)(\overline{\dbQ})$ such that the conjugate under $g$ of $j_{x_1}$ is $j_x$ and by considering the $l$-adic context for a prime $l\neq p$ we get that $g\in G(\overline{\dbQ_l})$ and thus that $g\in G(\overline{\dbQ})$; therefore $[j_{x_1}]=[j_x]$ (to be compared with [Mi7, Prop. 6.4 and Thm. 6.5]). As $(G_0,X_0)$ is compact, from Theorem 1.5 (f) we get that the natural map $\scrN(\dbF)=\scrN^{\text{cl}}(\dbF)\to\scrM(\dbF)$ is injective and thus loc. cit. can be applied to our present context. 
\endproof

\bigskip\noindent
{\bf 5.7. Proof of Theorem 1.5 (h).} If $m=1$, then the Tate--Hodge property holds for basic points $\Spec\dbF\to\scrN$ and thus from Theorem 1.5 (g) we get that the unconditional form of the Langlands--Rapoport conjecture holds for the basic locus of $\scrN_{k(v)}$. If $m>1$, then a similar argument shows that the unconditional form of the Langlands--Rapoport conjecture holds for the basic locus of $\scrN^m_{k(v)}$. But as we have $Z^0(G)=Z^0(G_1)$ and due to Theorem 1.5 (g), the fact that the unconditional form of the Langlands--Rapoport conjecture holds for the basic locus of $\scrN^m_{k(v)}$ implies that the unconditional form of the Langlands--Rapoport conjecture holds for the basic locus of $\scrN_{k(v)}$ (to be compared with [Mi5, p. 7]; the obstruction consists in a class in the cohomology group $H^1(\dbQ,Z^0(G_1))$). Thus Theorem 1.5 (h) holds. This ends the proof of Theorem 1.5.\endproof

\bigskip\smallskip
\noindent
{\boldsectionfont 6. Main results}
\bigskip 

In this section we prove our main results, as corollaries to the previous sections and to earlier works. Let $(G_1,X_1)$ be a Shimura pair. Let $p$ be a prime such that the group $G_{1,\dbQ_p}$ is unramified. In all that follows, each Shimura quadruple of the form $(G_1,X_1,H_1,v_1)$ will involve a hyperspecial subgroup $H_1$ of $G_1(\dbQ_p)$ and a prime $v_1$ of the reflex field $E(G_1,X_1)$ that divides $p$. Let $k(v_1)$ be the residue field of $v_1$ and let $\Theta_1$ be the Frobenius automorphism of $\dbF$ that fixes $k(v_1)$. 

Until the end we will assume that each simple factor of $(G_1^{\ad},X_1^{\ad})$ is of $A_n$, $B_n$, $C_n$, or $D_n^{\dbR}$ type (we emphasize that the group $G_1$ could be a torus; in such a case the adjoint group $G_1^{\ad}$ is trivial). Let $\scrN_1$ be the integral canonical model of $(G_1,X_1,H_1,v_1)$ (i.e., of $\Sh_{H_1}(G_1,X_1)$ over the local ring $O_{(v_1)}$ of $v_1$) and let $\scrN_1^{\ad}$ be the integral canonical model of $(G_1^{\ad},X_1^{\ad},H_1^{\ad},v_1^{\ad})$, cf. [Va1] and [Va5, Thm. 5.1] for $p\ge 5$ and cf. [Va14, Parts I and II] for $p\in\{2,3\}$. See [Va5, Subsubsect. 2.4.2] and [Va14, Part II] for the functorial morphisms between integral canonical models. We recall the following basic property (cf. [Va1] and [Va5, Thm. 5.1] for $p\ge 5$ and cf. [Va14, Part II, Thm. 1.6 (b)] for $p\in\{2,3\}$). 

\bigskip\noindent
{\bf 6.1. Theorem.} {\it The functorial morphism $\scrN_1\to\scrN^{\ad}_{1,O_{(v_1)}}$ is a pro-\'etale cover of its image.} 

\bigskip\noindent
{\bf 6.2. Main Corollary A.} {\it We assume that $(G_1,X_1)$ is compact. Then the Langlands--Rapoport conjecture holds for the integral canonical model $\scrN_1$ of the Shimura quadruple $(G_1,X_1,H_1,v_1)$ (i.e., for the set $\scrN_1(\dbF)$ acted upon by the group $\dbZ\Theta_1\times G_1(\dbA_f^{(p)})$).}

\medskip
\proof
We proceed by induction on the dimension of $X_1$. If $\dim(X_1)=0$, then $G_1$ is a torus and it is well known that the classical field theory implies that the Langlands--Rapoport conjecture holds for the set $\scrN_1(\dbF)$. The inductive step goes as follows. 

We can assume that $G_1^{\der}$ is simply connected, cf. [Mi3, Thm. 4.9]. Based on [Pf, Main Thm.], to prove the corollary we can replace the Shimura quadruple $(G_1,X_1,H_1,v_1)$ by another one $(G,X,H,v)$ with the properties that: (i) $G^{\der}=G_1^{\der}$ is simply connected and (ii) we have an identity $(G^{\ad},X^{\ad},H^{\ad},v^{\ad})=(G_1^{\ad},X_1^{\ad},H_1^{\ad},v_1^{\ad})$ between adjoint Shimura quadruples. Therefore (as in [Va5, proof of Thm. 5.1]) we can take $(G,X,H,v)$ to be a product of Shimura quadruples whose adjoints are simple. By induction, we can assume that $(G_0,X_0):=(G^{\ad},X^{\ad})$ is a simple, adjoint Shimura pair. Therefore we can assume that we are in the context of Theorem 1.5 and accordingly we will use $\scrN$ instead of $\scrN_1$. Based on Theorem 1.5 (g), it suffices to show that the endomorphism property and the unramified CM lift property hold for every point $y:\Spec\dbF\to\scrN$. To show this, we can assume that $y$ is not a basic point (cf. Theorem 1.5 (e)). 

Let $A_0$, $\dbF_{p^q}$, $\pi_0$, $(M_0,\phi_0,\psi_{M_0})$, and $\scrG_0$ be as in Subsection 2.1. As $y$ is not a basic point, we have $\pi_0\in \scrG(B(\dbF))\setminus Z(\scrG)(B(\dbF))$. Let $\tilde\scrC_{B(\dbF)}$ be the centralizer of $\pi_0$ in $\scrG_{B(\dbF)}$. The group $\tilde\scrC_{B(\dbF)}$ is a split reductive subgroup of $\scrG_{B(\dbF)}$ different from $\scrG_{B(\dbF)}$. Let $\tilde\scrC$ be the schematic closure of $\tilde\scrC_{B(\dbF)}$ in $\scrG$; it is a flat group scheme over $W(\dbF)$. We know that the isogeny property holds for $y$, cf. Theorem 1.5 (e). Thus to check that the endomorphism property and the unramified CM lift property hold for $y:\Spec\dbF\to\scrN$, we can replace $y$ by an arbitrary point in the set $o(y)$.

\medskip\noindent
{\bf Claim.} {\it By replacing $y$ with $y(h)$ for some element $h\in\grI(y)$, we can assume that $\tilde\scrC$ is a reductive group scheme and that there exists a Hodge cocharacter $\tilde\mu:\dbG_m\to\tilde\scrC$ of $(M,\phi,\tilde\scrC)$.}

\medskip
 We recall that by Hodge cocharacter we mean that we have a direct sum decomposition $M=F^1\oplus F^0$ such that $M=\phi({1\over p}F^1\oplus F^0)$  and $\dbG_m$ acts through $\tilde\mu$ on each $F^i$ via the weight $-i$. 
To prove the Claim, let $\scrT_{0,B(\dbF_{p^q})}$ be a maximal torus of $\scrG_{0,B(\dbF_{p^q})}$ of $\dbQ_p$-endomorphisms of $(M_0,\phi_0,\scrG_0)$ in the sense of [Va11, Def. 2.3 (c)] which splits over $B(\dbF)$; thus its Lie algebra is generated by elements fixed by both $\phi_0$ and $\pi_0$. By enlarging $\dbF_{p^q}$, we can assume that there exists a cocharacter $\mu_{0,B(\dbF_{p^q})}:\dbG_m\to \scrT_{0,B(\dbF_{p^q})}$ which produces a direct sum decomposition $M_0[{1\over p}]=V^1_0\oplus V^0_0$ with the properties that:

\medskip
{\bf (i)}  the torus $\dbG_m$ acts on each $V^i_0$ via the weight $-1$;

\smallskip
{\bf (ii)} the composite cocharacter $\mu_{0,B(\dbF)}:\dbG_m\to \scrG_{B(\dbF)}$ is $\scrG(B(\dbF))$-conjugate to the generic fibre of the cocharacter $\mu:\dbG_m\to\scrG$ of Subsubsection 1.1.4;

\smallskip
{\bf (iii)} the triple $(M_0[{1\over p}],V^1_0,\phi_0)$ is a filtered $F$-isocrystal over $\dbF_{p^q}$ (i.e., is an admissible filtered module over $B(\dbF_{p^q})$). 

\medskip
The existence of $\mu_{0,B(\dbF_{p^q})}$ follows from [Va13, Thm. 4.2 (b) and Cor. 4.1.1] if $(G_0,X_0)$ is of $B_n$ or $D_n^{\dbR}$ type, from [Va13, Thms. 8.4 (a) and Cor. 4.1.1] if $(G_0,X_0)$ is of $C_n$ type, and (to be compared with the proof of [Va13, Thms. 8.4 (a)]) from [Zi, Thm. 4.4] if $(G_0,X_0)$ is of $A_n$ type. From (1) and from [Va11, Cor. 4.5] we get that the $\dbZ_p$ structure of $(M_0,\phi,\scrG_0,(t_{\alpha})_{\alpha\in\scrJ})$ obtained as in [Va11, Subsect. 2.4], is isomorphic to $(L_p^{\vee},G_{\dbZ_p},(v_{\alpha})_{\alpha\in\scrJ})$. Thus, as $G^{\der}$ is simply connected, from [Va13, Lem. 2.4.5] we get that the Galois representation associated to $(M_0[{1\over p}],V^1_0,\phi_0)$ is of the form $\varrho:\Gal(B(\dbF_{p^q}))\to G_{\dbQ_p}(\dbQ_p)\leqslant \pmb{\GL}_{L_p^{\vee}}(\dbQ_p)$. The image $\text{Im}(\varrho)$ is a subgroup of the group of $\dbQ_p$-valued points of a maximal torus $T_{0,\dbQ_p}$ of $G_{\dbQ_p}$ which is the $\dbQ_p$-form of $\scrT_{0,B(\dbF_{p^q})}$ with respect to $(M_0[{1\over p}],\phi_0)$ introduced in [Va13, Ex. 2.2.1]; more precisely, $T_{0,\dbQ_p}$ corresponds to  $\scrT_{0,B(\dbF_{p^q})}$ via Fontaine comparison theory. As $T_{0,\dbQ_p}$ splits over $B(\dbF)$, it is unramified and thus it extends to a torus $T_{0,\dbZ_p}$ over $\dbZ_p$. Let $g\in G(\dbQ_p)$ be such that $T_{0,\dbZ_p}$ is a maximal torus of $g G_{\dbZ_p} g^{-1}$, cf. Lemma 2.4.2. 

As $\Gal(B(\dbF_{p^q}))$ is compact and as the hyperspecial subgroup $T_{0,\dbZ_p}(\dbZ_p)$ of $T_{0,\dbZ_p}(\dbQ_p)$ is the maximal compact subgroup of  $T_{0,\dbZ_p}(\dbQ_p)$, we have $\text{Im}(\varrho)\leqslant T_{0,\dbZ_p}(\dbZ_p)$. Thus to the $\dbZ_p$-lattice $g(L_p^{\vee})$ normalized by $\Gal(B(\dbF_{p^q}))$ corresponds a filtered $F$-crystal $(\tilde M_0,V^1_1\cap\tilde M_0,\phi_0)$ over $\dbF_{p^q}$, where $\tilde M_0$ is a $W(\dbF_{p^q})$-lattice of $M_0[{1\over p}]$. From [Va11, Thm. 1.2] applied in the \'etale context provided by $\varrho$ and from (1), we get that there exist isomorphisms $(\tilde M_0\otimes_W(\dbF_{p^q}) W(\dbF),(t_{\alpha})_{\alpha\in\scrJ})\arrowsim (M,(t_{\alpha})_{\alpha\in\scrJ})$ i.e., there exists an element $h\in\scrG(B(\dbF))$ such that $h(M)=\tilde M_0\otimes_W(\dbF_{p^q}) W(\dbF)$. Based on the property (ii) and the fact that $(h(M),\phi,h\scrG h^{-1})$ is a Shimura $F$-crystal, we get that $h\in\grI(y)$. From [Va11, Thm. 1.2] we get that $\scrT_{0,B(\dbF_{p^q})}$ extends to a maximal torus $\scrT_0$ of $\scrG_0$. Thus $\mu_{0,B(\dbF_{p^q})}:\dbG_m\to \scrT_{0,B(\dbF_{p^q})}$ extends to a cocharacter $\mu_0:\dbG_m\to \scrT_0$ which is a Hodge cocharacter of $(h(M),\phi,h\scrG h^{-1})$. The fact that the schematic closure $h\tilde\scrC h^{-1}$ of $\tilde\scrC_{B(\dbF)}$ in $h\scrG h^{-1}$ is a reductive group scheme follows from the fact that this schematic closure is the centralizer in $h\scrG h^{-1}$ of a subtorus of $\scrT_{0,W(\dbF)}$. Moreover, the cocharacter $\tilde\mu:\dbG_m\to h\tilde\scrC h^{-1}$ induced naturally by $\mu_{0,W(\dbF)}$ is a Hodge cocharacter of $(h(M),\phi,h\tilde\scrC h^{-1})$. Thus the Claim holds.

If $F^1$ is as above defined by $\tilde\mu:\dbG_m\to \tilde\scrC$, let $z:\Spec W(\dbF)\to\scrN$ be a lift of $y$ such that $F^1$ is the Hodge filtration of $M$ defined by the abelian scheme $\scrA:=z^*(\grA)$ over $W(\dbF)$ (cf. [Va14, Part I, Lem. 3.5.2]). If $p=2$ and the $2$-rank of $A$ is positive, then we choose $F^1$ and $z$ such that the $2$-divisible group of $\scrA$ is a direct sum of $2$-divisible groups over $\dbF$ whose special fibres are isoclinic (cf. [Va14, Part I, Lem. 3.5.2 and Thm. B7 (b)]). 

Let $\tilde\scrZ$ be the subtorus of $Z(\tilde\scrC)$ whose generic fibre is the extension to $B(\dbF)$ of the Frobenius torus of $\pi_0$. The Frobenius endomorphism $\pi_0$ leaves invariant $F^1$ and thus it lifts to an endomorphism of $\scrA$ (for $p=2$, cf. also the choice of $z$). This implies that the Mumford--Tate group of the generic fibre of $\scrA$ can be identified naturally with a reductive subgroup $G_{\scrA}$ of $G$ which is contained in the centralizer $\tilde C$ in $G$ of the generic fibre of a torus $Z_{\dbZ_{(p)}}$ of $G_{\dbZ_{(p)}}$ whose crystalline realization is precisely $\tilde\scrZ$. The centralizer $\tilde C_{\dbZ_{(p)}}$ of $Z_{\dbZ_{(p)}}$ in $G_{\dbZ_{(p)}}$ is a reductive group scheme whose generic fibre is $\tilde C$. Let $x_{\scrA}\in X$ be the composite monomorphism $\Res_{\dbC/\dbR} \dbG_m\hookrightarrow  G_{\scrA,\dbR}\hookrightarrow G_{\dbR}$, where $\Res_{\dbC/\dbR} \dbG_m\hookrightarrow  G_{\scrA,\dbR}$ is the Hodge monomorphism of $\scrA$. Let $\tilde G$ be the largest normal subgroup of $\tilde C$ such that the pair $(\tilde G,\tilde X)$ is a Shimura pair, where $\tilde X$ is the $\tilde G(\dbR)$-conjugacy class of $x_{\scrA}$ viewed as a monomorphism $\Res_{\dbC/\dbR} \dbG_m\hookrightarrow \tilde G_{\dbR}$. Thus $\tilde G^{\ad}$ is the product of those simple factors of $\tilde C^{\ad}$ which do not normalize $x_{\scrA}$ (i.e., of those simple factor whose extensions to $\dbR$ contain a non-trivial image of $x_{\scrA}$). As a conclusion, we have constructed an injective map $(\tilde G,\tilde X)\hookrightarrow (G,X)$ of Shimura pairs such that the following three properties hold:

\medskip
{\bf (iv)} there exists an element $\pi_0\in Z^0(\tilde G)(\dbQ)$ which is the Betti realization of the endomorphism $\pi_0$ of $\scrA$ and which does not belong to $Z(G)(\dbQ)$;

\smallskip
{\bf (v)} the group $\tilde G$ is a normal subgroup of the centralizer $\tilde C$ of $\pi_0$ in $G$ and the quotient group $\tilde C/\tilde G$ is semisimple;

\smallskip
{\bf (vi)} the schematic closures of $\tilde C$  and $\tilde G$ in $\pmb{\GL}_{L_{(p)}}$ (equivalently, in $G_{\dbZ_{(p)}}$) are reductive group schemes $\tilde C_{\dbZ_{(p)}}$ and $\tilde G_{\dbZ_{(p)}}$ (respectively) over $\dbZ_{(p)}$.

\medskip
The part of the property (vi) pertaining to $\tilde G_{\dbZ_{(p)}}$ follows from the part of the property (vi) pertaining to $\tilde C$ and from the fact that $\tilde G$ is a normal subgroup of $\tilde C$. Let $\tilde H:=\tilde G_{\dbZ_{(p)}}(\dbZ_p)$. Let $\tilde v$ be a prime of $E(\tilde G,\tilde X)$ that divides $v$. Let $\tilde\scrN$ be the integral canonical model of the Shimura quadruple $(\tilde G,\tilde X,\tilde H,\tilde v)$. We will choose $\tilde v$ such that $z:\Spec W(\dbF)\to\scrN$ factors naturally through the functorial morphism $\tilde\scrN\to\scrN$. We denote by $\tilde z:\Spec W(\dbF)\to\tilde\scrN$ and $\tilde y:\Spec\dbF\to\tilde\scrN$ the corresponding factorizations of $z:\Spec W(\dbF)\to\scrN$ and $y:\Spec\dbF\to\scrN$ (respectively). 

Due to the property (iv) and the fact that $G^{\der}$ is simply connected, the group $\tilde G^{\der}$ is simply connected and we have $\dim(\tilde X)<\dim(X)$. Thus each simple factor of $(\tilde G^{\ad},\tilde X^{\ad})$ is of $A_n$, $B_n$, $C_n$, or $D_n^{\dbR}$ type. Moreover, due to the property (iv), the group $G^{\der}\cap \tilde G$ is of the form $\Res_{F_0/\dbQ} \tilde J_0$, for some reductive subgroup $\tilde J_0$ of the simply connected semisimple group cover of $J_0$, where $F_0$ and $J_0$ are as in Subsection 5.1 (a). From this and the fact that $(G_0,X_0)$ is compact, we easily get that each simple factor of $(\tilde G^{\ad},\tilde X^{\ad})$ is compact. Based on this and the fact that $\dim(\tilde X)<\dim(X)$, by induction we get that the conditional Langlands--Rapoport conjecture holds for $\tilde\scrN$.

Due to the property (v), we have a direct sum decomposition of Lie algebras over $\dbQ$
$$\Lie(\tilde C)=\Lie(Z^0(\tilde G))\oplus\Lie(\tilde G^{\der})\oplus\Lie(\tilde C/\tilde G).\leqno (16)$$
Due to the property (v), the reductive subgroup of $\tilde C$ generated by $Z^0(\tilde G)$ and by the semisimple, normal subgroup of $\tilde C$ whose adjoint is $\tilde C/\tilde G$, is naturally a subgroup of $\scrE$. Due to this and (16), to prove that the endomorphism property holds for $y:\Spec\dbF\to\scrN$ it suffices to show that the endomorphism property holds for $\tilde y:\Spec\dbF\to\tilde\scrN$. But by induction, the Langlands--Rapoport conjecture holds for the set $\tilde\scrN(\dbF)$ and therefore the isogeny property holds for the basic point $\tilde y:\Spec\dbF\to\tilde\scrN$ (in fact we also get directly that the endomorphism property also holds for $\tilde y:\Spec\dbF\to\tilde\scrN$). Moreover, it is easy to see that Theorem 1.5 (e) implies that we have $o(\tilde y)=\tilde\grs_0(\dbF)$, where $\tilde\grs_0$ is the basic locus of $\tilde\scrN_{k(\tilde v)}$ (in other words, the changes of Shimura quadruples performed in the second paragraph of this proof preserve the property that the $\dbF$-valued points of basic loci form only one isogeny set). Thus to check in a simple way that the endomorphism property also holds for $\tilde y:\Spec\dbF\to\tilde\scrN$, let $\tilde z_1:\Spec W(\dbF)\to\tilde\scrN$ be a special point whose closed point maps to a pivotal point $\tilde y_1\in\tilde\grs_0(\dbF)=o(\tilde y)=o(\tilde y_1)$ (cf. Lemma 3.1.4). The endomorphism property holds for $\tilde y_1$ (cf. Lemma 3.5 (b)) and thus also for $\tilde y$. We conclude that the endomorphism property holds for $y:\Spec\dbF\to\scrN$. 

As the points in the set $o(\tilde y_1)=o(\tilde y)$ map (via the functorial morphism $\tilde\scrN\to\scrN$) to points in the set $o(y)$, we conclude that the composite of $\tilde z_1$ with the functorial morphism $\tilde\scrN\to\scrN$ is a special point $z_1:\Spec W(\dbF)\to\scrN$ that lifts a point $y_1\in o(y)$. Thus the unramified CM lift property holds for $y$. 

We conclude that the endomorphism property and the unramified CM lift property hold for all points $y:\Spec\dbF\to\scrN$. Thus the Langlands--Rapoport conjecture holds for $\scrN$, cf. Theorem 1.5 (g). This ends the induction and the proof of the corollary.\endproof 
 
\bigskip\noindent
{\bf 6.3. On Conjecture 1.2.} Let $(G_0,X_0)$ be an arbitrary simple, adjoint Shimura pair of $A_n$ type with $n>1$. We assume that the group $G_{0,\dbQ_p}$ is unramified and that $(G_0,X_0)$ is compact. In this subsection we will assume that the injective map $(G,X)\hookrightarrow (\pmb{GSp}(W,\psi),S)$ and the $\dbZ$-lattice $L$ of $W$ are as in Theorem 1.5 (i.e., as constructed in Subsection 5.1 (d)). In particular, we have a natural identity $(G^{\ad},X^{\ad})=(G_0,X_0)$. We will use the basic notations of the proof of Corollary 6.2. Let $y:\Spec\dbF\to\scrN$ be a point that is not basic. We will also use the notations of Subsection 6.2. We would like to get a criterion for when the Conjecture 1.2 holds for $y$. As the isogeny property holds for $y:\Spec\dbF\to\scrN$ (cf. Theorem 1.5 (e)), by replacing $y$ with $y(h)$ for some element $h\in\grI(y)$, we can assume that $\tilde\scrC$ is a reductive group scheme over $W(\dbF)$ and that there exists a Hodge cocharacter $\tilde\mu:\dbG_m\to\tilde\scrC$ which produces a direct sum decomposition $M=F^1\oplus F^0$ such that we have $M=\phi({1\over p}F^1\oplus F^0)$ (cf. proof of Corollary 6.2). Let $\tilde T$ be the smallest subtorus of $Z^0(\tilde G)$ with the property that each element of $\tilde X$ factors through the extension to $\dbR$ of the subgroup of $\tilde G$ generated by $\tilde G^{\der}$ and $\tilde T$.

\medskip\noindent
{\bf 6.3.1. Criterion.} {\it If $\tilde T=Z^0(\tilde G)$, then the Conjecture 1.2 holds for $y:\Spec\dbF\to\scrN$.}

\medskip
\proof
The proof of this is very much the same as the proof of Theorem 1.5 (d) for the $A_n$ type, with the assumption that $\tilde T$ is $Z^0(\tilde G)$ as a substitute for the technical condition (***). As in Subsection 5.3, we will allow the replacement of $A$ and $W$ with $A^m$ and $W^m$, with $m\in\dbN^*$. Not to complicate the notations, in what follows we will not mention each time that we are first performing such a replacement. One only has to show the existence of a maximal torus $\tilde U$ of $\tilde G$ that has the following four basic properties (to be compared with Lemma 2.4.1 and Subsection 5.3 (d) and (e)):

\medskip
{\bf (i)} the schematic closure of $\tilde U$ in $\pmb{\GL}_{L_{(p)}}$ is a torus over $\dbZ_{(p)}$;

\smallskip
{\bf (ii)} the image of $\tilde U_{\dbQ_p}$ in $\tilde G^{\ad}_{\dbQ_p}$ is anisotropic;

\smallskip
{\bf (iii)} there exists a point $\tilde x\in \tilde X$ such that the smallest subgroup $\circleddash$ of $\tilde G$ with the property that $\tilde x$ factors through $\circleddash_{\dbR}$, is $\tilde U$;

\smallskip
{\bf (iv)} the torus $\tilde U$ is the subgroup of $\pmb{GSp}(W,\psi)$ that fixes a semisimple $\dbQ$--subalgebra of $\End(W)$. 

\medskip
Let $s$ be the number of simple factors of $(\tilde G^{\ad},\tilde X^{\ad})$. If $s=0$ (i.e., if $\tilde G$ is a torus), then we can take $\tilde U=\tilde T=\tilde G=Z^0(\tilde G)$. Property (iii) holds (cf. the identity $\tilde T=Z^0(\tilde G)$), property (ii) is vacuous, property (i) holds due to the property 6.2 (vi), and property (iv) follows from Subsection 5.1 (d) and the property 6.2 (v). 

Suppose now that $s>0$. Then by induction on $t\in \{0,\ldots,s\}$ we show that there exists a reductive subgroup $\tilde G_t$ of $\tilde G$ such that the following five properties hold:

\medskip
{\bf (v)} the schematic closure of $\tilde G_t$ in $\pmb{\GL}_{L_{(p)}}$ is a reductive group scheme over $\dbZ_{(p)}$;

\smallskip
{\bf (vi)} we have $Z^0(\tilde G)\leqslant \tilde G_t$ and the image of $Z^0(\tilde G_t)_{\dbQ_p}$ in $\tilde G^{\ad}_{\dbQ_p}$ is anisotropic;

\smallskip
{\bf (vii)} the derived subgroup $\tilde G_t^{\der}$ is a normal subgroup of $\tilde G^{\der}$ and $\tilde G_t^{\ad}$ is a product of $s-t$ simple factors of $\tilde G^{\ad}$;
  
\smallskip
{\bf (viii)}  there exists a point $\tilde x_t\in \tilde X$ such that the smallest subgroup $\circleddash$ of $\tilde G_t$ with the property that $\tilde x$ factors through $\circleddash_{\dbR}$, is $\tilde G_t$;

\smallskip
{\bf (ix)} the group $\tilde G_t$ is the subgroup of $\tilde G$ that fixes $Z^0(\tilde G_t)$ (thus $\tilde G_t$ is the subgroup of $\pmb{GSp}(W,\psi)$ that fixes a semisimple $\dbQ$--subalgebra of $\End(W)$). 

\medskip
The existence of $\tilde G_t$ is argued as in Subsection 5.3 (f). Taking $\tilde U=\tilde G_s$, we get that all the four properties (i) to (iv) hold (cf. properties (v) to (ix)).\endproof

\bigskip\noindent
{\bf 6.4. Main Corollary B.} {\it We assume that $(G_1,X_1)$ is compact and we have an injective map $(G_1,X_1,H_1,v_1)\hookrightarrow (\pmb{GSp}(W,\psi),S,K_p,p)$ of Shimura quadruples, where  $(\pmb{GSp}(W,\psi),S,K_p,p)$ is the Shimura quadruple of a Siegel modular pair $(\pmb{GSp}(W,\psi),S)$. Let $\scrM$ be the integral canonical model of $(\pmb{GSp}(W,\psi),S,K_p,p)$ and let $\scrN_1\to\scrM_{O_{(v_1)}}$ be the functorial morphism. We also assume that there exists a $\dbZ_{(p)}$-lattice $L_{(p)}$ which is self dual with respect to $\psi$, for which the schematic closure of $G_1$ in $\pmb{\GL}_{L_{(p)}}$ is a reductive group scheme, and for which we have $\pmb{GSp}(L_{(p)},\psi)(\dbZ_p)=K_p$. Then we have:

\medskip
{\bf (a)} The isogeny property, the endomorphism property, and the unramified CM lift property hold for all points $\Spec \dbF\to\scrN_1$. 

\smallskip
{\bf (b)} The $\dbF$-valued points of the basic locus $\grs_{1,0}$ of $\scrN_{1,k(v_1)}$ form one isogeny set. 

\smallskip
{\bf (c)} The functorial morphism $\scrN_1\to\scrM_{O_{(v_1)}}$ is a closed embedding.}

\medskip
\proof
Part (a) follows from the proof of Corollary 6.2. To prove (b) we have to show that for each point $y_1\in\grs_{1,0}(\dbF)$, we have an identity $\grs_{1,0}(\dbF)=o(y_1)$. To check this identity we can assume that $(G_0,X_0):=(G_1^{\ad},X_1^{\ad})$ is a simple, adjoint Shimura pair and that the injective map $(G_1,X_1,H_1,v_1)\hookrightarrow (\pmb{GSp}(W,\psi),S,K_p,p)$ is obtained as in Theorem 1.5 (cf. proof of Corollary 6.2 and Theorem 6.1). Therefore the identity $\grs_{1,0}(\dbF)=o(y_1)$ follows from Theorem 1.5 (e). Thus (b) holds.

To check part (c), we can assume that $G_1$ is not a torus (as this is trivial). It is known that $\scrN_1$ is the normalization of the schematic closure $\scrN_1^{\text{cl}}$ of $\Sh_{H_1}(G_1,X_1)$ in $\Sh_{K_p}(\pmb{GSp}(W,\psi),S)_{E(G_1,X_1)}$, cf. end of Subsection 1.1.2.  The inverse image via the map $\scrN_1(\dbF)\to\scrN_1^{\text{cl}}(\dbF)$ of the set $\text{Im}(\grs_{1,0}(\dbF)\to \scrN_1^{\text{cl}}(\dbF))$  is the set $\grs_{1,0}(\dbF)$ itself, cf. Proposition 3.2. Based on this and Lemma 3.6, the proof of (c) is the same as of Subsection 5.5, provided we know that the map $\grs_{1,0}(\dbF)\to\scrM(\dbF)=\scrM_{O_{(v_1)}}(\dbF)$ is injective. 

Let $o^{\text{big}}(y_1)$ be the isogeny set of $y_1$ but viewed as an $\dbF$-valued point of $\scrM$. As $\grs_{1,0}(\dbF)=o(y_1)$, to end the proof of (c) it suffices to show that the natural map of sets
$o(y)\to o^{\text{big}}(y)$
is injective. But this follows from Corollary 2.5.5. Thus (c) also holds. \endproof

\bigskip\noindent
{\bf 6.5. Corollary.} {\it Let $(G_1,X_1,H_1,v_1)\hookrightarrow (G_2,X_2,H_2,v_2)$ be an injective map of Shimura quadruples. We assume that each simple factor of $(G_2^{\ad},X_2^{\ad})$ is compact and of $A_n$, $B_n$, $C_n$, or $D_n^{\dbR}$ type (this implies that $(G_1,X_1)$ is also compact). We also assume that the monomorphism $H_1\hookrightarrow H_2$ is the $\dbZ_p$-valued points monomorphism associated to a monomorphism $G_{1,\dbZ_{(p)}}\hookrightarrow G_{2,\dbZ_{(p)}}$ between reductive group schemes over $\dbZ_{(p)}$ that extends the monomorphism $G_1\hookrightarrow G_2$. Let $\scrN_1\to\scrN_{2,O_{(v_1)}}$ be the functorial morphism, where $\scrN_2$ is the integral canonical model of $(G_2,X_2,H_2,v_2)$. Then the morphism $\scrN_1\to\scrN_{2,O_{(v_1)}}$ is a closed embedding.}

\medskip
\proof
We can assume that $G_2$ is not a torus (as this case is trivial). The morphism $\scrN_{1,E(G_1,X_1)}\to\scrN_{2,E(G_1,X_1)}$ is a closed embedding. Thus to prove that $\scrN_1\to\scrN_{2,O_{(v_1)}}$ is a closed embedding, we can consider pull-backs to $W(\dbF)$ and we can work in the pro-\'etale cover topology of  a fixed connected component of $\scrN_{2,W(\dbF)}$. Based on this and Theorem 6.1, as in the proof of Corollary 6.2 we argue that we can assume that we have an injective map $(G_2,X_2,H_2,v_2)\hookrightarrow (G_3,X_3,H_3,v_3)$ such that the following three properties hold:

\medskip
\item
{\bf (i)} We have an identity $G_2^{\der}=G_3^{\der}$ between simply connected semisimple groups.

\smallskip
\item
{\bf (ii)} The quadruple $(G_3,X_3,H_3,v_3)$ is a product $\prod_{i\in I}(G_{3,i},X_{3,i},H_{3,i},v_{3,i})$ of Shimura quadruples whose adjoints $(G_{3,i}^{\ad},X_{3,i}^{\ad},H_{3,i}^{\ad},v_{3,i}^{\ad})$ are simple and which admit injective maps $(G_{3,i},X_{3,i},H_{3,i},v_{3,i})\hookrightarrow (\pmb{GSp}(W_i,\psi_i),S_i,K_{i,p},p)$ into Shimura quadruples of Siegel modular varieties. Moreover, for each $i\in I$ there exists a $\dbZ_{(p)}$-lattice $L_{i,(p)}$ of $W_i$ which is self dual with respect to $\psi_i$ and with the properties that the schematic closure of $G_{3,i}$ in $\pmb{\GL}_{L_{i,(p)}}$ is a reductive group scheme and we have $K_{i,p}=\pmb{GSp}(L_{i,(p)},\psi_i)(\dbZ_p)$.

\smallskip
\item
{\bf (iii)} The torus $Z^0(G_2)$ is the extension of a torus which over $\dbR$ is compact by $\dbG_m$.

\medskip
Thus using a Hodge quasi product $\prod_{i\in I}^{\scrH} (G_{3,i},X_{3,i},H_{3,i},v_{3,i})$ as in [Va3, Subsect. 2.4] and an injective map of the form 
$$\prod_{i\in I}^{\scrH} (G_{3,i},X_{3,i},H_{3,i},v_{3,i})\hookrightarrow (\pmb{GSp}(\oplus_{i\in I} W_i,\oplus_{i\in I}\eps_i\psi_i),S_I,K_{I,p},p),$$
where we have $\eps_i\in\{-1,1\}$ for all $i\in I$, based on the properties (i) and (ii) we can assume that the injective map $(G_2,X_2,H_2,v_2)\hookrightarrow (G_3,X_3,H_3,v_3)$ factors through such a Hodge quasi product $\prod_{i\in I}^{\scrH} (G_{3,i},X_{3,i},H_{3,i},v_{3,i})$. We conclude that we can assume that we have an injective map $(G_2,X_2,H_2,v_2)\hookrightarrow (\pmb{GSp}(W,\psi),S,K_p,p)$ of Shimura quadruples such that there exists a $\dbZ_{(p)}$-lattice $L_{(p)}$ of $W$ which is self dual with respect to $\psi$, for which we have $K_p=\pmb{GSp}(L_{(p)},\psi)(\dbZ_p)$, and which has the property that the schematic closure of $G_2$ in $\pmb{\GL}_{L_{(p)}}$ is a reductive group scheme. One can make all the necessary arrangements so that the last reductive group scheme is $G_{2,\dbZ_{(p)}}$ itself (for instance, cf. [Va14, Part I, Lem. 4.2.1]). Therefore the schematic closure of $G_1$ in $\pmb{\GL}_{L_{(p)}}$ is the reductive group scheme $G_{1,\dbZ_{(p)}}$. Thus the corollary follows from Corollary 6.4 (c) applied to the functorial composite morphism $\scrN_1\to\scrN_{2,O_{(v_1)}}\to\scrM_{O_{(v_1)}}$.\endproof 

\medskip\noindent
{\bf 6.5.1. Remark.} Based on Corollary 6.5, one can easily get a functorial version of Corollary 6.2. 

\bigskip\noindent
{\bf 6.6. Main Corollary C.} {\it We assume that the isogeny property holds for all points $y:\Spec \dbF\to\scrN$ with $\scrN$ as in Subsubsection 1.1.2 associated to a Shimura pair $(G,X)$ with the property that the adjoint Shimura pair $(G^{\ad},X^{\ad})$ is a product of simple, adjoint Shimura pairs of $A_n$, $B_n$, $C_n$, or $D_n^{\dbR}$ type. Then Corollaries 6.2, 6.4 (a), and 6.4 (b) continue to hold without any compactness assumption.}

\medskip
\proof
The same proofs apply (the compactness assumption was only used to get a more elementary proof of Lemma 3.6 and indirectly of the isogeny property for all points $y:\Spec \dbF\to\scrN$ with $\scrN$ as in the corollary; to be compared with the footnote of Lemma 3.6).\endproof

\bigskip\smallskip
\noindent
{\boldsectionfont Appendix: Shimura $F$-crystals}
\bigskip 

The below notations are independent from the previous notations of the main text. Let $k$ be an algebraically closed field of positive characteristic $p>0$. Let $W(k)$ be the ring of Witt vectors with coefficients in $k$ and let $B(k)$ be the field of fractions of $W(k)$. Let $\sigma=\sigma_k$ be the Frobenius automorphism of $k$, $W(k)$, and $B(k)$. 

\bigskip\noindent
{\bf A1. Basic definitions.} A Shimura $F$-crystal over $k$ is a triple $C=(M,\phi,G)$, where $M$ is a free $W(k)$-module of finite rank, $\phi:M\to M$ is a $\sigma$-linear endomorphism such that $pM\subseteq \phi(M)$, and $G$ is a reductive, closed subgroup scheme of $\GL_M$ such that there exists a direct sum decomposition $M=F^1\oplus F^0$ for which the following two axioms hold (cf. [Va8,10]):

\medskip
{\bf (i)} we have identities $\phi(M+{1\over p}F^1)=M$ and $\phi(\Lie(G_{B(k)}))=\Lie(G_{B(k)})$, and 

\smallskip
{\bf (ii)} the cocharacter $\mu:\dbG_m\to \pmb{\GL}_M$ that acts trivially on $F^0$ and as the inverse of the identical character of $\dbG_m$ on $F^1$ (i.e., with weight $-1$ on $F^1$), factors through $G$. 

\medskip
The quadruple $C(F^1):=(M,F^1,\phi,G)$ is called a Shimura filtered $F$-crystal over $k$. Either $C(F^1)$ or $F^1$ is called a lift of $C$. Let $M^{\vee}:=\Hom(M,W(k)$ and let $\scrT(M):=\oplus_{s,t\in \db N} M^{\otimes s}\otimes_{W(k)} M^{\vee\otimes t}$ be the essential tensor algebra of $M\oplus M^{\vee}$.

We recall that $\phi$ act naturally on $\End(M)[{1\over p}]$ via the rule: $e\in\End(M)[{1\over p}]$ is mapped to $\phi(e):=\phi\circ e\circ \phi^{-1}$. Thus $\phi$ acts in the natural tensorial way on $\scrT(M)[{1\over p}]$ and therefore it makes sense to speak about the identity $\phi(\Lie(G_{B(k)}))=\Lie(G_{B(k)})$. 

Let $(t_{\alpha})_{\alpha\in\scrJ}$ be a family of tensors of $\scrT(M)$ fixed by $\phi$ and such that $G_{B(k)}$ is the subgroup of $\pmb{GL}_{M[{1\over p}]}$ that fixes $t_{\alpha}$ for all $\alpha\in\scrJ$ (cf. [Va10, proof of Claim 2.2.2] or [Va11, Lem. 2.5.3]). 

\bigskip\noindent
{\bf A2. Some group schemes associated to $C$.} Let $P^-$ be the unique parabolic subgroup scheme of $G$ such that the following two properties hold (cf. [Va10, Subsect. 2.3]): 

\medskip
{\bf (i)} the Lie algebra $\Lie(P^-)[{1\over p}]$ is left invariant by $\phi$;

\smallskip
{\bf  (ii)} all Newton polygon slopes of $(\Lie(P^-)[{1\over p}],\phi)$ are non-positive and all Newton polygon slopes of $(\Lie(G)[{1\over p}]/\Lie(P^-)[{1\over p}],\phi)$ are positive. 

\medskip
Let $U^-$ be the unipotent radical of $P^-$. Let $L_{B(k)}$ be the unique Levi subgroup of $P^-_{B(k)}$ such that its Lie algebra $\Lie(L_{B(k)})$ is normalized by $\phi$; all the Newton polygon slopes of $(\Lie(L_{B(k)}),\phi)$ are $0$. Let $L$ be the schematic closure of $P_{B(k)}$ in $G$; it is an integral, flat subgroup scheme of either $G$ or $P^-$. If $L$ is a Levi subgroup scheme of $P^{-1}$, then we call it the Levi subgroup scheme of $C$.

\bigskip\noindent
{\bf A3. Definitions.} We say that the Shimura $F$-crystal $C$ is:

\medskip
{\bf (a)} basic, if all the Newton polygon slopes of $(\Lie(G_{B(k)}),\phi)$ are $0$ (i.e., we have $L=P^-=G$);

\smallskip
{\bf (b)}  Levi, if $L$ is a Levi subgroup scheme of $P^{-}$;

\smallskip
{\bf (c)}  pivotal, if the following two properties hold (to be compared with [Va8, Def. 8.2 and Thm. 8.3]):

\medskip\noindent
{\bf (c.i)} for each element $g\in \Ker(G(W(k))\to G(k))$ there exists an element $g_1\in G(W(k))$ such that we have $g\phi=g_1\phi g_1^{-1}$;

\smallskip\noindent
{\bf (c.ii)} we have inclusions $p\Lie(G)\subseteq \{x\in\Lie(G)|\phi(x)=x\}\otimes_{\dbZ_p} W(k)\subseteq \Lie(G)$.

\medskip\noindent
{\bf A3.1. Remark.} If $C$ is pivotal, then there exists a maximal torus $T$ of $G$ such that there exists a cocharacter $\mu:\dbG_m\to G$ as in Subsection A1 that factors through $T$ and moreover we have $\phi(\Lie(T))=\Lie(T)$ (cf. [Va8, Basic Thm. C and Thm. 8.3]). 

\bigskip\noindent
{\bf A4. The formal deformation space of $C$.} We fix a direct sum decomposition $M=F^1\oplus F^0$ as in Subsection A1; thus $F^1$ is a lift of $C$. Let $P$ be the parabolic subgroup scheme of $G$ which is the normalizer of $F^1$ in $G$ (cf. [Va8, Subsubsect. 3.3.4]). Let $U$ be the unipotent radical of $P$; it is a commutative group scheme isomorphic to $\dbG_a^d$, where $b$ is the rank of $\Lie(G)\cap \Hom(F^1,F^0)$ (cf. [Va8, Subsubsects. 3.3.3 and 3.3.4]). Here we use the canonical direct sum decomposition $\End(M)=\End(F^1)\oplus \End(F^0)\oplus \Hom(F^1,F^0)\oplus \Hom(F^0,F^1)$ of $W(k)$-modules. Let $R$ be the completion of the local ring of $U$ at the identity element of $U_k$. We fix an identification $R=W(k)[[x_1,\ldots,x_b]]$.  Let $\Phi_R$ be the Frobenius endomorphism of $R$ which is compatible with $\sigma$ and maps each $x_i$ to $x_i^p$. Let $M_R:=M\otimes_{W(k)} R$ and $F^1_R=F^1\otimes_{W(k)} R$. Let $G_R=G\times_{\Spec W(k)} \Spec R$; it is a closed, reductive subgroup scheme of $\pmb{GL}_{M_R}$. Let $S:=R/pR=k[[x_1,\ldots,x_b]]$.

We consider the universal element
$$u_{\text{univ}}\in U(R)\leqslant G(R)$$ 
defined by the natural morphism $\Spec R\to U$ of $W(k)$-schemes. It is known that there exists a $p$-divisible group $D$ over $W(k)$ whose filtered Dieudonn\'e module is $(M,F^1,\phi)$. If $p>2$, then $D$ is unique and if $p=2$, then $D_k$ is unique. The existence of $D$ allows us to apply [Fa, Sect. 7, Thm. 10]: there exists a unique connection $\nabla:M_R\to M_R\otimes_{R} \oplus_{i=1}^b Rdx_i$ such that we have an identity
$$\nabla\circ u_{\text{univ}}(\phi\otimes\Phi_R)=(u_{\text{univ}}(\phi\otimes\Phi_R)\otimes d\Phi_R)\circ \nabla\leqno (17)$$ 
of maps from $M_R$ to $M_R\otimes_R \oplus_{i=1}^b Rdx_i$. The connection $\nabla$ is integrable and topologically nilpotent and $G$-invariant in the sense of [Va11, Def. 3.1.1 (d)] (thus the connection on $\scrT(M)\otimes_{W(k)} R$ induced naturally by $\nabla$ annihilates $t_{\alpha}$ for all $\alpha\in\scrJ$), cf. [Fa, Sect. 7] or [Va11, Thm. 3.2 and Cor. 3.3.2]. We emphasize that [Fa, Sect. 7] pertains to all primes including $2$ (odd primes are used in [Fa] only until [Fa, Sect. 7], cf. [Fa, Sect. 1]). Let $\nabla_0$ be the flat connection on $M_R$ that annihilates $M\otimes 1$. Due to Formula (17), it is easy to check that the connections $\nabla_0+u_{\text{univ}}^{-1}du_{\text{univ}}$ and $\nabla_y$ are congruent modulo $(p,x_1,\ldots,x_b)$. From this, the fact that $\Lie(U)$ is naturally a direct summand of $\Hom(F^1,M/F^1)$, and the definition of $u_{\text{univ}}\in U(R)$, we get that the Kodaira--Spencer map of $\nabla$ is injective modulo $(p,x_1,\ldots,x_b)$ and thus its image is a direct summand of $\Hom(F^1,M/F^1)\otimes_{W(k)} R$ of rank $b$.

Let $\scrD$ be the formal deformation space of $D$. It is isomorphic to 
$\Spf W(k)[[x_1,\ldots,x_m]]$ where $m$ is the product of the dimension and the codimension of $D_k$. We recall that the categories of $p$-divisible groups over $\Spf S$ and respectively over $\Spec S$, are canonically isomorphic (cf. [Me, Ch. II, Lem. 4.16]). From this and [Fa, Sect. 7, Thm. 10] we get the existence of a $p$-divisible group $E$ over $R$ whose filtered $F$-crystal over $R$ is $(M_R,F^1_R,u_{\text{univ}}(\phi\otimes\Phi_R),\nabla)$ and whose reduction modulo the ideal $(x_1,\ldots,x_b)$ of $R$ is $D$.  Let $E_S$ be the restriction of $E$ to $S$; it is a versal deformation of $D_k$ (cf. last paragraph) and therefore we can view naturally $\Spf S$ as a closed formal subscheme of $\scrD$. The below theorem shows that the closed formal subscheme $\Spf S$ of $\scrD$ is intrinsically associated to $C$. 

\medskip\noindent
 {\bf A4.1. Theorem.} {\it  The closed formal subscheme $\Spf S$ of $\scrD$ depends only on $C$ and not on the choices of the direct sum decomposition $M=F^1\oplus F^0$ as in Subsection A1 and of the identification $R=W(k)[[x_1,\ldots,x_b]]$.}
 
\medskip
\proof 
Let $M=\tilde F^1\oplus \tilde F^0$ be another direct sum decomposition  as in Subsection A1, let $\tilde F^1_R=\tilde F^1\otimes_{W(k)} R$, let $\tilde U$ be the closed, flat, smooth, commutative group subscheme of $G$ whose Lie algebra is $\Lie(G)\cap \Hom(\tilde F^1,\tilde F^0)$, let $\tilde R$ be the completion of the local ring of $\tilde U$ at the identity element of $\tilde U_k$, and let $\tilde u_{\text{univ}}\in\tilde U(\tilde R)$ be defined by the natural morphism $\Spec \tilde R\to\tilde U$ of $W(k)$-schemes. Let $\tilde S=\tilde R/p\tilde R$ and let $\Spf \tilde S$ be the closed formal subscheme of $\scrD$ obtained via the choice of an isomorphism $\tilde R=R=W(k)[[x_1,\ldots,x_b]]$ and thus via a $p$-divisible group $\tilde E$ over $R$ whose filtered $F$-crystal over $\tilde S$ is of the form $(M_R,\tilde F^1_R,\tilde u_{\text{univ}}(\phi\otimes\Phi_R),\tilde\nabla)$.  

As in [Va12, Prop. 6.4.6] we argue that there exists a morphism $z:\Spec W(k)\to \Spec \tilde R$ such that the Shimura filtered $F$-crystal associated to $z^*(\tilde E,G_R,(t_{\alpha})_{\alpha\in\scrJ})$ is isomorphic to $(M,F^1,\phi,G,(t_{\alpha})_{\alpha\in\scrJ})$. This implies that $(M_R,\tilde F^1_R,\tilde u_{\text{univ}}(\phi\otimes\Phi_R),\tilde\nabla,G_R,(t_{\alpha})_{\alpha\in\scrJ})$ is isomorphic to $(M_R,F^1_R,\tilde u(\phi\otimes\Phi_R),\tilde\nabla',G_R,(t_{\alpha})_{\alpha\in\scrJ})$ for a suitable element $\tilde u\in G(R)$ whose reduction modulo the ideal $(x_1,\ldots,x_b)$ of $R$ is the identity element of $G(W(k))$. Thus as in [Va12, Thm. 6.2.4] we argue that there exists a morphism $\psi:\Spec \tilde R\to\Spec R$ of $W(k)$-schemes such that $(M_R,\tilde F^1_R,\tilde u(\phi\otimes\Phi_R),\tilde\nabla',G,(t_{\alpha})_{\alpha\in\scrJ})$ is naturally identified with $\psi^*(M_R,F^1_R,u_{\text{univ}}(\phi\otimes\Phi_R),\nabla,G,(t_{\alpha})_{\alpha\in\scrJ})$ and which at the level of $W(k)$-algebras maps the ideal $(x_1,\ldots,x_b)$ of $R$ onto the ideal of $(x_1,\ldots,x_b)$ of $\tilde R=R$. By reasons of Kodaira--Spencer maps we easily get that $\psi$ is an isomorphism. This implies that $\Spf \tilde S=\Spf S$ as formal subschemes of $\scrD$.\endproof

\bigskip\noindent
{\bf A5. Stratifications of the formal deformation space $\Spf S$ of $C$.} Let $K$ be an algebraically closed field containing $k$. For a morphism $y:\Spec K\to \Spec S$, let $C_y$ be the pull-back of $(M_R,u_{\text{univ}}(\phi\otimes\Phi_R),\nabla,G_R)$ via $y$. As $\nabla$ respects the $G$-action, $C_y$ is a Shimura $F$-crystal over $K$ of the form $(M\otimes_{W(k)} W(K),g_y(\phi\otimes\sigma_K),G_{W(K)})$; here $G_{W(K)}:=G\times_{\Spec W(k)} \Spec W(K)$.

Let $\grc=\Spf S_{\grc}$, $\grr=\Spf S_{\grr}$, and $\grl=\Spf S_{\grl}$ be the reduced, closed formal subschemes of the formal deformation space $\Spf S$ of $C$ defined as follows:

\medskip
$\bullet$  A geometric point $y:\Spec K\to \Spec S$ factors through $\Spec S_{\grc}$ if and only if $C_y$ is isomorphic to $C\otimes K:=(M\otimes_{W(k)} W(K),\phi\otimes\sigma_K,G_{W(K)})$ under an isomorphism defined by an element $h_y\in G(W(K))$ (i.e., we have $g_y(\phi\otimes\sigma_K)=h_y(\phi\otimes\sigma_K)h_y^{-1}$).

$\bullet$  A geometric point $y:\Spec K\to \Spec S$ factors through $\Spec S_{\grr}$ if and only if $C_y$ is isogeneous to $C\otimes K$ under an isogeny defined by an element $l_y\in G(B(K))$ (i.e., we have $g_y(\phi\otimes\sigma_K)=r_y(\phi\otimes\sigma_K)r_y^{-1}$).

$\bullet$  A geometric point $y:\Spec K\to \Spec S$ factors through $\Spec S_{\grl}$ if and only if $C_y$ mod $p$ is isomorphic to $C\otimes K$ mod $p$ via an element of $l_y\in G(W(K))$ (i.e., we have $g_y(\phi\otimes\sigma_K)=c_y l_y(\phi\otimes\sigma_K)l_y^{-1}$ for some element $c_y\in \Ker(G(W(K))\to G(K))$).

\medskip
The existence of $\grr$ is a particular case of [RR, Thm. 3.6 (ii)]. The existence of $\grc$ and $\grl$ are a direct consequence of the global deformations of [Va1, Thm. 3.1] and of the existence of the level $m$ stratifications in [Va7, Subsects. 4.2 and 4.3]. 

We call $\grc$ as the formal Traverso stratum of $C$, $\grr$  as the formal rational stratum of $C$, and $\grl$ as the formal level $1$ stratum of $C$. It is known  that $\grc$ and $\grl$ are formally smooth (to be compared with [Va7, Cor. 4.3 (a)]). 

We say that $\grl$ has a non-trivial intersection with the basic locus, if there exists $c\in \Ker(G(W(k))\to G(k))$ such that $(M,c\phi,G)$ is basic. 

\bigskip\noindent
{\bf A6. Theorem.} {\it The following three properties hold:}

\medskip
{\bf (a)} {\it If $C$ is not basic, then $\dim(\grc)>0$.}

\smallskip
{\bf (b)} {\it If $C$ is not basic and $\grl$ has a non-trivial intersection with the basic locus, then we have $\dim(\grr\cap\grl)\ge 2$.}

\smallskip
{\bf (c)} {\it If $C$ is basic but not pivotal, then we have $\dim(\grr\cap\grl)\ge 1$.}

\medskip
\proof
We divide the proof into seven steps as follows.

\smallskip
{\bf Step 1.} From the non-positive standard form of $C$ (see [Va10, Subsect. 3.2]) we get the existence:

\medskip\noindent
{\bf (i)} of a maximal torus $T$ of $P^-$ such that there exists a cocharacter $\mu:\dbG_m\to G$ as in Subsection A1 that factors through $T$;

\smallskip\noindent
{\bf (ii)} of an element $g\in P^-(W(k))$ such that $\varphi:=g\phi$ normalizes the Lie algebra of $T$.

\medskip
Based on Theorem A4.1, we can assume that the cocharacter $\mu:\dbG_m\to G$ of the property (i) is the one that defines the direct sum decomposition $M=F^1\oplus F^0$ we used to construct $\Spf S$ in Subsection A4. 

Let $H$ be the unique Levi subgroup scheme of $P^-$ that contains $T$. We can write $g=u\ell$, where $u\in U^-(W(k))$ and $\ell\in H(W(k))$. We have a direct sum decomposition 
$\Lie(U^-)=\oplus_{i=-1}^1 \gru_i$ such that $\dbG_m$ acts on $\gru_i$ through $\mu$ and through inner conjugation via the weight $-i$. Let $c_{i}\in\dbN$ be the rank of $\gru_i$. We have an identity
$$\Lie(U^-)=\oplus_{i=-1}^1 \phi(p^{-i}\gru_i).\leqno (18)$$
\indent
{\bf Step 2.} Let $U^{-}_M$ be the analogue of $U^-$ but working with $\pmb{\GL}_M$ instead of with $G$. Thus $\Lie(U^{-}_M)$ is the largest direct summand of $\End(M)$ such that $\phi$ normalizes $\Lie(U^{-}_M)[{1\over p}]$ and all Newton polygon slopes of $(\Lie(U^{-}_M)[{1\over p}],\phi)$ are negative. Let $O_M$ be the largest $W(k)$-submodule of $\Lie(U^{-}_M)$ with the property that $(O_M,p\phi)$ is a Dieudonn\'e module. Let $O:=\Lie(U^-)\cap O_M$. The pair $(O,p\phi)$ is  the largest Dieudonn\'e module contained in $(\Lie(U^-),p\phi)$. The codimension of the $p$-divisible group over $k$ whose Dieudonn\'e module is $(O,p\phi)$, is $c_{-1}$ (cf. [Va9, Lem. 3.2 (a) and (b)]). It is known that $O_M$ is a (nilpotent) subalgebra of $\Lie(U^{-}_M)$, cf. [Va9, Subsect. 4.1]. Let $V_M$ be the affine, smooth group scheme over $W(k)$ defined by the rule: if $\star$ is a commutative $W(k)$-algebra, then $V_M(\star):=1_{M\otimes_{W(k)} \star}+O_M\otimes_{W(k)} \star$. We have a natural homomorphism $V_M\to U^{-}_M$ whose generic fibre is an isomorphism, to be viewed as a canonical identification. Let $V$ be the schematic closure of $U^-_{B(k)}$ in $V_M$. Let $W$ be a closed, smooth subscheme of $V$ which contains the identity section of $V$, which is isomorphic to $\Spec W(k)[x_1,\ldots,x_e]$ for some non-negative integer $e\in\{0,\ldots,c_{-1}\}$, and whose tangent space $\Lie(W)$ at the identity section is contained in $O$ and  we have a $k$-linear monomorphism
$\Lie(W)/p\Lie(W)\hookrightarrow \Lie(U^-)/[p\Lie(U^-)+\gru_0+\gru_1]$. 

\smallskip
{\bf Step 3.} We check that $\dim(\grc)\ge e$. 
Let $\Spec Q$ be the completion of the local ring of $W$ at the identity element of $W_k$. Let $w_{\text{univ}}:\Spec Q\to W$ be the universal (natural) morphism of $W(k)$-schemes. We fix an isomorphism  $Q=W(k)[[x_1,\ldots,x_e]]$ and let $\Phi_Q$ be the Frobenius endomorphism of $Q$ which is compatible with $\sigma$ and maps each $x_i$ to $x_i^p$. Using $W$, $w_{\text{univ}}$, $\Phi_Q$ (instead of $U$, $u_{\text{univ}}$, $\Phi_R$), as above we construct a deformation of $(M,\phi,(t_{\alpha})_{\alpha\in\scrJ},\psi_M)$ over $Q/pQ$ which is versal and which defines a formally smooth subscheme $\Spf Q/pQ$ of $\Spf S$ isomorphic to $\Spf k[[x_1,\ldots,x_e]]$. As in [Va9, Lem. 3.1 (b)] we argue that $\Spf Q/pQ$ is contained in $\grc$. This implies that $\dim(\grc)\ge e$. 

\smallskip
{\bf Step 4.} If $\Lie(U^-)^p=0$ inside $\End(M)$, then it is easy to see that the exponential map from $O$ to $V_M$ allows us to conclude that $V$ is smooth; this implies that we can choose $W$ such that we have $e=c_{-1}$. It seems to us that we can always choose $W$ such that we have $e=c_{-1}$. Here we will only check that in general we can choose $W$ such that we have the following inequality $e\ge\min\{2,c_{-1}\}$. 

To check this we can assume that $c_{-1}>0$. Let $Z(U^-)$ be the largest smooth, closed subgroup scheme of $U^-$ that centralizes $U^-$. It has positive relative dimension, it is isomorphic to a product of $\dbG_a$ group schemes, and the quotient group scheme $U^-/Z(U^-)$ exists and it is smooth. Let $\grz:=\Lie(Z(U^-))$. For $i\in\{-1,0,1\}$, let $e_{i}$ be the rank of $\grz_i:=\grz\cap \gru_i$. From the property (i) we get that $\mu$ normalizes $Z(U^-)$ and thus we have $\grz=\oplus_{i=-1}^1 \grz_i=\oplus_{i=-1}^1 \phi(p^{-i}\grz_i)$. If $e_{-1}\ge 2$, then we can choose $W$ to be a direct summand of $Z(U^-)$ which is of relative dimension $e_{-1}$ and which is contained in $V$. Thus to prove that $e\ge\min\{2,c_{-1}\}$, we can assume that $e_{-1}=1$ and that $c_{-1}\ge 2$.  

As $e_{-1}=1$ and as all the Newton polygon slopes of $(\grz[{1\over p}],\phi)$ are negative, from the identity $\grz=\oplus_{i=-1}^1 \phi(p^{-i}\grz_i)$ we get that $e_{1}=0$ (i.e., $\grz_1=0$) and that all Newton polygon slopes of $(\grz[{1\over p}],\phi)$ are equal to ${{-1}\over {1+e_{0}}}$. This implies that $(\grz,p\phi)$ is a $p$-divisible group of codimension $1$. Thus $\grz\subseteq O$, $Z(U^-)\leqslant V$, and moreover $(O/\grz,p\phi)$ is the largest Dieudonn\'e module contained in $(\Lie(U^-)/\grz,p\phi)$. Thus by working with the largest smooth, closed subgroup scheme $Z(U^-/Z(U^-))$ of $U^-/Z(U^-)$ that centralizes $U^-/Z(U^-)$, as in the previous paragraph we argue that $V/Z(U^-)$ contains a direct summand of $Z(U^-/Z(U^-))$ which is isomorphic to $\dbG_a$ and whose special fibre has a Lie algebra that has a non-trivial image in $\Lie(U^-)/[p\Lie(U^-)+\grz+\gru_0+\gru_1]$. Therefore we can choose $W$ to have a relative dimension $e$ at least equal to $e_{-1}+1\ge 2$. This ends the argument that we can choose $W$ such that $e\ge \min\{2,c_{-1}\}$. 

\smallskip
{\bf Step 5.} We choose $W$ such that $e\ge \min\{2,c_{-1}\}$, cf. Step 4. If $(M,\phi,G)$ is not basic, then $P^-\neq G$, $U^-$ is of positive dimension, $H\neq G$, the Newton polygon slopes of $(\Lie(U^-),p\phi)$ belong to the interval $[0,1)$, and we have $c_{-1}>0$ and $\dim(\grc)\ge e\ge 1$, cf. Step 3. Thus (a) holds.

\smallskip
{\bf Step 6.}  We check that (b) implies (c). Thus we will assume that $C$ is basic but not pivotal. As $C$ is not pivotal, we have $\dim(\grl)\ge 1$ (cf. [Va8, Cor. 11.1 (c) and Basic Thm. D]). If $\grl$ is contained in $\grr$, then $\dim(\grl\cap\grr)=\dim(\grl)\ge 1$. Thus we can assume that $\grl$ is not contained in $\grr$. We write $\grl\cap \grr=\Spf S_{\grl\cap\grr}$ with $S_{\grl\cap\grr}$ as a quotient ring of both $S_{\grl}$ and $S_{\grr}$. 

From (b) applied to suitable pull-backs $C_y$ obtained as in Subsection A5 using geometric points $\Spec K\to\Spec S_{\grl}\setminus \Spec S_{\grl\cap\grr}$ we get that each irreducible component of a stratum of the rational stratification of $\Spec S_{\grl}\setminus \Spec S_{\grl\cap\grr}$ associated naturally to $(M_R,u_{\text{univ}}(\phi\otimes\Phi_R),\nabla,G_R,(t_{\alpha})_{\alpha\in\scrJ})$ based on [RR, Thm. 3.6 (ii)] has dimension at least $2$. From this and the fact that the rational stratifications are pure (i.e., the open subscheme $\Spec S_{\grl}\setminus \Spec S_{\grl\cap\grr}$ of $\Spec S_{\grl}$ is an affine scheme, cf. [Va2, Main Thm. B]), we get that $\dim(\grr\cap\grl)\ge 1$. Thus (b) implies (c). 

\smallskip
{\bf Step 7.} Next we will prove (b) and (c) by induction on $\dim(G_k)$. Obviously (b) and (c) hold if $\dim(G_k)\le 2$. Let $f\in\dbN^*\setminus\{1,2\}$. We assume that (b) and (c) hold if $\dim(G_k)<f$ and we check that (b) and (c) hold even if $\dim(G_k)=f$. It suffices to check that (b) holds, cf. Step 6. Thus $(M,\phi,G)$ is not basic and therefore $\dim(H_k)<f=\dim(G_k)$. To prove (b) we can assume that $c_{-1}=1$ (as otherwise we have $\dim(\grr\cap\grl)\ge\dim(\grc)\ge e\ge\min\{2,c_{-1}\}\ge 2$). As $c_{-1}=1$ and as all the Newton polygon slopes of $(\Lie(U^-)[{1\over p}],\phi)$ are negative, from (18) we get that $c_1=0$ (i.e., $\gru_1=0$) and that all Newton polygon slopes of $(\Lie(U^-)[{1\over p}],\phi)$ are equal to ${{-1}\over {c_0+1}}$. This implies that $(\Lie(U^-),p\phi)$ is a $p$-divisible group of codimension $1$. Thus, by replacing $(M,\phi,G)$ with a Shimura $F$-crystal $(M,u_1\phi u_1^{-1},G)$ isomorphic to it for a suitable $u_1\in U^-(W(k))$, we can assume that $u=1_M$. Therefore the triple $(M,\phi,H)$ is a basic Shimura $F$-crystal. Let $\grc_0$, $\grr_0$, and $\grl_0$ be the analogues of $\grc$, $\grr$, and $\grl$ (respectively) but for $(M,\phi,H)$ instead of for $(M,\phi,G)$. Thus $\grc_0$, $\grr_0$, and $\grl_0$ are closed formal subschemes of $\grc$, $\grr$, and $\grl$ (respectively) and $\grc_0$ and $\grl_0$ are formally smooth. If $(M,\phi,H)$ is not pivotal, then by induction we have $\dim(\grr_0\cap \grl_0)\ge 1$. From this and the equality $c_{-1}=1$ we get that $\dim(\grr\cap \grl)\ge 1+c_{-1}=2$ (in other words, one can put together the two $1$ dimensional formal smooth subschemes of $\Spf S$ that come from $W$ and from $H$ to form a $2$ dimensional formal smooth subscheme of $\Spf S$ contained in $\grr\cap\grl$). Thus (b) holds if $(M,\phi,H)$ is not pivotal.  

We are left to show that the assumption that $(M,\phi,H)$ is pivotal leads to a contradiction. The assumption that $(M,\phi,H)$ is pivotal implies that we can assume that $g=\ell$ is the identity element (to be compared with Remark A3.1). Thus $\phi=\varphi$ normalizes $\Lie(T)$. This implies that there exist a $W(k)$-basis $\{w_1,\ldots,w_f\}$ for $\Lie(G)$ contained in $\Lie(U^-)\cup \Lie(H)\cup\Lie(U^-)^{\text{opp}}$ (with $\Lie(U^-)^{\text{opp}}$ as the opposite of $\Lie(U^-)$ with respect to $T$) and a permutation $\tau$ of the set $\{1,2,\ldots,f\}$, such that for all $i\in\{1,\ldots,f\}$ we have
$$\phi(w_i)=p^{\eps_i}w_{\tau(i)}$$ 
for some number $\eps_i\in\{-1,0,1\}$ and the $W(k)$-span of $w_i$ is normalized by $T$ (to be compared with [Va8, Subsubsect. 4.1.1]). As $(\Lie(U^-),p\phi)$ is a Dieudonn\'e module of codimension $1$, the pair $(\Lie(U^-)^{\text{opp}},\phi)$ is a Dieudonn\'e module of dimension $1$. From this and the fact that $(M,\phi,H)$ is pivotal, we get that for each cycle $(i_1\,\ldots\,i_s)$ of $\tau$ the sequence $(\eps_{i_1},\ldots,\eps_{i_s})$ does not contain any subsequence of the form $(-1,0,\ldots,0,-1)$ (to be compared with [Va8, Thm. 8.3] applied to $(M,\phi,H$)). This implies that $\dim(\grl)=0$ (cf. [Va8, Basic Thms. A and D (d) and Subsubsect. 4.1.1]) and therefore $(M,\phi,G)$ is pivotal (cf. [Va8, Cor. 11.1 (c) and Basic Thm. D]). This contradicts the fact that $(M,\phi,G)$ is not basic. Thus (b) holds. This ends the induction and thus also the proof of the theorem.

\bigskip\noindent
{\bf Acknowledgments.} We would like to thank Binghamton, Northwestern, and Bielefeld Universities, IAS--Princeton, and TIFR--Mumbai for good working conditions. We would also like to thank Thomas Zink for some comments and Laurent Fargues for pointing out a mistake in an earlier version of the Definition 1.4.4. 
\bigskip
\noindent
\references{37}
{\nspace{

\bigskip

\Ref[BHC]
A. Borel and Harish-Chandra,
\sl Arithmetic subgroups of algebraic groups,
\rm Ann. of Math. (2) {\bf 75} (1962), no. 3, 485--535.

\Ref[BLR]
S. Bosch, W. L\"utkebohmert, and M. Raynaud,
\sl N\'eron models,
\rm Ergebnisse der Mathematik und ihrer Grenzgebiete (3), Vol. {\bf 21}, Springer-Verlag, Berlin, 1990.

\Ref[BM]  P. Berthelot and W. Messing, 
\sl Th\'eorie de Dieudonn\'e cristalline. III, 
\rm The Grothendieck Festschrift, Vol. I, 173--247, Progr. Math., Vol. {\bf 86}, Birkh\"auser Boston, Boston, MA, 1990.

\Ref[BT]
F. Bruhat and J. Tits, 
\sl Groupes r\'eductifs sur un corps local. II. Sch\'emas en groupes. Existence d'une donn\'ee radicielle valu\'ee,
\rm Inst. Hautes \'Etudes Sci. Publ. Math., Vol. {\bf 60},  5--184, 1984.

\Ref[De1]
P. Deligne,
\sl Travaux de Shimura,
\rm S\'eminaire  Bourbaki, 23\`eme ann\'ee (1970/71), Exp. No. 389, Lecture Notes in Math., Vol. {\bf 244}, 123--165, Springer-Verlag, Berlin, 1971.

\Ref[De2]
P. Deligne,
\sl Vari\'et\'es de Shimura: interpr\'etation modulaire, et
techniques de construction de mod\`eles canoniques,
\rm Automorphic forms, representations and $L$-functions (Oregon State Univ., Corvallis, OR, 1977), Part 2,  247--289, Proc. Sympos. Pure Math., Vol. {\bf 33}, Amer. Math. Soc., Providence, RI, 1979.

\Ref[De3]
P. Deligne,
\sl Hodge cycles on abelian varieties,
\rm Hodge cycles, motives, and Shimura varieties, Lecture Notes in Math., Vol.  {\bf 900}, 9--100, Springer-Verlag, Berlin-New York, 1982.

\Ref[dJ]
J. de Jong, 
\sl Homomorphisms of Barsotti--Tate groups and crystals in positive characteristic, 
\rm Invent. Math. {\bf 134} (1998), no. 2,  301--333. Erratum: Invent. Math. {\bf 138}  (1999),  no. 1, 225.

\Ref[DG]
M. Demazure, A. Grothendieck, et al.,
\sl Sch\'emas en groupes. Vol. {\bf III},
\rm S\'eminaire de G\'eom\'etrie Alg\'ebrique du Bois Marie 1962/64 (SGA 3), Lecture Notes in Math., Vol. {\bf 153}, Springer-Verlag, Berlin-New York, 1970. 

\Ref[Fa]
G. Faltings,
\sl Integral crystalline cohomology over very ramified valuation rings, 
\rm J. Amer. Math. Soc. {\bf 12} (1999), no. 1, 117--144.

\Ref[Ha]
G. Harder,
\sl \"Uber die Galoiskohomologie halbeinfacher Matrizengruppen II,
\rm  Math. Z. {\bf 92} (1966), 396--415.

\Ref[Ja]
J. C. Jantzen,
\sl Representations of algebraic groups,
\rm  Pure and Applied Mathematics, Vol. {\bf 131}, Academic Press, Inc., Boston, MA, 1987.

\Ref[Ki] M. Kisin,
\sl Integral canonical models of Shimura varieties of abelian type, 
\rm J. Amer. Math. Soc. {\bf 23} (2010), no. 4,  967--1012.

\Ref[Ko] R. E. Kottwitz,
\sl Shimura varieties and $\lambda$-adic representations,
\rm Automorphic forms, Shimura varieties, and $L$-functions (Univ. of Michigan, Ann Arbor, MI, 1988), Vol. I, 161--209, Perspect. Math. {\bf 10}, Academic Press, Boston, MA, 1990. 

\Ref[LR]
R. Langlands and M. Rapoport,
\sl Shimuravariet\"aten und Gerben, 
\rm J. Reine Angew. Math. {\bf 378} (1987), 113--220.

\Ref[Mi1]
J. S. Milne,
\sl Points on Shimura varieties mod $p$,
\rm Automorphic forms, representations and $L$-functions (Oregon State Univ., Corvallis, OR, 1977), Part 2,  165--184, Proc. Sympos. Pure Math., Vol. {\bf 33}, Amer. Math. Soc., Providence, RI, 1979.

\Ref[Mi2]
J. S. Milne,
\sl The conjecture of Langlands and Rapoport for Siegel modular varieties,
\rm Bull. Amer. Math. Soc. (N.S.) {\bf 24} (1991),  no. 2, 335--341.

\Ref[Mi3]
J. S. Milne,
\sl The points on a Shimura variety modulo a prime of good
reduction,
\rm The Zeta functions of Picard modular surfaces, 153--255, Univ. Montr\'eal, Montreal, Quebec, 1992.

\Ref[Mi4]
J. S. Milne,
\sl Shimura varieties and motives,
\rm Motives (Seattle, WA, 1991), Part 2, 447--523, Proc. Sympos. Pure Math., Vol. {\bf 55}, Amer. Math. Soc., Providence, RI, 1994.

\Ref[Mi5]
J. S. Milne,
\sl Towards a proof of the conjecture of Langlands and Rapoport,

\rm http://arxiv.org/abs/0707.3177. 

\Ref[Mi6]
J. S. Milne,
\sl The Tate conjecture over finite fields,
\rm http://arxiv.org/abs/0709.3040. 

\Ref[Mi7]
J. S. Milne,
\sl Points on Shimura varieties over finite fields: the conjecture of Langlands and Rapoport,
\rm  40 pages manuscript, November 11, 2009, 

http://arxiv.org/abs/0707.3173.

\Ref[MFK]
D. Mumford, J. Fogarty, and F. Kirwan, 
\sl Geometric invariant theory. Third edition, 
\rm Ergebnisse der Math. und ihrer Grenzgebiete (2), Vol. {\bf 34}, Springer-Verlag, Berlin, 1994.

\Ref[Oo1]
F. Oort, 
\sl Newton polygons and formal groups: conjectures by Manin and Grothendieck, 
\rm Ann. of Math. (2) {\bf 152} (2000), no. 1, 183--206.

\Ref[Oo2]
F. Oort, 
\sl Newton polygon strata in the moduli of abelian varieties, 
\rm Moduli of abelian varieties (Texel Island, 1999), 417--440, Progr. Math., {\bf 195}, Birkh\"auser, Basel, 2001.

\Ref[Pf]
M. Pfau,
\sl The conjecture of Langlands and Rapoport for certain Shimura varieties of non-rational weight,
\rm J. Reine Angew. Math. {\bf 471} (1996), 165--199. 

\Ref[Pi]
R. Pink,
\sl $l$-adic algebraic monodromy groups, cocharacters, and the Mumford--Tate conjecture,
\rm J. Reine Angew. Math. {\bf 495} (1998), 187--237.

\Ref[Re]
H. Reimann,
\sl The semi-simple zeta function of quaternionic Shimura varieties,
\rm Lecture Notes in Math., Vol. {\bf 1657}, Springer-Verlag, Berlin, 1997.

\Ref[RR] 
M. Rapoport and M. Richartz, 
\sl On the classification and specialization of $F$-isocrystals with additional structure, 
\rm Compositio Math. {\bf 103} (1996), no. 2, 153--181.

\Ref[Va1]
A. Vasiu, 
\sl Integral canonical models of Shimura varieties of preabelian
 type, 
\rm Asian J. Math., Vol. {\bf 3} (1999), no. 2, 401--518.

\Ref[Va2]
A. Vasiu,
\sl Crystalline boundedness principle,
\rm Ann. Sci. \'Ecole Norm. Sup. {\bf 39} (2006), no. 2, 245--300.

\Ref[Va3]
A. Vasiu, 
\sl Some cases of the Mumford-Tate conjecture and Shimura varieties,     
\rm Indiana Univ. Math. J. {\bf 57} (2008), no. 1, 1--76.

\Ref[Va4]
A. Vasiu, 
\sl Projective integral models of Shimura varieties of Hodge type with compact factors,
\rm J. Reine Angew. Math. {\bf  618} (2008), 51--75. 

\Ref[Va5]
A. Vasiu, 
\sl Integral canonical models of unitary Shimura varieties, 
\rm Asian J. Math., Vol. {\bf 12} (2008), no. 2, 151--176.

\Ref[Va6]
A. Vasiu, 
\sl Geometry of Shimura varieties of Hodge type over finite fields,
\rm Higher-dimensional geometry over finite fields, 197--243, NATO Sci. Peace Secur. Ser. D Inf. Commun. Secur., {\bf 16}, IOS, Amsterdam, 2008.

\Ref[Va7]
A. Vasiu,
\sl Level $m$ stratifications of versal deformations of $p$-divisible groups,
\rm J. Alg. Geom. {\bf 17} (2008), no. 4, 599--641.

\Ref[Va8]
A. Vasiu,
\sl Mod $p$ classification of Shimura $F$-crystals,
\rm Math. Nachr. {\bf 283} (2010), no. 8, 1068--1113.

\Ref[Va9]
A. Vasiu,
\sl Deformation subspaces of $p$-divisible groups as formal Lie group structures associated to $p$-divisible groups,
\rm J. Alg. Geom. {\bf 20} (2011), no. 1, 1--45.

\Ref[Va10]
A. Vasiu,
\sl Manin problems for Shimura varieties of Hodge type,
\rm J. Ramanujan Math. Soc. {\bf 26} (2011), no. 1, 31--84.

\Ref[Va11]
A. Vasiu,
\sl A motivic conjecture of Milne,
\rm  67 pages to appear in J. Reine Angew. Math., http://www.degruyter.com/view/j/\break crelle.ahead-of-print/crelle-2012-0009/crelle-2012-0009.xml?format=INT.

\Ref[Va12]
A. Vasiu,
\sl Generalized Serre--Tate ordinary theory,
\rm 196 pages (including contents and index), to be published in 2013 by International Press, Inc., available at http://www.math.binghamton.edu/adrian/\#reductive.

\Ref[Va13]
A. Vasiu,
\sl CM lifts for isogeny classes of Shimura $F$-crystals over finite fields,
\rm 62 pages manuscript, June 18, 2012, http://www.math.binghamton.edu/adrian/.

\Ref[Va14]
A. Vasiu
\sl Good reductions of Shimura varieties of Hodge type in arbitrary unramified mixed characteristic, Parts I and II,
\rm 53 and 29 pages manuscripts, dated July 24, 2012, http://xxx.arxiv.org/abs/0707.1668 and http://arxiv.org/abs/0712.1572.

\Ref[VZ]
A. Vasiu and T. Zink, 
\sl Purity results for $p$-divisible groups and abelian schemes over regular bases of mixed characteristic, 
\rm Documenta Math. {\bf 15} (2010), 571--599. 

\Ref[Wi]
J.-P. Wintenberger,
\sl Un scindage de la filtration de Hodge pour certaines
variet\'es algebriques sur les corps locaux, 
\rm Ann. of Math. (2) {\bf 119} (1984), no. 3, 511--548.

\Ref[Zi]
T. Zink,
\sl Isogenieklassen von Punkten von Shimuramannigfaltigkeiten mit Werten in einem endlichen K\"orper,
\rm Math. Nachr. {\bf 112} (1983), 103--124.

}}

\bigskip
\hbox{Adrian Vasiu}
\hbox{Department of Mathematical Sciences, Binghamton University}
\hbox{P. O. Box 6000, Binghamton, New York 13902-6000, U.S.A.}
\hbox{e-mail: adrian\@math.binghamton.edu,\;\;fax: 1-607-777-2450,\;\;phone 1-607-777-6036}

\enddocument